\documentclass[smallextended,referee,envcountsect,]{svjour3}
\smartqed
\usepackage{comment}
\usepackage{graphicx}
\usepackage{amsfonts}
\usepackage{amssymb}
\usepackage{amsmath,amsfonts,amssymb}
\usepackage{mathtools}
\usepackage{algorithmic}
\usepackage{algorithm}
\usepackage{array}
\usepackage{multirow}
\usepackage[table]{xcolor}
\usepackage{caption}
\usepackage{subcaption}
\usepackage{multirow}
\usepackage{array}
\usepackage{longtable}
\usepackage{hyperref}
\usepackage{rotating,tabularx}
\usepackage{geometry}
\def\EMAIL#1{\href{mailto:#1}{#1}}

\begin{document}

\title{Nonconvex quasi-variational inequalities: stability analysis and application\\ to numerical optimization}

\titlerunning{Nonconvex quasi-variational inequalities}        

\author{Joydeep Dutta \and Lahoussine Lafhim \and Alain Zemkoho \and Shenglong Zhou}

\institute{Joydeep Dutta \at
              Department of Economic Sciences, Indian Institute of Technology Kanpur, India,           \EMAIL{jdutta@iitk.ac.in}   
           \and
          Lahoussine Lafhim \at
           Department of Mathematics, Sidi Mohammed Ben Abdellah University, Morocco,
           \EMAIL{lahoussine.lafhim@usmba.ac.ma}
           \and 
         Alain Zemkoho \at
             School of Mathematical Sciences, University of Southampton, United Kingdom,
              \EMAIL{a.b.zemkoho@soton.ac.uk}
              \and 
         Shenglong Zhou \at
           School of Mathematics and Statistics, Beijing Jiaotong University, Beijing, China,
           \EMAIL{shlzhou@bjtu.edu.cn}
}

\date{Received: date / Accepted: date}

\maketitle

\begin{abstract}
We consider a parametric quasi-variational inequality (QVI) without any convexity assumption. Using the concept of \emph{optimal value function}, we transform the problem into that of solving a nonsmooth system of inequalities.  Based on this reformulation, new coderivative estimates as well as robust stability conditions for the optimal solution map of this QVI are developed. Also, for an optimization problem with QVI constraint, necessary optimality conditions are constructed and subsequently, a tailored semismooth Newton-type method is designed, implemented, and tested on a wide range of optimization examples from the literature. In addition to the fact that our approach does not require convexity, its coderivative and stability analysis do not involve second order derivatives, and subsequently, the proposed Newton scheme does not need third order derivatives, as it is the case for some previous works in the literature. 

\end{abstract}


\section{Introduction}
Quasi-variational inequalities have been widely studied in the literature,
 considering the large number of applications and the mathematical challenges involved in solving them;
 see, e.g., \cite{ChanPang1982,PangFukushima2005,OutrataZowe1995} and references therein, for some numerical
 methods to solve different classes of the problem, as well as a number of applications.  
 
In this paper, we consider a parametric quasi-variational inequality (QVI) with the first primary goal to construct estimates of the coderivative of the solution set-valued mapping (in the sense of Mordukhovich \cite{MordukhovichBook2006})  and sufficient conditions ensuring that it is Lipschitz-like (in the sense of Aubin \cite{Aubin1984}).
Some of these questions have been addressed in the literature (see, e.g., \cite{MOC07,GfrererYe2020NewNecc,GfrererYe2020NewConst}). However, as it is common in the literature, the constraint set of the QVI problem is assumed in the existing works to be convex. Hence, as usual, the analysis typically relies on a generalized equation reformulation based on the normal cone to this constraint set, in the sense of convex analysis.   In the absence of this convexity assumption, this approach is not applicable. Secondly, the stability analysis of this solution set-valued mapping (computation of its coderivative and construction of conditions ensuring that it is Lipschitz-like) conducted through \eqref{S-KKT} would require second order information, which is not available for various applications (see, e.g., \cite{HobbsPang2007}). 

The core of the analysis in this paper is based on a \emph{value function reformulation} of the  QVI problem.
Unlike in the existing literature, our transformation does not require any convexity assumption. We then use this value function reformulation we develop completely new results for solution set-valued mapping of our parametric QVI without any convexity assumption. Additionally, no second order derivatives are involved in our analysis, as it would be the case when using the aforementioned generalized equation reformulation of a QVI. 

As second main aim in this paper, we consider an optimization problem, partly constrained by a parametric QVI; i.e., an optimization problem with a quasi-variational inequality constraint (OPQVI). For this problem, we also introduce a \textit{value function reformulation}, which does not seem to have been considered in the existing literature. We use this reformulation to develop necessary optimality conditions for our OPQVI problem under suitable \emph{calmness} conditions. Subsequently, we develop a version of the semismooth Newton method tailored to these necessary optimality conditions and establish its convergence. Extensive experiments on over 124 examples from the literature \cite{MOC07,BOLIB2017} are then conducted to demonstrate the efficiency potential of our method.

Special classes of the OPQVI problem have been studied in many papers in the literature; see, e.g., \cite{AdamHenrionOutrata2018,MOC07,HenrionSurowiecCalmness2011,HenrionOutrataSurowiec2012,GfrererYe2020NewNecc,GfrererYe2020NewConst}. However, the main focus has usually been on the derivation of necessary optimality conditions based on the generalized (GE) reformulation, with a special attention on the construction of suitable qualification conditions. To the best of our knowledge, not much is available in the literature in terms of solution algorithms for the OPQVI problem studied in this paper. Moreover, we are not aware of any work where the problem has been studied in the absence of convexity. Also, an attempt to develop a Newton-type method for the problem from the perspective of the GE reformulation would require third order derivatives for constraint functions involved in the QVI. This is another reason why it is attractive to develop such methods based on the value function reformulation, as done in this paper. Additionally, the extent of the numerical experiments conducted in this paper could inspire more work on the development of numerical methods for the OPQVI problem.

For the { {remainder}} of the paper, in the next section, we present the QVI and OPQVI problems studied in this paper. Subsequently, in Section \ref{Basic concepts and background material}, we first provide some preliminary tools from the variational analysis that will be needed in the subsequent sections. In particular, the focus in Section \ref{Basic concepts and background material} is on generalized differentiation tools for nonsmooth functions and set-valued mappings. In Section \ref{Coderivative and robust stability of solution maps}, we develop stability results for the solution set-valued mapping for our parametric QVI;  the presence of the value function constraint in this solution set-valued mapping motivates the introduction and study of a version of the \emph{uniform weak sharp minimum} concept tailored to our QVI and provide sufficient conditions ensuring that it holds. Finally, in Section \ref{Optimization problems with a QVI constraint}, a \emph{partial calmness} concept tailored to the value function reformulation of the OPQVI problem is introduced to build necessary optimality conditions. Subsequently, suitable second order conditions are then introduced to establish the convergence of a semismooth Newton scheme introduced to solve the aforementioned necessary conditions.

\section{Description of the problems studied in this paper}
Our basic problem of interest in this paper is the quasi-variational inequality (QVI) problem to find $y\in K(x, y)$, for a given parameter $x\in \mathbb{R}^n$, such that we have 
\begin{equation}\tag{QVI}\label{VI}
\left\langle f_0(x,y), \;\, \varsigma - y \right\rangle \geq 0 \;\;\;
 \forall \varsigma\in  K(x, y):= \left\{\varsigma \in \mathbb{R}^m\left|\; g_0(x, y, \varsigma)\leq 0\right.\right\},
\end{equation}
where $f_0 : \mathbb{R}^n\times \mathbb{R}^m \rightarrow \mathbb{R}^m$ is a continuously differentiable function
and $K :\mathbb{R}^n\times \mathbb{R}^m \rightrightarrows \mathbb{R}^m$ denotes the feasibility
set-valued map defined by a function $g_0\, : \mathbb{R}^n\times \mathbb{R}^m\times \mathbb{R}^m \rightarrow \mathbb{R}^q$,
 which is also assumed to be continuously differentiable.
 This paper has two main goals, with  the first
 one being to study the stability of the solution set-valued mapping $S : \mathbb{R}^n \rightrightarrows \mathbb{R}^m$
 \begin{equation}\label{S-Soultion-Map}
\begin{array}{c}
  S(x):=\left\{y\in K(x,y)\left|\; \left\langle f_0(x,y), \;\, \varsigma - y \right\rangle \geq 0  \;\;\; \forall \varsigma\in K(x, y)\right.\right\}
\end{array}
\end{equation}
associated to the \eqref{VI}. More precisely, we aim to construct estimates of the coderivative of the set-valued mapping $S$ \eqref{S-Soultion-Map}), in the sense
of Mordukhovich \cite{MordukhovichBook2006},   and sufficient conditions ensuring that it is Lipschitz-like, in the sense of Aubin \cite{Aubin1984}. 
If we assume that the set $K(x, y)$ is convex for all $(x,y)$, $S$ \eqref{S-Soultion-Map} can be rewritten as  
\begin{equation}\label{S-KKT}
S(x)= \left\{y\in \mathbb{R}^m\left|\, 0\in f_0(x,y) + N_{K(x, y)}(y)\right.\right\},
\end{equation}
where $N_{K(x, y)}(y)$ represents the normal cone, in the sense of convex analysis, to the set $K(x, y)$ at the point $y$;
cf. next section for the definition of this concept. 

The expression of $S$ in \eqref{S-KKT}, known as the generalized equation (GE) reformulation of \eqref{S-Soultion-Map}, has been widely used in the literature for stability analysis purposes or to develop numerical methods for variational inequalities when stability analysis is not necessarily under consideration; see, e.g., \cite{MOC07,GfrererYe2020NewNecc,GfrererYe2020NewConst} and references therein. Obviously, the GE reformulation \eqref{S-KKT}  is only possible if the set-valued mapping $K$ \eqref{VI} is convex-valued. In the absence of this convexity assumption, this approach is not applicable. Secondly, as mentioned in the introduction, the stability analysis of $S$ (computation of the coderivative and construction of conditions ensuring that it is Lipschitz-like) conducted through \eqref{S-KKT} would require second order information (which is not available for various applications; see, e.g., \cite{HobbsPang2007}) for the function $g_0$, which describes $K$ \eqref{VI}.  
Our analysis in this paper is based on the fact that the set-valued mapping $S$ \eqref{S-Soultion-Map} can easily be rewritten, without any assumption, as
\begin{equation}\label{VI-1}
 S(x):=\left\{y\in K(x,y)\left|\;\, y^\top f_0(x,y)-\varphi(x,y)\leq 0\right.\right\}
\end{equation}
for all $x\in \mathbb{R}^n$. Here, $\varphi$ denotes the optimal value function
\begin{equation}\label{varphi}
   \varphi(x,y):= \underset{\varsigma}\min \{\varsigma^\top f_0(x, y)|\;\; \varsigma\in K(x, y)\}.
\end{equation}
The core of our analysis is based on this \emph{value function reformulation} of the  problem  \eqref{VI}.
Clearly, unlike \eqref{S-KKT}, transformation \eqref{VI-1} does not require any convexity assumption. Hence, based on \eqref{VI-1}--\eqref{varphi}, we develop completely new results for $S$ \eqref{S-Soultion-Map} without any convexity assumption. Additionally, no second order derivatives are involved in our analysis, as it would be the case when using the GE reformulation in \eqref{S-KKT}.

After conducted the stability analysis of solution set-valued mapping $S$ via the value function reformulation \eqref{VI-1} (in Section \ref{Basic concepts and background material}), as second main objective of this paper, we will focus our attention on the following optimization problem with a quasi-variational inequality constraint (OPQVI):
\begin{equation}\tag{OPQVI}\label{OPQVIC}
\underset{x, y}{\min}~F(x,y) \;\; \mbox{ s.t. } \;\; G(x,y)\leq 0, \;\; y\in S(x).
\end{equation}
{ {Here, the functions $F\,:\mathbb{R}^n\times \mathbb{R}^m \rightarrow \mathbb{R}$ and $G\,:\mathbb{R}^n\times \mathbb{R}^m \rightarrow \mathbb{R}^p$ are also assumed to be continuously differentiable}},
and the set-valued mapping $S$ is defined as in \eqref{S-Soultion-Map}. Problem \eqref{OPQVIC} has a two-level optimization structure, with $F$ and $G$ representing the upper-level objective and constraint functions, respectively.
In the same vein, $S$ corresponds to the optimal solution set-valued mapping of the lower-level problem, which,
unlike in standard two-level/bilevel optimization (see, e.g., \cite{DempeFoundations2002,DempeZemkohoBook}),
is a quasi-variational inequality defined as in \eqref{VI}. In Section \ref{Optimization problems with a QVI constraint}, we will also use the value function reformulation \eqref{VI-1} to transform problem \eqref{OPQVIC} into a single-level problem, and then while using relevant results from the analysis of $S$ in Section \ref{Basic concepts and background material}, we will build necessary optimality conditions for the problem. These optimality conditions  will then be deployed for the development of a semismooth Newton-type method for which the convergence will be studied and numerical experiments conducted to show its effectiveness. 

\section{Basic concepts and background material}\label{Basic concepts and background material}
We start this section with some basic notation to be mostly used in the following main sections of the paper. First, for a vector $x\in\mathbb{R}^{n}$ and a scalar $\epsilon  >0$, we denote by
\[
\mathbb{U}_{\epsilon}(x):= \left\{y\in\mathbb{R}^{n}|~\|y-x\|_{\infty} < \epsilon \right\} \; \mbox{ and } \; \mathbb{B}_{\epsilon}(x):= \{y\in\mathbb{R}^{n}|~\|y-x\|_{\infty} \leq \epsilon \}
\]
the open and closed $\epsilon$-balls around $x$,  respectively. { {For a matrix $M\in \mathbb{R}^{m\times n}$ and an index set $I \subset \{1,\cdots ,m\}$, $M_{I}\in \mathbb{R}^{\mid I\mid \times n}$ denotes the matrix obtained from $M$ by deleting all rows with indices that do not belong to $I$.}}

For the points $(\bar x, \bar y)$, $(\bar x, \bar y, \bar y)$, and $(\bar x, \bar y, \bar\varsigma)$, the following notation will be used to collect { {active index sets}} of constraints for the feasible sets associated to problems \eqref{VI} and \eqref{OPQVIC}, which are defined by the functions $(x, y) \mapsto G(x, y)$, $(x, y) \mapsto g_0(x, y, y)$, and $(x, y, \varsigma) \mapsto g_0(x, y, \varsigma)$, respectively:
\begin{equation}\label{I1}
\begin{array}{lclll}
     I^1&:=& I^G(\bar x, \bar y)                   &:=&  \left\{i\in \{1, \ldots, p\}\left|~G_i(\bar x, \bar y)=0\right.\right\},\\[1ex]
    I^2 &:=& I^{g_0}(\bar x, \bar y, \bar y)       &:=& \left\{j\in \{1, \ldots, q\}|~g_{0_j}(\bar x, \bar y, \bar y)=0\right\},\\[1ex]
    I^3 &:=& I^{g_0}(\bar x, \bar y, \bar\varsigma)&:=& \left\{j\in \{1, \ldots, q\}|~g_{0_j}(\bar x, \bar y, \bar\varsigma)=0\right\}.
\end{array}
\end{equation}
If we associate to a point $(\bar x, \bar y)$, feasible for the constraint defined by $G$ \eqref{OPQVIC}, a Lagrange multiplier $\bar u$, we can proceed with the following standard partition of the indices:
\begin{equation}\label{multiplier sets}
\begin{array}{lllll}
 \eta^1  &:=& \eta^G(\bar x, \bar y, \bar u) &:=& \{i\in \{1, \ldots, p\}|\; \bar u_i =0, \, G_i(\bar x, \bar y)<0\},\\[1ex]
\theta^1  &:=&  \theta^G(\bar x, \bar y, \bar u)&:=& \{i\in \{1, \ldots, p\}|\; \bar u_i =0, \, G_i(\bar x, \bar y)=0\},\\[1ex]
\nu^1  &:=& \nu^G(\bar x, \bar y, \bar u)&:=& \{i\in \{1, \ldots, p\}|\; \bar u_i >0, \, G_i(\bar x, \bar y)=0\}.
\end{array}
\end{equation}
Similarly to \eqref{multiplier sets}, for a feasible point $(\bar x, \bar y)$ (resp. $(\bar x, \bar y, \bar\varsigma)$) to the constraint defined by the function $(x, y) \mapsto g_0(x, y, y)$ (resp. $(x, y, \varsigma) \mapsto g_0(x, y, \varsigma)$) and a corresponding Lagrange multiplier vector $\bar v$ (resp. $\bar w$),  we can analogously define the following partitions:
\begin{equation}\label{nu2nu3}
    \begin{array}{lllllllll}
\eta^2&:=& \eta^{g_0}(\bar x, \bar y, \bar y, \bar v),                   & \theta^2&:=&\theta^{g_0}(\bar x, \bar y, \bar y, \bar v),                   & \nu^2&:=&\nu^{g_0}(\bar x, \bar y, \bar y, \bar v),\\[1ex]
\eta^3&:=& \eta^{g_0}(\bar x, \bar y, \bar\varsigma, \bar w),& \theta^3&:=&\theta^{g_0}(\bar x, \bar y, \bar\varsigma, \bar w),& \nu^3&:=&\nu^{g_0}(\bar x, \bar y, \bar\varsigma, \bar w).
\end{array}
\end{equation}

{ {Although symbols that we use here are very standard, it might be useful to make it clear how we use the symbol $\nabla$ interchangeably between the gradient of a scalar function and the Jacobian of a vector-valued function. For a continuously differentiable function $\psi : \mathbb{R}^n\times \mathbb{R}^m \times \mathbb{R}^p \rightarrow \mathbb{R}^q$, in case there is ambiguity in the use of the usual partial derivatives, we will utilize $\nabla_1 \psi(x, y, \varsigma)$, (resp. $\nabla_2 \psi(x, y, \varsigma)$, $\nabla_3 \psi(x, y, \varsigma)$) to represent the Jacobian of the function $x \mapsto \psi(x, y, \varsigma)$ (resp. $y \mapsto \psi(x, y, \varsigma)$, $\varsigma \mapsto \psi(x, y, \varsigma)$) when the vector $(y, \varsigma)$ (resp. $(x, \varsigma)$, $(y, \varsigma)$) is fixed. A similar notation will be used for the gradient of $\psi$ if it is a real-valued function. Common scenarios for $\psi$ in this paper will be $\psi(x, y, \varsigma):= g_0(x, y, \varsigma)$ and $\psi(x, y, \varsigma):= \varsigma^\top f_0(x, y)$ when studying our \eqref{VI}.
}}

\subsection{Nonsmooth functions and normal cones}\label{Nonsmooth functions and normal cones}
Given that the optimal value function $\varphi$ \eqref{varphi} is nonsmooth in general, we need generalized concepts of
 differentiability to deal with it. We start the discussion on this by recalling the \emph{generalized directional derivative} in the sense of
 Clarke \cite{ClarkeBook1983}, which, for a function $\psi : \mathbb{R}^n \rightarrow \mathbb{R}$, is defined at a point $\bar{x}\in \mathbb{R}^n$ in direction $d\in \mathbb{R}^n$ by
 \begin{equation*}\label{Clarke Directional Derivative}
    \psi^0(\bar x; d):=\underset{t\downarrow 0}{\underset{x \rightarrow \bar x}\limsup} \frac{1}{t}\left[\psi(x + td)-\psi(x)\right].
\end{equation*}
 This quantity exists if $\psi$ is a Lipschitz continuous function around $\bar x$ \cite[Proposition 2.1.1]{ClarkeBook1983}.
  Utilizing this notion, the \textit{Clarke subdifferential} of $\psi$ at the point $\bar{x}$ can be defined by 
\begin{equation}\label{Clarke Subdifferential}
    \bar\partial \psi(\bar x):= \left\{\zeta\in \mathbb{R}^n|\; \psi^o(\bar x; d)\geq \langle\zeta, d\rangle, \; \forall d\in \mathbb{R}^n  \right\}.
\end{equation}
Note that {{$\bar\partial \psi(\bar x) = \left\{\nabla \psi(\bar x)\right\}$ }} if $\psi$ is continuously differentiable at $\bar x$, and for a convex function,
this concept coincides with the subdifferential in the sense of convex analysis.
Furthermore, this concept can be extended to a vector-valued function $\psi : \mathbb{R}^n \rightarrow \mathbb{R}^m$ that
is Lipschitz continuous around a point $\bar x$, as such a function is differentiable almost everywhere around this point. Hence,
 the \emph{Clarke generalized Jacobian} of the function $\psi$ at the point $\bar{x}$ can be written as 
\begin{equation}\label{Clarke Subdifferential-Derivative}
    \bar\partial \psi(\bar x):= \mbox{co}\left\{\lim\nabla\psi(x^n)\left|~x^n \rightarrow \bar x, \;\, x^n\in D_\psi \right.\right\},
\end{equation}
where ``co'' stands for the convex hull and  $D_\psi$ represents the set of points where $\psi : \mathbb{R}^n \rightarrow \mathbb{R}^m$ is
differentiable. { {Note that the quantity defined in  \eqref{Clarke Subdifferential-Derivative} coincides with the one in \eqref{Clarke Subdifferential}
 for real-valued functions}}. Following  the expression in \eqref{Clarke Subdifferential-Derivative},
the $B$-\emph{subdifferential} (see, e.g., \cite{Qi1993convergence}) can be defined by
\begin{equation}\label{B-Subdifferential}
    \partial_B \psi(\bar x):= \left\{\lim~\nabla \psi(x^n)\left|~x^n \rightarrow \bar x, \;\, x^n\in D_\psi\right.\right\}.
\end{equation}

The concept of semismoothness will play an important role in the convergence analysis of the Newton method
that will be introduced in Section \ref{Optimization problems with a QVI constraint} of this paper. Let a
function  $\psi : \mathbb{R}^{n} \rightarrow \mathbb{R}^{m}$ be Lipschitz continuous around the point $\bar x$. Then 
$\psi$ will be said to be \emph{semismooth} at $\bar x$ if the limit
\[
\lim \left\{V^\prime d^\prime\left|~V^\prime\in \bar\partial \psi (\bar x + td^\prime), \; d^\prime \rightarrow d, \; t\downarrow 0\right.\right\}
\]
exists for all $d\in \mathbb{R}^n$. If in addition,
\[
Vd - \lim \left\{ V^\prime d^\prime\left|~V^\prime\in \bar\partial \psi (\bar x+ td'), \; d' \rightarrow d, \; t\downarrow 0 \right.\right\} = O\left(\|d\|^2\right)
\]
holds for all $V\in \bar\partial \psi (\bar x+d)$ with $d \rightarrow 0$, then $\psi$ is said to be \emph{strongly semismooth} at $\bar x$. The function $\psi$ will be said to be \textit{SC$^1$} if it is continuously differentiable and $\nabla \psi$ is semismooth. Also,  $\psi$ is \textit{LC$^2$} if it is twice continuously differentiable and $\nabla^2 \psi$ is locally Lipschitzian. These semismoothness concepts for vector-valued functions introduced in  \cite{QiSunANonsmoothVersion1993} are extensions of the original real-valued functions-based ones from \cite{Mifflin1977}.

{ {In the sequel, the concept of basic normal cone will also be necessary. 
 For a closed subset $C$ of $\mathbb{R}^n$, the \textit{basic} (also known as \emph{Mordukhovich} or \textit{limiting}) \textit{normal cone}
 to $C$ at one of its points $\bar x$ is the set
\begin{equation}\label{basic normal cone}
 N_C(\bar x):= \left\{v\in \mathbb{R}^n\left|\; \exists v_k \rightarrow v, \, x_k \rightarrow \bar x \,(x_k\in C): \, v_k\in \widehat{N}_C(x_k)\right.\right\},
\end{equation}
where $\widehat{N}_C$ denotes the polar of the contingent/Boulingand tangent cone to $C$}}:
\begin{equation}\label{polar normal cone}
 \widehat{N}_C(\bar x):=\left\{v\in \mathbb{R}^n \left |\; \langle v, u-\bar x\rangle \leq o(\|u-\bar x\|)\;\; \forall u\in C\right. \right\}.
\end{equation}
$\widehat{N}_C(\bar{x})$ is also known as the \textit{regular} or \textit{Fr\'{e}chet} normal cone. A set $C \subset \mathbb{R}^n$ is said to be \textit{normally regular} at $\bar{x} \in C$ if $ N_C(\bar{x}) = \widehat{N}_C(\bar{x})$. If $C$ is a convex set, it is normally regular at all points of $C$, given that both of the above cones coincide with the normal cone in the sense of convex analysis.

The regular normal cone allows us to introduce the notion of a \textit{regular subdifferential} also known as the Fr\'{e}chet subdifferential. Let the function $\psi : \mathbb{R}^n \rightarrow \mathbb{R}$ be Lipschitz continuous around $\bar x$. Then its regular subdifferential at $\bar{x} \in \mathbb{R}^n$ is given as 
\[
    \hat{\partial}\psi(\bar{x}) = \left\{v\in \mathbb{R}^n\left|~(v, -1) \in \widehat{N}_{\mbox{epi}\psi} (\bar{x}, \psi(\bar{x}))\right.\right\},
\]
where epi$\psi$ denotes the epigraph of the function $\psi$.
Let us now introduce the notion of a \textit{basic subdifferential} which is also known as the \textit{limiting} subdifferential or the \textit{Mordukhovich} subdifferential. For $\psi : \mathbb{R}^n \rightarrow \mathbb{R}$, which Lipschitz continous around $\bar x$, its basic subdifferential at this point is given by 
\[
    \partial\psi(\bar{x}) = \left\{v\in \mathbb{R}^n\left|~(v, -1) \in N_{\mbox{epi}\psi} (\bar{x}, \psi(\bar{x}))\right.\right\}.
\]
Considering the interplay between \eqref{basic normal cone} and \eqref{polar normal cone}, $v \in \partial \psi(\bar{x})$ if and only if there exist sequences $\{x^k \}$ with $x^k \rightarrow \bar{x}$ and $\{v^k\}$ with $ v^k \rightarrow v$ such that $v^k \in \hat{\partial} \psi(x^k)$. The basic subdifferential for a locally Lipschitz function is always non-empty and compact. If $\psi$ is a convex function, then both its regular and basic subdifferentials coincides with the subdifferential in the sense of convex analysis. Furthermore, if a function $\psi : \mathbb{R}^n \rightarrow \mathbb{R}$ is Lipschitz continous around a point $\bar x$, then we have the nice relationship
\begin{equation}\label{LimitingClarke}
    \bar{\partial}\psi = \mbox{cl} \mbox{co}\partial \psi (\bar{x})
\end{equation}
between the Clarke and basic subdifferentials; cf. \cite[Theorem 3.57 of Volume I]{MordukhovichBook2006}, where $\mbox{cl} \mbox{co}$ stands for the closure of the convex hull of the corresponding set. Hence, the Clarke subdifferential is often referred to in some literature as the \textit{convexified subdifferential}.


Next, we provide a formula for the computation of the basic normal cone of a set defined by some functional constraint, which will be used for some key aspects of the analysis in this paper. 
\begin{lemma}[see Theorem 6.10 in \cite{MordukhovichGeneralizedDifferential1994} or Theorem 6.14 in  \cite{RockafellarWetsBook1998}]\label{normal cone estimate}
Let $C:=\Omega \cap \psi^{-1}(\Xi)$, where $\Omega \subseteq \mathbb{R}^n$ and $\Xi\subseteq \mathbb{R}^m$
are closed sets and the function $\psi :\mathbb{R}^n \rightarrow \mathbb{R}^m$ is  Lipschitz  around $\bar{x}$. Furthermore, assume that 
\begin{equation}\label{BCQ}
\left.
\begin{array}{r}
   0\in  \partial \langle v, \psi \rangle (\bar{x})+N_{\Omega}(\bar x)\\
   v\in   N_{\Xi}(\psi(\bar{x}))
\end{array}
\right\} \Longrightarrow v=0.    
\end{equation}
Then the following estimate holds:
\begin{equation}\label{normal cone to operator constraint}
   N_C(\bar x) \subseteq     \bigcup\left\{\partial \langle v, \psi \rangle (\bar x) + N_{\Omega}(\bar x)\left|\; v\in N_\Xi(\psi(\bar x)) \right.\right\}.
\end{equation}
Equality holds in \eqref{normal cone to operator constraint}, provided that $\Xi$ is normally regular at $\psi(\bar x) $ and $ \Omega$ is normally regular at $\bar{x}$.
\end{lemma}
{ {Note that in \eqref{normal cone to operator constraint}, the term $\partial \langle v, \psi \rangle (\bar x)$ refers to the basic subdifferential of the function $x \mapsto \sum^m_{i=1} v_i \psi_i(x)$ at $\bar x$. The qualification condition \eqref{BCQ} is usually referred to as {\it basic constraint qualification} (\textit{BCQ} for short) in the literature. In fact, if $\Xi = -\mathbb{R}^m_+$, then the BCQ is actually equivalent to the Mangasarian-Fromovitz constraint qualification if $\psi$ is a continuously differentiable function.}}
\subsection{Set-valued mappings and their generalized differentiation}\label{Tools from var}
Set-valued mappings and some related concepts will play a key role in the derivation of key results in this paper. To proceed, we consider a set-valued mapping $\Psi :\mathbb{R}^n\rightrightarrows\mathbb{R}^m$. 
{ {The domain of $\Psi$ is denoted as $\text{dom}\Psi$ and defined by
$\text{dom}\Psi := \{ x \in \mathbb{R}^n  \left| \ \Psi(x) \neq \emptyset\right\}$,
while its graph, denoted as $\text{gph} \Psi$, is given by
\begin{equation*}
    \text{gph} \Psi  :=  \{ (x, y) \in \mathbb{R}^m \times \mathbb{R}^m \left | \ y \in \Psi(x) \right \} = \{ (x, y) \in \text{dom} \Psi \times \mathbb{R}^m \left|\ y \in \Psi(x)\right\}.
\end{equation*}
}} 
\noindent$\Psi$ will be said to be  \textit{inner semicompact} at some point $\bar{x}\in \text{dom}\Psi$, if for every sequence $x_k\rightarrow\bar{x}$ with $\Psi(x_k)\neq\emptyset$, there is a sequence of $y_k\in\Psi(x_k)$ that contains a convergent subsequence as $k\rightarrow\infty$. It follows that $\Psi$ is inner semicompactness at $\bar x$ whenever this set-valued mapping is uniformly bounded around $\bar{x}$; i.e., there exists a neighborhood $U$ of $\bar{x}$ and a bounded set $\Theta\subset\mathbb{R}^m$ such that $\Psi(x)\subset\Theta$ for all $x\in U$.

 { {The set-valued mapping $\Psi :\mathbb{R}^n\rightrightarrows\mathbb{R}^m$ is said to be 
 { {inner semicontinuous at $(\bar{x},\bar{y})\in \mbox{gph}\Psi$ w.r.t. $\Omega\in\mathbb{R}^{n}$ whenever for each neighborhood $V\subset \mathbb{R}^{m}$ of $\bar{y}$, there exists a neighborhood $U\subset \mathbb{R}^{n}$ of $\bar{x}$ such that $\Psi (x)\cap V\neq \emptyset$ for all $x \in \Omega \cap U$. If $\Omega = \mathbb{R}^{n}$ can be chosen, $\Psi$ is called inner semicontinuous at $(\bar{x},\bar{y})$ for brevity; i.e., if for every sequence $x_k\rightarrow\bar{x}$, there is a sequence of $y_k\in\Psi(x_k)$ that converges to $\bar{y}$ as $k\rightarrow\infty$. For single-valued mappings $\Psi :\mathbb{R}^n\to\mathbb{R}^m$, this property obviously reduces to the continuity of $\Psi$ at $\bar x$.}} 
 Obviously, the inner semicontinuity property is more restrictive than the inner semicompactness while bringing us to more precise results of the coderivative calculus, as it will be clear in the next section.}} 

As it will become clear in the subsequent sections, some Lipschitz-type properties will also play a central role in the main results. Hence, we introduce some useful concepts.  A set-valued mapping $\Psi : \mathbb{R}^{n}\rightrightarrows \mathbb{R}^{m}$
is said to be \emph{Lipschitz-like} (or satisfies the \textit{Aubin property} \cite{Aubin1984}) around $\left(\bar{x},\bar{y}\right)\in \mbox{gph}\Psi$ if
 there exist neighborhoods $U$ of $\bar{x}$, $V$ of $\bar{y}$, and a constant $l > 0$ such that
\begin{equation}\label{exact-estim}
	\Psi\left( x\right) \cap V \subseteq \Psi\left( u\right) +l\|u-x\|\mathbb{B} \ \ \text{for all} \ x,u \in U
\end{equation}
with $\mathbb{B}$ denoting the unit ball in $\mathbb{R}^{m}$ (see; e.g., \cite{HenrionJouraniOutrataCalmness2002}). This property can be equivalently expressed as 
\begin{equation}\label{Aubin-Property}
    d(y, \Psi(u)) \le l \| x - u \| \,\mbox{ for all }\, x, u \in U \;\, \mbox{and}\;\, y \in V \cap\Psi(x).
\end{equation}
The weaker calmness property can be obtained from \eqref{Aubin-Property} by setting $ u = \bar{x}$. Thus, the set-valued mapping $\Psi :\mathbb{R}^n \rightrightarrows \mathbb{R}^m$ will be said to be \textit{calm} (or satisfy the {\em calmness} property) at the point $(\bar x, \bar y)$ with $\bar y\in \Psi(\bar x)$, if there exist neighborhoods $U$ and $V$ of $\bar x$ and $\bar y$, respectively, and a constant $\ell >0$ such that
\begin{equation}\label{Definition Calmness}
d(y, \Psi(\bar x))\leq \ell \|x-\bar x\| \,\mbox{ for all }\, y\in V\cap \Psi(x), \;\, x\in U.
\end{equation}
In the case where $V=\mathbb{R}^m$, this property goes back to \cite{Robinson1981} who called it
the ``upper Lipschitz property'' of $\Psi$ at $\bar x$. It is proved in \cite{Robinson1981} that the upper Lipschitz (and hence calmness)
 property holds at every point of the graph of a set-valued mapping $\Psi$ if it is polyhedral, i.e., expressible as the union
  of finitely many polyhedral sets. Conditions for the validity of the calmness property and its broad applications to variational analysis and optimization have been developed in \cite{HenrionJouraniOutrataCalmness2002,hen-out01,MOC07} among many other publications. 
  As an important application of calmness, note that based on \cite[Theorem 4.1]{HenrionJouraniOutrataCalmness2002}, inclusion \eqref{normal cone to operator constraint} also holds if the following set-valued mapping is calm at the point $(0, \bar x)$ with $\bar x\in \Omega \cap \psi^{-1}(\Xi)$: 
\begin{equation}\label{CalmnesMapAgain}
\Psi(\varsigma):=\left\{x\in \Omega\left|\;\, \psi(x) + \varsigma \in \Xi\right.\right\}.
\end{equation}

  Next, we use the concept of basic normal cone introduced in the previous subsection to define the notion of {\em coderivative}, which will allow us to present the most powerful characterization of the Lipschitz-likeness property.  For a given set-valued mapping $\Psi :\mathbb{R}^n \rightrightarrows \mathbb{R}^m$,  its coderivative at some point $(\bar{x}, \bar{y})\in \text{gph}\, \Psi$ corresponds to another set-valued mapping 
$D^*\Psi(\bar{x},\bar{y}): \mathbb{R}^m \rightrightarrows \mathbb{R}^n$ that is defined at any $y^*\in \mathbb{R}^m$ by
\begin{equation*}\label{cod-definition}
   D^*\Psi(\bar{x},\bar{y})(y^*):= \left\{x^*\in \mathbb{R}^n|\; (x^*, -y^*)\in N_{\text{gph}\, \Psi}(\bar{x}, \bar{y})\right\}.
\end{equation*}
For a single-valued locally Lipschitz continuous function $\psi : \mathbb{R}^n \rightarrow \mathbb{R}^m$, the coderivative reduces to  
\begin{equation*}\label{coderivative-sub}
    D^* \psi (\bar{x}) (y^*) = \partial \langle y^*,  \psi \rangle (\bar{x}),
\end{equation*}
where similarly to \eqref{normal cone to operator constraint}, the term $\partial \langle y^*, \psi \rangle (\bar x)$ denotes the basic subdifferential of the function $x \mapsto \sum^m_{i=1} y^*_i \psi_i(x)$ at $\bar x$.
The coderivative plays a key role in identifying whether a set-valued mapping is Lipschitz-like; i.e., the well-known \textit{Murdokhovich} criterion (see, e.g., \cite{RockafellarWetsBook1998,MordukhovichGeneralizedDifferential1994}), which says that for a set-valued mapping $\Psi : \mathbb{R}^{n}\rightrightarrows \mathbb{R}^{m}$, which has a closed graph around $(\bar x, \bar y)\in \text{gph}\Psi$, it is Lipschitz-like around this point if and only if 
\begin{equation*} \label{Aubin}
    D^* \Psi ( \bar{x}, \bar{y})(0) = \{0\}.
\end{equation*}
Thanks to the concept of coderivative, we can also characterize the \textit{Lipschitz modulus} of a set-valued mapping when it is Lipschitz-like \eqref{exact-estim}. To present this characterization, we first introduce the notion a Lipschitz modulus of a set-valued mapping $\Psi : \mathbb{R}^{n}\rightrightarrows \mathbb{R}^{m}$. For this, we need  the notion of \textit{outer norm}  
of a positively homogeneous set-valued mapping. { {Note that $\Psi$ is said to be positively homogeneous if $0\in \Psi (0)$ and $\Psi (\lambda x) = \lambda \Psi (x) $ for all $\lambda > 0$ and $ x \in \mathbb{R}^n$.}}
If $\Psi$ is a positively homogeneous
mapping, its \emph{outer norm} is defined by
\begin{equation*}
	\|\Psi\|^{+} :={\displaystyle \sup_{x\in\mathbb{B}}} \
{\displaystyle \sup_{u\in\Psi\left( x\right) }}\|u\|.  
\end{equation*}
Furthermore, the infinimum of all $l > 0$ for which \eqref{exact-estim} holds, also known as the \emph{Lipschitz modulus} of $\Psi$ at $\left(\bar{x},\bar{y}\right)$, is given as follows (cf. \cite[Theorem 9.40]{RockafellarWetsBook1998}),  via the outer norm of the coderivative of the set-valued mapping $\Psi$:
\begin{equation}\label{LipModulus}
	\text{lip} \Psi \left( \bar{x},\bar{y}\right) := \inf \left\{l\in ]0, +\infty[ \ |
\eqref{exact-estim} \ \text{holds for some} \  U \ \text{and} \ V  \right\}=\left\|D^{\ast}\Psi\left( \bar{x},\bar{y}\right)\right\|^{+}.
\end{equation}

Finally, to close this section, we consider some continuous functions $\psi:\mathbb{R}^{n}\times\mathbb{R}^{m}\rightarrow \mathbb{R}^p$ 
and $\phi:\mathbb{R}^{n}\times\mathbb{R}^{m}\rightarrow \mathbb{R}^q$ and associate the set-valued mapping  $\Psi$, which is defined by
		\begin{equation}\label{regularityRCPLD}
			\Psi\left(x\right) := \left\{y\in\mathbb{R}^{m} \left|~\psi\left(x, y\right)  \leq 0, \;\, \phi\left(x, y\right)  = 0\right.\right\}.
		\end{equation}
	Given a point $\left( \bar{x},\bar{y}\right)\in \text{gph}\Psi$ and set $\Omega\subset \mathbb{R}^{n}$, $\Psi$ will
be said to be \emph{R-regular} at $\left( \bar{x},\bar{y}\right)$ w.r.t. $\Omega\subseteq \mathbb{R}^{n}$
if there are some positive numbers $L$, $\epsilon$, and $\delta$ such that
\begin{equation}\label{R-regularityMap}
\begin{array}{rcl}
d\left(y, \, \Psi\left(x\right) \right) 
& \leq & \ L \max \left\{0, \;\,\max\left\{\psi_{i}\left(x, y\right)|\;\, i=1, \ldots, p\right\}, \;\;
\max\left\{\left|\phi_{j}\left(x, y\right)\right|\left|\;\, \right. j=1, \ldots, q\right\}\right\}
\end{array}
\end{equation}
for all $x\in \mathbb{U}_{\delta}\left(\bar{x}\right) \cap \Omega$ and $y\in
 \mathbb{U}_{\epsilon}\left(\bar{y}\right)$. For more details on R-regularity, see \cite{MehlitzMinchenko2020} and references therein.

\section{Coderivative and robust stability of solution maps}\label{Coderivative and robust stability of solution maps}
The plan here is to use the value function reformulation \eqref{VI-1} of  $S$ \eqref{S-Soultion-Map} to construct an estimate for its coderivative  in terms of the derivatives of $f_0$ and $g_0$, and subsequently, to leverage on this result to provide new rules for the Lipschitz-like property of this set-valued mapping. To proceed,  some further definitions and notations are in order. First, the Mangasarian-Fromowitz constraint qualification (MFCQ) for the constraint describing $K$ \eqref{VI},  which will be said to hold at $(\bar x, \bar y, \bar\varsigma)$ if there exists a vector $d\in \mathbb{R}^m$ such that
\begin{equation}\label{LL regularity}
    \nabla_3 g_{0j}(\bar x, \bar y, \bar\varsigma)^\top d < 0\;  
    \mbox{ for all }\; j\in I^3.  
\end{equation}
Recall that $I^3$ is defined in \eqref{I1} and $\nabla_3 g_{0j}(x, y, \varsigma)$ represents the Jacobian of the function $\varsigma \mapsto g_{0j}(x, y, \varsigma)$ for a fixed vector $(x, y)\in \mathbb{R}^n\times \mathbb{R}^m$.
Next, we associate to a vector $(\bar x, \bar y, \bar\varsigma, \varsigma^*)$ a set of  Lagrange multipliers
\begin{equation*}
  \Lambda^{0}(\bar x, \bar y, \bar\varsigma, \varsigma^*)
   :=  \left\{\left(\begin{array}{c}
        v\\ w\\ \lambda
\end{array}\right)\in \mathbb{R}^{2q + 1} \left|
\left.\begin{array}{l}
   \lambda \geq 0,\; v\geq 0, \; g_0(\bar x, \bar y, \bar{y})\leq 0, \; v^\top g_0(\bar x, \bar y, \bar{y})=0,\;  w\in \Lambda(\bar x, \bar y, \bar \varsigma)\\[1ex]
     \varsigma^*  + \lambda\left.\left(f_0(\bar x, \bar y) + \nabla_2f_0(\bar x, \bar y)^\top (\bar y - \bar\varsigma) - \nabla_2 g_{0}(\bar x, \bar y, \bar\varsigma)^{\top}w  \right)\right.\\[1ex]
      \qquad \qquad \qquad \qquad    \qquad \qquad \qquad \qquad + \nabla_2 g_{0}(\bar x, \bar y, \bar y)^{\top}v=0
      \end{array}\right.\right.\right\} 
\end{equation*}
related to the  problem  \eqref{VI}, with $\nabla_2 f_{0}(x, y)$ denoting the Jacobian of the function $f_{0}(x, .)$ for a fixed vector $x\in \mathbb{R}^n$, while $\nabla_2 g_{0}(x, y, \varsigma)$ represents the Jacobian of the function $g_{0}(x, ., \varsigma)$ with $(x, \varsigma)\in \mathbb{R}^n\times \mathbb{R}^m$ fixed. Furthermore, recall that the other involved set of Lagrange multipliers 
\begin{equation}\label{Lambda(x,y)}
\begin{array}{rll}
\Lambda(x, y, \varsigma) &:=& \left\{w\in \mathbb{R}^q\left|\;\nabla_3\ell(x, y, \varsigma, w)=0,\right.\right.\left. \left. w\geq 0, \; g_0(x, y, \varsigma)\leq 0,\; w^{\top}g_0(x, y, \varsigma)= 0 \right.\right\}
\end{array}
\end{equation}
is defined with $x\in \mathbb{R}^n$ and $(y,\varsigma)\in \mathbb{R}^{2m}$. In \eqref{Lambda(x,y)}, the function $\ell$ represents the Lagrangian function associated with the underlying parametric optimization problem \eqref{varphi}:
 \begin{equation}\label{ell}
    \ell(x, y, \varsigma, w):=\varsigma^\top f_0(x, y)+ w^\top g_0(x, y, \varsigma).
 \end{equation}
Furthermore, in \eqref{Lambda(x,y)},  $\nabla_3 \ell(x, y, \varsigma, w) := f_0(x, y) + \nabla_3 g_{0}(x, y, \varsigma)^\top w$ with $\nabla_3 g_{0}(x, y, \varsigma)$ denoting the Jacobian of the function $\varsigma \mapsto g_{0j}(x, y, \varsigma)$
when the vectors $x\in \mathbb{R}^n$ and $y\in \mathbb{R}^m$ are fixed.

 Throughout, the linear independence constraint qualification (LICQ) is satisfied at $(\bar x, \bar y, \bar\varsigma)$ if it holds that
 \begin{equation}\label{LICQ}
 \mbox{the family of vectors }\;\, \left\{\nabla_3 g_{0j}(\bar x, \bar y, \bar\varsigma)\left|\;\, j\in I^3\right.\right\} \;\, \mbox{ is linearly independent},
 \end{equation}
recalling that $I^3$ is given in \eqref{I1} and $\nabla_3 g_0(x, y,\varsigma)$ stands for the Jacobian  of the function $g_0$ w.r.t.
  $\varsigma$. And before we can state the first result of this section, note that for a fixed parameter $x\in \mathbb{R}^n$, the set-valued mapping $S$ \eqref{S-Soultion-Map} can be written as
$S(x)=\left\{y\in \mathbb{R}^m|\, y\in S_0(x, y)\right\}$,
where  $S_0 : \mathbb{R}^n\times \mathbb{R}^m \rightrightarrows \mathbb{R}^m$ is given as
\begin{equation}\label{S-v}
 S_0(x,y)  := \arg\underset{\varsigma}\min~\left\{\varsigma^\top f_0(x,y)\left|\; \varsigma\in K(x, y)\right.\right\}.
\end{equation}
\begin{theorem}\label{I}
Consider a point $(\bar x, \bar y)$ such that $\bar y\in S(\bar x)$ and let $\bar{\varsigma} \in S_0(\bar x, \bar y)$.
\begin{itemize}
  \item[\emph{(i)}]Suppose that the MFCQ \eqref{LL regularity} holds at $(\bar x, \bar y, \bar y)$ and $(\bar x, \bar y, \bar{\varsigma})$.
  Furthermore, assume that $S_0$ \eqref{S-v} is inner semicontinuous at $(\bar x, \bar y, \bar{\varsigma})$ and let the set-valued mapping
\begin{equation}\label{PsiCalm}
\Psi(\theta):=\left\{(x,y)\in \mathbb{R}^n \times \mathbb{R}^m\left|\; g_0(x, y, y)\leq 0, \;\; y^\top f_0(x, y) - \varphi(x, y)\leq \theta \right.\right\}
\end{equation}
 be calm at $(0, \bar x, \bar y)$. Then for all $y^*\in \mathbb{R}^m$, it holds that 
{\small{
\begin{equation}\label{CoD S}
\begin{array}{l}
  D^*S(\bar x, \bar y)(y^*) \subseteq
   \underset{(v, w, \lambda) \in \Lambda^0 (\bar x, \bar y, \bar{\varsigma}, y^*)} \bigcup
 \left\{\nabla_1 g_0(\bar x, \bar y, \bar y)^\top v + \lambda \left(\displaystyle 
  \nabla_1 f_{0}(\bar x, \bar y)^\top  (\bar y- \bar{\varsigma})  - \nabla_1 g_0(\bar x, \bar y, \bar{\varsigma})^\top w\right)\right\}.
\end{array}
\end{equation}}}
\item[\emph{(ii)}] Suppose that in addition to the assumptions in \emph{(i)}, the following qualification condition holds:
{\small{
\begin{equation}\label{Stab CQ-1}
\begin{array}{l}
(v, w, \lambda) \in \Lambda^0(\bar x, \bar y, \bar{\varsigma}, 0) \;\; \Longrightarrow \;\;  
  \left[\nabla_1 g_0(\bar x, \bar y, \bar y)^\top v
 + \lambda\left(\displaystyle \nabla_1 f_{0}(\bar x, \bar y)^\top  (\bar y- \bar{\varsigma})
    - \nabla_1 g_0(\bar x, \bar y, \bar{\varsigma})^\top w\right) = 0\right].
\end{array}
\end{equation}}}
${}$\\
Then the set-valued mapping $S$ \eqref{S-Soultion-Map}  is Lipschitz-like around $(\bar x, \bar y)$ and its Lipschitz modulus at this point has the following upper estimate:
\[
\text{\emph{lip}}\, S(\bar x, \bar y) \leq \underset{\underset{(v, w, \lambda)\in \Lambda^0(\bar x, \bar y, \bar{\varsigma}, y^*)}{\|y^*\|\leq 1}}{\sup}
\left\{\left\|\nabla_1 g_0(\bar x, \bar y, \bar y)^\top v + \lambda\left(\nabla_1 f_{0}(\bar x, \bar y)^\top  (\bar y- \bar{\varsigma})
-\nabla_1 g_0(\bar x,\bar y,\bar{\varsigma})^\top w\right)\right\|\right\}.
\]
\end{itemize}

\end{theorem}
\begin{proof}
To begin, note that with the expression of the set-valued mapping $S$ \eqref{VI-1}, we have  
\[
\begin{array}{l}
\mbox{gph}S = \Omega \cap \psi^{-1}(\mathbb{R}_-)\,
\mbox{ with } \, \Omega:=\left\{(x, y)\in \mathbb{R}^n\times \mathbb{R}^m\left|\, \tilde{g}(x, y)\leq 0 \right.\right\},\, \psi(x, y):= \tilde{f}(x, y) - \varphi(x, y),
\end{array}
\]
where $\tilde{g}(x,y):= g_0(x, y, y)$ and $\tilde{f}(x, y):=y^\top f_0(x, y)$. Then observe that the function $\varphi$ is Lipshitz continuous around  $(\bar x, \bar y)$ given that the set-valued mapping $S_0$ \eqref{S-v} is inner semicontinuous at $(\bar x, \bar y, \bar{\varsigma})$ and the MFCQ is satisfied at this same point. Hence, from Lemma \ref{normal cone estimate}, with the calmness of the set-valued mapping \eqref{PsiCalm}, as corresponding expression of \eqref{CalmnesMapAgain}, it follows that 
\begin{equation} \label{estimate-eq-1}
    N_{\text{gph}S} ( \bar{x}, \bar{y}) \subseteq \bigcup \left \{\partial \langle\lambda, \psi\rangle (\bar{x}, \bar{y}) + N_\Omega (\bar{x}, \bar{y}) \left |\; \lambda \ge 0\right.\right\},
\end{equation}
while considering the fact that $\psi(\bar x, \bar y)=0$, as $\bar y\in S(\bar x)$, and the functions $\tilde{f}$ and $\tilde{g}$ continuously differentiable. Now, considering the definition of the coderivative and inclusion \eqref{estimate-eq-1}, it follows that for all $y^*\in \mathbb{R}^m$, 
\begin{equation} \label{estimate-eq-2}
  D^*S(\bar{x}, \bar{y}) (y^*) \subseteq \{ x^* \in \mathbb{R}^n \left |\; \exists \lambda \ge 0:\quad  \left(x^*, -y^*\right)\in \partial \langle\lambda, \psi\rangle (\bar{x}, \bar{y}) + N_\Omega(\bar{x}, \bar{y}) \right \}.
\end{equation}
Since the MFCQ also holds at $(\bar{x}, \bar{y}, \bar{y})$, using Lemma \ref{normal cone estimate} one more time, 
\begin{equation*}\label{NOmega}
    N_\Omega(\bar{x}, \bar{y}) \subseteq \left\{\nabla g(\bar{x}, \bar{y})^\top v \left|\; v \ge 0,\; g(\bar{x}, \bar{y}) \le 0, \; v^T g(\bar{x}, \bar{y})=0\right.\right\}.
\end{equation*}

On the other hand,  observe that 
\begin{equation*}\label{partialLP}
 \partial \langle\lambda, \psi\rangle (\bar{x}, \bar{y}) \subseteq \lambda \nabla f(\bar{x}, \bar{y}) + \partial \langle \lambda, \, -\varphi\rangle (\bar{x}, \bar{y}) \subseteq \lambda \nabla f(\bar{x}, \bar{y}) -\lambda \bar{\partial} \varphi (\bar{x}, \bar{y})  
\end{equation*}
given that as $\varphi$ is Lipschitz continuous around $(\bar x, \bar y)$, its limiting subdifferential is included in the Clarke subdifferential, according to \eqref{LimitingClarke}. Furthermore, it follows from \cite[Corollary 5.3]{MordukhovichNamPhanVarAnalMargBlP} that the combination of the inner semicontinuity of $S_0$ and MFCQ \eqref{LL regularity} both holding at the point $(\bar x, \bar y, \bar{\varsigma})$ leads to
\begin{equation}\label{SubPhi}
\begin{array}{c}
\bar\partial \varphi(\bar x, \bar y) \;\, \subseteq \;\, \left\{\left.\nabla f_0(\bar x, \bar y)^\top\bar{\varsigma}  +  \nabla_{1, 2} g_0(\bar x, \bar y, \bar{\varsigma})^\top w\right|\;\, w\in \Lambda (\bar x, \bar y, \bar{\varsigma})\right\}.
\end{array}
\end{equation}

Also note that as $\tilde{f}(x, y) = y^T f_0(x, y)$, it holds that
\begin{equation}\label{NablaF}
  \nabla \tilde{f}(x,y)= f_0(x,y) + \nabla f_0(x, y)^\top y.
\end{equation}
Now, let $x^* \in D^*S(\bar{x}, \bar{y})(y^*)$. Then it follows from \eqref{estimate-eq-2}--\eqref{NablaF} that we can find $\lambda$, $v$, $w$ such that 
\[
\begin{array}{r}
\lambda \geq 0,\; v\geq 0, \; g_0(\bar x, \bar y, \bar{y})\leq 0, \; v^\top g_0(\bar x, \bar y, \bar{y})=0,\;  w\in \Lambda(\bar x, \bar y, \bar\varsigma),\\[1.5ex]
\left[\begin{array}{r}
     x^*\\
     -y^*
\end{array}\right]=\left[\begin{array}{r}
 \lambda \left(
  \nabla_1 f_{0}(\bar x, \bar y)^\top  (\bar y- \bar{\varsigma})  - \nabla_1 g_0(\bar x, \bar y, \bar{\varsigma})^\top w\right) + \nabla_1 g_0(\bar x, \bar y, \bar y)^\top v\\[1.5ex]
  \lambda\left.\left(f_0(\bar x, \bar y) + \nabla_2f_0(\bar x, \bar y)^\top (\bar y - \bar\varsigma) - \nabla_2 g_{0}(\bar x, \bar y, \bar\varsigma)^{\top}w  \right)\right. + \nabla_2 g_{0}(\bar x, \bar y, \bar y)^{\top}v
\end{array}\right].
\end{array}
\]
Considering this and the definition of $\Lambda^0 (\bar x, \bar y, \bar{\varsigma}, y^*)$, we have inclusion \eqref{CoD S}.

For assertion (ii), observe that from \eqref{CoD S}, we have the following inclusion:
\[
\begin{array}{ll}
  D^*S(\bar x, \bar y)(0) \subseteq
  & \underset{(v, w, \lambda) \in \Lambda^0 (\bar x, \bar y, \bar\varsigma, 0)}  \bigcup
  \left\{\lambda\left(\nabla_1 f_{0}(\bar x, \bar y)^\top (\bar y - \bar\varsigma) - \nabla_1 g_0(\bar x, \bar y, \bar\varsigma)^\top w\right) + \nabla_1 g_0(\bar x, \bar y, \bar y)^\top v\right\}.
\end{array}
\]
Hence, it is clear that condition \eqref{Stab CQ-1} is sufficient for $D^*S(\bar x, \bar y)(0) \subseteq \{0\}$. Thus, since the coderivative mapping is positively homogeneous (i.e., in particular, we have $0\in D^*S(\bar x, \bar y)(0)$),  the coderivative/Mordukhovich criterion $D^*S(\bar x, \bar y)(0) = \{0\}$ is satisfied. This implies that $S$ is Lipschitz-like around $(\bar x, \bar y)$, and it follows from \eqref{LipModulus} that the exact Lipschitzian bound of $S$ around $(\bar x, \bar y)$ can be obtained as
$$
\text{{lip}}\, S(\bar x, \bar y) = \|D^*S(\bar x, \bar y)\|^+:=\sup\left\{\|u\|\left|\, u\in D^*S(\bar x, \bar y)(y^*), \, \|y^*\|\leq 1\right.\right\}.
$$
Considering inclusion \eqref{CoD S} one more time, we have the stated upper bound on the quantity  $\text{{lip}}\, S(\bar x, \bar y)$. \qed
\end{proof}

As it is clear from the proof of this theorem, the calmness of the set-valued mapping $\Psi$ \eqref{PsiCalm} is crucial.
 Hence, we are now going to provide some sufficient conditions to ensure that it holds. The first sufficient condition
  is based on the concept of \emph{uniform weak sharp minimum}, which is well-known in parametric optimization;
  see, e.g., \cite{DempeZemkohoGenMFCQ,Ye1998New}. Here, we introduced a version of the concept tailored to the parametric  problem  \eqref{VI}.
\begin{definition}\label{DefLUWSMC}
	Problem \eqref{VI} will be said to satisfy a \emph{local uniform weak sharp minimum condition} (LUWSMC)
	at $(\bar x, \bar y)\in \mbox{gph} \ S$ if there exist some $\alpha>0$ and $\varepsilon> 0$ such that
	\begin{equation}\label{LUWSMC}
 \begin{array}{l}
\forall \left( x,y\right) \in \mathbb{U}_{\epsilon}\left( \bar{x},\bar{y}\right): \quad \left[
y\in K\left( x,y\right)\; \Longrightarrow \; 	d(y, \, S(x))\leq \alpha \left(f(x,y) - \varphi(x, y)\right)\right].
 \end{array}
	\end{equation}
\end{definition}

{ {We can now state the following result, first reported in \cite{HenrionSurowiecCalmness2011} when $K$ is a fixed mapping. However, the proof here is essentially adapted from this paper.}}
\begin{theorem}\label{sharp-one}If the LUWSMC holds at  $(\bar x, \bar y)\in \mbox{gph}S$, then the set-valued mapping $\Psi$ \eqref{PsiCalm} is calm at  $(0, \bar x, \bar y)$.
\end{theorem}
\begin{proof} 
Obviously, in this context, we have $\Psi(0) = \mbox{gph}S$. For some numbers $\alpha>0$ and $\varepsilon> 0$ such that \eqref{LUWSMC} holds, consider any $(x, y, \theta)\in \mathbb{U}_{\epsilon}\left(\bar x, \bar y ,0\right)$ such that $(x, y)\in \Psi(\theta)$. Then we have 
	\[
 \begin{array}{c}
 	d\left((x, y), \Psi(0)\right)  =  d\left((x, y),  \mbox{gph}S\right)
                                      \leq  d(y, S(x))\leq \alpha \left(f(x,y) - \varphi(x, y)\right)
                                     \leq  \alpha \theta = \alpha |\theta -0|,
 \end{array}
	\]
	where the last equality results from the fact that  $(x, y)\in \Psi(\theta)$ implies $\theta \geq 0$, as, by the definition of $\varphi$, it holds that $f(x, y) \geq \varphi(x, y)$ for all $y\in K\left( x,y\right) $. \qed
\end{proof}

Next, we provide another sufficient condition based on the R-regularity concept introduced in Subsection \ref{Tools from var}, but now
 for the fulfillment of the LUWSMC. { {To proceed, consider the expression of $S$ in \eqref{VI-1}, and note that
the R-regularity constraint qualification (RRCQ) will be said to hold at the point $\left(\bar{x}, \bar{y}\right)\in \text{gph}\,S$ if $S$ is R-regular \eqref{R-regularityMap} at $\left( \bar{x}, \bar{y}\right)$ w.r.t. $\text{dom}S$}}. Given that the notion of R-regularity for a parametric system is a generalization of that  of uniform error bound, the proof of this result is inspired by \cite{Ye1998New}.  

\begin{proposition}\label{RRCQimpliesLUWSMC}
	If the RRCQ holds at { {$\left( \bar{x}, \bar{y}\right)\in \text{dom}S$}} and there are some neighborhoods 
$U\subset\mathbb{R}^{n}$ of $\bar{x}$ and $V\subset\mathbb{R}^{n}$ of $\bar{y}$ such that
$\text{dom}K\cap \left(U\times V\right) =  \left(\text{dom}S\cap U\right) \times V$, then the
LUWSMC is satisfied at $\left( \bar{x}, \bar{y}\right)$.
\end{proposition}
\begin{proof}
	Fix $\left( \bar{x}, \bar{y}\right) \in  \text{gph}S$. Since the mapping $S$ is R-regular at $\left( \bar{x}, \bar{y}\right)$ w.r.t. $\text{dom} \ S$, there exist  $\sigma > 0$ and $\epsilon > 0$ such that for all $\left( x,y\right) \in \mathbb{U}_{\epsilon}\left(\bar{x}, \bar{y}\right)\cap \left( \text{dom}S \times V\right)$, we have the inequality
	\begin{equation*}
		d\left( y,S\left( x\right) \right) \ \leq \ \sigma \ \max \{0, \, \tilde{g}\left(x,y\right), \, \tilde{f}\left(x,y\right)-\varphi\left( x,y\right)\},
	\end{equation*}
 where similarly to the notation used in the proof of Theorem \ref{I}, $\tilde{g}(x,y):= g_0(x, y, y)$ and $\tilde{f}(x, y):=y^\top f_0(x, y)$. 
	From the definition of the optimal value function, for any $\left( x,y\right) \in \mathbb{U}_{\epsilon}\left( \bar{x}, \bar{y}\right)$ with $\tilde{g}\left( x,y\right)\leq 0$, we have the inequality $\tilde{f}\left( x,y\right)-\varphi\left(x, y\right) \geq 0$. Hence, for all $\left( x,y\right) \in \mathbb{U}_{\epsilon}\left( \bar{x}, \bar{y}\right)\cap \left(\text{dom}S\cap U\right) \times V$, one gets
	\begin{equation}\label{QueqRF}
	y\in K\left( x,y\right) \Longrightarrow d\left(y, S\left(x\right) \right) \ \leq \ \sigma\left(\tilde{f}\left(x,y\right)-\varphi\left( x,y\right)\right).
	\end{equation}
Since $\text{dom}K\cap \left( U\times V\right) =  \left(\text{dom}S\cap U\right) \times V$, condition \eqref{QueqRF} holds for all $\left( x,y\right) \in \mathbb{U}_{\epsilon}\left( \bar{x}, \bar{y}\right)\cap \left( U\times V\right)$. \qed
\end{proof}

Note that R-regularity for a set-valued mapping as stated in \eqref{regularityRCPLD} is guaranteed under the validity of the
corresponding MFCQ. However, based on \cite{DempeZemkohoGenMFCQ,YeZhuOptCondForBilevel1995}, the MFCQ will automatically fail for the constraints describing the set-valued mappings $S$ \eqref{VI-1} and $S_{0}$ \eqref{S-v} because of the presence of the with the optimal value function. To overcome this failure, we next provide a more tractable framework for the fulfilment of the condition, based on the \textit{relaxed constant rank constraint qualification (RCRCQ)}, introduced in \cite{Andreani-2012} and studied in various papers, including \cite{MehlitzMinchenko2020-100}, which inspired the next result. 
To proceed, note that if we assume that there exists a neighborhood $U$ of the point $\left(\bar{x},\bar{y}\right)$ in which the value function $\varphi$ \eqref{varphi} is differentiable, we can introduce the matrix 
\begin{equation*}
	\mathcal{B}\left( x,y\right)  :=
	\left[
	\begin{array}{c}
		\nabla_2 \tilde{f}(\bar x, \bar y) \\[1ex]
  \nabla_2 \tilde{g}(x,y)_{I^2}
	\end{array}
	\right],
\end{equation*}
where the functions $\tilde{g}(x,y):= g_0(x, y, y)$ and $\tilde{f}(x, y):=y^\top f_0(x, y)$ are ones with notation first used in the proof of Theorem \ref{I}. Here, $I^2$ denotes the index set defined in \eqref{I1}. The RCRCQ will be said to hold at a point $\left(\bar{x},\bar{y}\right)\in \mbox{gph}S$ if there is a neighborhood $V$ of $\left( \bar{x},\bar{y}\right)$ such that for each index set $J\subset \{1,\cdots , |I^2|+1 \}$,  the matrix $\mathcal{B}\left( x,y\right)_{J}$ has constant row rank for all $\left( x,y\right) \in V$. 
 
\begin{proposition}\label{vertix} { {Consider problem \eqref{VI} and let its solution set-valued mapping $S$ \eqref{VI-1} be inner semicontinuous at $\left(\bar{x},\bar{y}\right)\in \mbox{gph}S$. Furthermore, suppose that there exists a neighborhood $U$ of the point $\left( \bar{x},\bar{y}\right)$ in which $\varphi$ \eqref{varphi} is differentiable, and let the RCRCQ be satisfied at $\left(\bar{x},\bar{y}\right)$. Then $S$ \eqref{VI-1} is R-regular at the point $\left(\bar{x}, \bar{y}\right)$.}}
\end{proposition}
\begin{proof}  Follows immediately from 
   \cite[Theorem 4.2]{MehlitzMinchenko2020-100}; also see \cite{MehlitzMinchenko2020} for futher details and discussions. \qed
\end{proof}

Next, we provide an example where all the assumptions of Proposition \ref{vertix} are satisfied at a suitable selection of points on the graph of the solution set-valued mapping of the corresponding QVI. 
\begin{example}
Consider a  problem  \eqref{VI} with $f_0(x,y):=Cy$ and $g_0(x,y, \varsigma):=A(x,y)\varsigma - b(x,y)$, where  
\begin{equation}\label{QVI-Ex}
	C := \left. \left( \begin{array}{lcl}
		       1 & \;\, & 0 \\
		       0 &  \;\, & 2
		       \end{array}
		       \right.
		       \right),  \ \ \ \  A\left(x,y\right) := \left. \left(\begin{array}{ccc}
		       	y_{1} & \ & 0 \\
		       	0      &  \ & 3y_{2}
		       \end{array}
		       \right.
		       \right), \ \ \mbox{ and } \ \ b\left( x,y\right) = \left. \left( \begin{array}{l}
		       	x^{2}_{1} \\[1ex]
		       	x^{2}_{2}
		       \end{array}
		       \right.
		       \right).
\end{equation}
The primal and dual optimization problems associated to the QVI described in \eqref{QVI-Ex} can be written as
\begin{equation*}
	\left(P_{x,y}\right) \ : \ \left\{
	\begin{array}{l}
	{\displaystyle \min_{\varsigma}} \ \varsigma^{T}f_{0}\left( x,y\right) =\varsigma_{1}y_{1}+2\varsigma_{2}y_{2} \\
	\mbox{s.t.} \quad -y_{1}\varsigma_{1} \geq - x^{2}_{1}, \;\,-3y_{2}\varsigma_{2} \geq -x^{2}_{2},
    \end{array}
    \right. \quad \mbox{ and } \quad \left(D_{x,y}\right) \ : \ \left\{
    \begin{array}{l}
    	{\displaystyle \max_{\sigma}} \ -x^{2}_{1}\sigma_1 - x^{2}_{2} \sigma_2 \\
    	\mbox{s.t.} \quad -\sigma_{1}y_{1}=y_{1},\;\, -3\sigma_{2}y_{2}=2y_{2}, \;\, y\geq 0,
    \end{array}
    \right.
\end{equation*}
respectively. We can easily check that for any $(x, y)$ such that $y >0$, the points
\[
\varsigma = \left(\dfrac{x^{2}_{1}}{y_{1}},\;\, \dfrac{x^{2}_{2}}{3y_{2}}\right)
\;\, \mbox{ and }
\;\,
\sigma=\left(-1, \;\, -\dfrac{2}{3} \right)
\]
are optimal for  $\left(P_{x,y}\right)$ and $\left(D_{x,y}\right)$, respectively.
It is therefore clear that for any $(x, y)$ such that $y>0$, the corresponding optimal value function  $\varphi$ \eqref{varphi} is continuously differentiable and we have  
\begin{equation*}
\varphi\left( x,y\right) =
x^{2}_{1} + \dfrac{2}{3}x^{2}_{2} \;\; \mbox{ and }\;\;	S\left( x\right) =\left. \left\{ y\in \mathbb{R}^{2}\left|
	                \begin{array}{lcl}
		 &  & y^{2}_{1}+2y^{2}_{2} - x^{2}_{1}  - \dfrac{2}{3}x^{2}_{2}\leq 0 \\
		 &  & y^{2}_{1}-x^{2}_{1} \leq 0, \;\; 3y^{2}_{2}-x^{2}_{2} \leq 0
	\end{array}
		       \right.
               \right.
		       \right\}.
\end{equation*}
Moreover, for any $(x, y)$ such that $y>0$ the following matrix has a constant rank. 
\begin{equation*}
	\mathcal{B}\left( x,y\right)  :=
	\left[
	\begin{array}{cc}
		2y_1 & 4y_2 \\[1ex]
        -2y_1 & 0\\[1ex]
        0 & -6y_2
	\end{array}
	\right].
\end{equation*}
Hence, the RCRCQ will be said to hold at a point $\left(x, y\right)\in \mbox{gph}S$ such that $y>0$.

Finally, we can also check that the family of points 
\[
x:=(x_1, \; x_2), \;\, y:=\left(x_1, \;\, \frac{x_2}{\sqrt{3}}\right) \;\mbox{ wit }\; x>0
\]
is such that $y\in S(x)$ and $y>0$ and the set-valued mapping $S$ \eqref{VI-1} is inner semicontinuous at all such points. \qed
\end{example}		

Note that all the relationships discussed above can be summarized in the following diagram:
\begin{figure*}[htp]
\begin{center}
\[
\text{RCRCQ}_{\Omega} \quad \xRightarrow{ \text{(Proposition \ref{vertix})} } \text{RRCQ}_{\Omega} \quad \xRightarrow{ \text{(Proposition \ref{RRCQimpliesLUWSMC})} } \quad \text{LUWSMC}_{\Omega} \quad \xRightarrow{ \text{(Theorem \ref{sharp-one})} }\quad \text{CALM}_{\Omega}
\]
${}$\\[-5ex]
\caption{\emph{Here, RCRCQ represents the regularity condition and corresponding framework in Proposition \ref{vertix}. For the definition of the RRCQ, see the discussion {{just before}} Proposition \ref{RRCQimpliesLUWSMC}, while the LUWSMC is described in Definition \ref{DefLUWSMC}. As for CALM, it refers to the calmness of the set-valued mapping $\Psi$ \eqref{PsiCalm} imposed in Theorem \ref{I}. Note that the index $\Omega$ is used here just to differentiate the set-valued mapping $\Psi$ \eqref{PsiCalm} (based on $\Omega$) from a slightly different mapping that will be introduced in the next section.}}\label{Figure1} 
\end{center}
\end{figure*}

A few final comments on the assumptions in Theorem \ref{I} are in order. First, it is important to recall that there are various works in the literature where sufficient conditions to ensure the continuous differentiability of an optimal value function in the form $\varphi$ \eqref{varphi} are given; see, e.g., \cite{BonnansShapiro2000,Fiacco1983,MehlitzZemkoho2021} and references therein. 
With regards to the calmness assumption, note that a parallel assumption can be made in the context of the generalized equation reformulation in \eqref{S-KKT}; see \cite[Theorem 3.1(ii)]{MOC07}.  However, it was shown in \cite{HenrionSurowiecCalmness2011} that in the context of the feasible set of a scenario of a bilevel optimization problem with the unperturbed lower-level feasible set, which can have some similarity with a case of the  problem  in \eqref{VI} (see discussion in next section), the corresponding version of the calmness condition for the reformulation in \eqref{S-KKT} shows a behavior better than the calmness of its version of the set-valued mapping $\Psi$ \eqref{PsiCalm}. Nevertheless, our value function reformulation \eqref{VI-1} of the  problem  \eqref{VI} has at least two key advantages over its generalized equation model in \eqref{S-KKT}: (1) it does not require convexity and (2) it does not require second order derivatives for the stability analysis in Theorem \ref{I}. We will see in the next section that these two points are  important in the design of a second order method for \eqref{OPQVIC}.

To close this section, we discuss the inner semicontinuity the set-valued mapping $S_0$ \eqref{S-v} also required in Theorem \ref{I}. First, note that this assumption can be weakened to the \emph{inner semicompactness}. However, applying this weaker assumption will lead to a very loose upper bound on the coderivative of $S$ in \eqref{CoD S} due the convex hull operator that will appear { {in the upper}} estimate of the subdifferential of  $\varphi$  \eqref{SubPhi}; see, e.g., \cite{DempeDuttaMordukhovichNewNece,DempeZemkohoGenMFCQ}. This would subsequently lead to a stronger assumption in \eqref{Stab CQ-1} as well as a looser upper bound for $\text{lip}\,S(\bar x, \bar y)$ obtained in Theorem \ref{I}. 
The following adaptation of \cite[Lemma 2.2]{MehlitzMinchenko2020} provides a scenario where $S_0$ \eqref{S-v} is inner semicontinuous; further discussions on the topic can be found in the latter paper and references therein. 

\begin{proposition}\label{InnerSemS}
Let the set-valued mapping $S_{0}$ \eqref{S-v}, written in the form 
\[
S(x, y):=\left\{\varsigma\in \mathbb{R}^m\left|\;\, \varsigma^\top f_0(x,y)-\varphi(x,y)\leq 0, \;\, g_0(x, y, \varsigma)\leq 0\right.\right\},
\]
be R-regular at the point $\left( \bar{x},\bar{y},\bar{\varsigma}\right) \in  \text{gph}S_{0}$ w.r.t. $\text{dom}S_{0}$. Moreover, let the set-valued mapping $K$, as defined in \eqref{VI}, be locally bounded at $\left(\bar{x},\bar{y}\right) \in  \text{dom}S_{0}$ and inner semicontinuous at $\left(\bar{x},\bar{y},\bar{\varsigma}\right)$ w.r.t. $\text{dom}K$. Then the set-valued mapping $S_{0}$ is inner semicontinuous at $\left(\bar{x},\bar{y},\bar{\varsigma}\right)$ w.r.t. $\text{dom}S_{0}$.
\end{proposition}
\begin{proof}
With the notation
$
\overline{g}_{0}\left( x,y,\varsigma\right) :=\varsigma^{T}f_{0}\left( x,y\right)-\varphi\left( x,y\right)
$
and
$
\overline{g}_{j}\left( x,y,\varsigma\right):= g_{j}\left( x,y,\varsigma\right)$
 for {  {$j\in \{1, \ldots, q\} $}}, the set-valued mapping $S_0$ can be rewritten as
\begin{equation*}
	S_{0}\left( x,y\right) =\{\varsigma \in \mathbb{R}^{m}|\;\, \overline{g}_{0_{j}}\left( x,y,\varsigma\right)\leq 0, \ j\in J\cup \{0\}\}.
\end{equation*}
Let $\left( x_{k},y_{k}\right)_{k}\subset \text{dom} \ S_{0}$ such that $x_{k}\rightarrow \bar{x}$ and $y_{k}\rightarrow \bar{y}$.
First, we claim that $\varphi$ is continuous at $\left(\bar{x},\bar{y}\right)$. Indeed, by definition, the functions  $\overline{g}_{1}$, \ldots, $\overline{g}_{p}$ are continuous (as they are continuously differentiable), then $K$ is upper semicontinuous at $\left( \bar{x},\bar{y}\right)$. Since $f_{0}$ is continuous, it follows from \cite[Theorem 4.2.1]{bank-1983} that $\varphi$ is {  {lower semicontinuous}} at $\left( \bar{x},\bar{y}\right)$. Now, under {  {the inner semicontinuity of $K$ at $\left(\bar{x},\bar{y},\bar{\varsigma}\right)$}}, we get from \cite[Theorem 4.2.1]{bank-1983}, that $\varphi$ is upper semicontinuous at $\left( \bar{x},\bar{y}\right)$. Consequently, $\varphi$ is continuous at $\left( \bar{x},\bar{y}\right)$.
On the other hand, the assumption of the proposition guarantees the existence of a constant $L> 0$, $\delta > 0$, and some $\epsilon >0$ such that
	\begin{equation*}
	d\left(\varsigma, \, S_{0}\left(x,y\right) \right) \ \leq \ L \max \left\{ 0, \;\,\max\left\{\overline{g}_{j}\left(x, y,\varsigma\right)|\;\, j\in J \cup \{0\} \right\} \right\}
	\end{equation*}
for all $\left( x,y\right) \in \mathbb{U}_{\epsilon}\left(\bar{x},\bar{y}\right)$ and all $\varsigma \in \mathbb{U}_{\delta}\left(\bar{\varsigma}\right)$. Consequently,
	\begin{equation*}
		d\left(\bar{\varsigma}, \, S_{0}\left(x,y\right) \right) \ \leq \ L \max \left\{ 0, \;\,\max\left\{\overline{g}_{j}\left(x, y,\bar{\varsigma}\right)|\;\, j\in J \cup \{0\} \right\} \right\} \;\mbox{ for all }\; \left( x,y\right) \in \mathbb{U}_{\epsilon}\left(\bar{x},\bar{y}\right).
	\end{equation*}
Therefore, by continuity of $\overline{g}_{0}\left( x_{k},y_{k},\cdot\right),\,\overline{g}_{1}\left( x_{k},y_{k},\cdot\right), \ldots , \overline{g}_{q}\left(x_{k},y_{k},\cdot\right)$ and the choice of $\left( x_{k},y_{k}\right)_{k}\subset \text{dom}S_{0}$, the set $S_{0}\left(x_{k},y_{k}\right) $ is nonempty and closed. Hence, there exist $\varsigma_{k}\in \Pi\left(\bar{\varsigma},S_{0}\left(x_{k},y_{k}\right)\right)$ for sufficiently large $k$ such that	
\begin{equation*}
	\|\bar{\varsigma}-\varsigma_{k} \| \ \leq \ L \max \left\{ 0, \;\,\max\left\{\overline{g}_{j}\left(x_{k}, y_{k},\bar{\varsigma}\right)|\;\, j\in J \cup \{0\} \right\} \right\}.
\end{equation*}
{{Note that $\Pi\left( x,A\right)$ stands for the projection of $x$ onto $A$.}}
Hence, the continuity of {$\overline{g}_{0},\overline{g}_{1}, \cdots, \overline{g}_{q}$} ensures that ${\displaystyle \lim_{k\rightarrow +\infty}} \ \|\bar{\varsigma}-\varsigma_{k} \|=0$. Consequently $S_{0}$ is inner semicontinuous at  $\left(\bar{x},\bar{y},\bar{\varsigma}\right)$ w.r.t. $\text{dom}S_{0}$.	\qed
\end{proof}

\section{Optimization problems with a QVI constraint}
\label{Optimization problems with a QVI constraint}
Our primary goal in this section is to develop necessary optimality conditions and a solution algorithm for problem \eqref{OPQVIC}. Considering the value function reformulation \eqref{VI-1} of   \eqref{VI}, the problem can be rewritten as
\begin{equation}\label{OPQVIC-LLVF-K}\tag{OPQVI-R}
 \min~F(z)  \;\; \mbox{ s.t. } \; G(z)\leq 0, \; g(z)\leq 0, \;  f(z) - \varphi(z)\leq 0
\end{equation}
with the notation $z:=(x, y)$, $g(x,y):= g_0(x, y, y)$, and $f(x,y) := y^\top f_0  (x,y)$, 
as well as the optimal value function $\varphi$ defined in \eqref{varphi}.
In the next subsection, we start with the derivation of necessary optimality conditions for problem \eqref{OPQVIC-LLVF-K}. Subsequently,  a semismooth Newton method to solve the problem is proposed and tested on the BOLIB \cite{BOLIB2017} library of bilevel optimization problems cast as problem  \eqref{OPQVIC}. Throughout this section, we assume that  $F$, $G$, $f_0$, and $g_0$ are twice continuously differentiable.

\subsection{Necessary and sufficient optimality conditions}\label{First order necessary optimality conditions}
We start here by recalling that a standard constraint qualification, such as the MFCQ, cannot hold for problem \eqref{OPQVIC-LLVF-K} \cite{DempeZemkohoGenMFCQ,YeZhuOptCondForBilevel1995}. Hence, to derive its necessary optimality conditions, we instead consider the partial penalization
\begin{equation}\label{OPQVIC-LLVF-PEN}\tag{OPQVI-R$_\lambda$}
  \underset{z}\min~F(z) + \lambda \left(f(z) - \varphi(z)\right) \;\; \mbox{ s.t. } \;\; G(z)\leq 0,  \;\;  g(z)\leq 0,
\end{equation}
 which removes the value function constraint $f(z) - \varphi(z)\leq 0$, responsible for the failure of constraint qualifications, from the feasible set of { {the problem}}. In problem \eqref{OPQVIC-LLVF-PEN}, $\lambda \in (0, \infty)$ corresponds to the penalization parameter. A close connection can be established between problems \eqref{OPQVIC-LLVF-K} and \eqref{OPQVIC-LLVF-PEN}, based on the so-called  \emph{partial calmness} concept \cite{YeZhuOptCondForBilevel1995}. But before we present the corresponding result, note that problem \eqref{OPQVIC-LLVF-K} will be said to be partially calm at one of its feasible points $\bar z$ if there is $\lambda \in (0, \infty)$ and a neighborhood $U$ of $(0, \bar z)$ such that 
\begin{equation}\label{partial calmness}
\begin{array}{l}
  F(z)-F(\bar z )+\lambda |\varsigma|\geq 0,\;\; \forall (\varsigma, z)\in U: \, G(z)\leq 0,\,\,  g(z)\leq 0, \;  f(z)- \varphi(z)+\varsigma=0.
\end{array}
\end{equation}
\begin{theorem}[{ {\cite{YeZhuOptCondForBilevel1995}}}]\label{equivalen}
 Let  $ \bar z $ be locally optimal for problem \eqref{OPQVIC-LLVF-K}. Then the problem is partially calm at $\bar z$ if and only if there exists $\lambda \in (0, \infty)$ such that this point is also locally optimal for problem \eqref{OPQVIC-LLVF-PEN}.
\end{theorem}

For a moment now, we are going to discuss sufficient conditions to ensure that problem \eqref{OPQVIC-LLVF-K} is partially calm. The first result of this series is a bit like a counterpart of Theorem \ref{I}.
\begin{theorem}\label{OptimalityConditions}Let $\bar z$ be a local optimal solution of problem \eqref{OPQVIC-LLVF-K}. The problem is partially calm at $\bar z$ if a set-valued mapping $\bar\Psi$, obtained by replacing $\Omega$ in \eqref{PsiCalm} by the following set, is calm at $(0, \bar z)$:
\begin{equation*}
\bar\Omega := \left\{z:=(x,y)\in \mathbb{R}^n \times \mathbb{R}^m\left|\, g(z)\leq 0, \;\, G(z)\leq 0\right.\right\}.
\end{equation*}
\end{theorem}

\begin{proof}From the definition of the calmness of {  {$\bar\Psi$}} at $(0, \bar z)$, as given in \eqref{Definition Calmness}, there exists a neighborhood $U$ of $(0, \bar z)$ and a number $\ell>0$ such that
\begin{equation}\label{CalmEneq}
d(z, \bar{\Psi}(0))\leq \ell |\varsigma| \;\; \forall (\varsigma, z)\in U, \;\; z\in {  {\bar\Psi(\varsigma)}}.
\end{equation}
As $\bar z$ is a local optimal solution of problem  \eqref{OPQVIC-LLVF-K}, {  {$\bar\Psi(0)\neq \emptyset$}}. Hence, consider a point {  {$z^*\in \bar\Psi(0)$}} such that { {{{$d(z, \bar\Psi(0))=\|z-z^*\|$}} for some $z\in \bar{\Psi}(\varsigma)$ with $(\varsigma, z) \in U$}}. Furthermore, also based on the fact that $\bar z$ is a local optimal solution of problem \eqref{OPQVIC-LLVF-K}, assume without loss of generality that $U$ is small enough such that $F(\bar z) \leq F(z^*)$, considering the fact that {  {$\bar\Psi(0)$}} coincides with the feasible set of \eqref{OPQVIC-LLVF-K}. Denoting by $\lambda>0$ a {  {Lipschitz}} constant of $F$ near $\bar z$, it follows that for all  $(\varsigma, z)\in U$ and {  {$z\in \bar\Psi(\varsigma)$}},
\[
F(\bar z) - F(z) \leq F(z^*) - F(z) \leq \lambda \|z-z^*\|\leq \ell\lambda |\varsigma|
\]
by considering \eqref{CalmEneq}. The proof ends by comparing this condition with \eqref{partial calmness}. \qed
\end{proof}

Clearly, the framework of sufficient conditions ensuring that the set-valued mapping $\Psi$ \eqref{PsiCalm} is calm can
straightforwardly be extended to $\bar\Psi$, with an adjustment similar to the one done in Theorem \ref{OptimalityConditions}. Without unnecessarily repeating corresponding results from the previous section here, we can summarize them in the following extended counterpart of Figure \ref{Figure1}:

\begin{figure*}[htp]
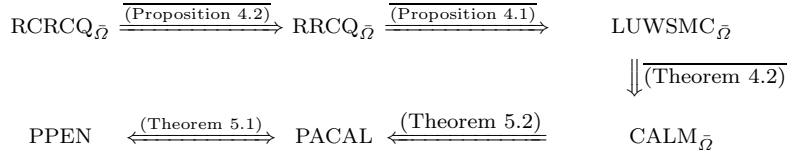

\begin{center}
\[
\begin{array}{ccccc}
   \text{RCRCQ}_{\bar\Omega} & \xRightarrow{ \overline{\text{(Proposition \ref{vertix})}} } & \text{RRCQ}_{\bar\Omega} & \xRightarrow{ \overline{\text{(Proposition \ref{RRCQimpliesLUWSMC})}} } & \text{LUWSMC}_{\bar\Omega}\\[2ex]
   & \qquad \qquad  &  & \qquad \qquad  & \qquad \quad \Big\Downarrow  \overline{\text{\footnotesize{\text{(Theorem \ref{sharp-one})}}}} \\[2ex]
  \text{PPEN} & \xLeftrightarrow{ \text{(Theorem \ref{equivalen})} }  & \text{PACAL} & \xLeftarrow{\text{\footnotesize{(Theorem \ref{OptimalityConditions})}} } & \text{CALM}_{\bar \Omega}
\end{array}
\]
${}$\\[-3ex]
\caption{\emph{The index $\bar\Omega$ here distinguishes the corresponding relationships with the ones in Figure \ref{Figure1}, which are instead based on $\Psi$ \eqref{PsiCalm} associated to $\Omega$. The lines over  Proposition \ref{vertix}, Proposition \ref{RRCQimpliesLUWSMC}, and Theorem \ref{sharp-one} specify that we are referring here to the versions of these results for the set-valued mapping $\bar\Psi$ introduced in Theorem \ref{OptimalityConditions}. As for PACAL, it represents the partial calmness condition \eqref{partial calmness}, while PPEN is used for the existence of $\lambda \in (0, \infty)$ such that a given point $\bar z$ is locally optimal for the partially penalized problem \eqref{OPQVIC-LLVF-PEN}. Finally, note that  $\text{RCRCQ}_{\bar\Omega}$ represents the framework in Proposition \ref{vertix} with $\bar \Omega$ being the counterpart of $\Omega$ given in Theorem \ref{OptimalityConditions}.}}\label{Figure2}
\end{center}
\end{figure*}
Now, based on the partial penalization of problem \eqref{OPQVIC-LLVF-K}, we are going to establish the necessary optimality conditions that will be at the center of the analysis in this section. To proceed, we need the \textit{upper-level regularity} condition, which will be said to hold for problem \eqref{OPQVIC-LLVF-K} at a point $\bar z$  if there exists a vector  $d\in \mathbb{R}^{n+m}$   such that 
\begin{equation}\label{MFCQ-UL}
   \nabla G_i (\bar z)^\top d < 0 \mbox{ for } i\in I^1 \; \mbox{ and } \; \nabla g_j (\bar z)^\top d < 0 \mbox{ for } j\in I^2
\end{equation}
with $I^1$ and $I^2$ defined in \eqref{I1}.  
Obviously, \eqref{MFCQ-UL} corresponds to the MFCQ for the joint upper- and lower-level constraints of problem \eqref{OPQVIC-LLVF-K}. In the same vein,  the \textit{lower-level regularity} condition will be said to hold at $(\bar x, \bar y, \bar\varsigma)$  if the qualification condition \eqref{LL regularity} is satisfied at this point.  
\begin{theorem}\label{KN stationarity partial calm 0}  Let $z:=(x, y)$ be a local optimal solution of problem \eqref{OPQVIC-LLVF-K}.
Suppose that $S_0$ is inner semicontinuous at the point $(x, y, y)$ with $y\in S_0(z)$, where the lower-level regularity condition \eqref{LL regularity} also holds. Furthermore, suppose that problem \eqref{OPQVIC-LLVF-K} is partially calm at the point $z$, where we also assume that the upper-level regularity \eqref{MFCQ-UL} is satisfied. 
 Then there exist $\lambda \in ]0, \; \infty[$, $u\in \mathbb{R}^p$, and $(v, w)\in \mathbb{R}^{2q}$ such that
 \begin{subequations}
 \begin{eqnarray}\label{optimality conditions r0}
 \nabla F(z) + \nabla  G(z)^{\top}u +\nabla  g(z)^{\top}v+\lambda \nabla  f(z) - \lambda \nabla_z \ell(z, y, w) =0,\label{KS-r1}\\
 f_0(z)+\nabla_3 g_0(x, y, y)^{\top}w=0,\label{KS-r2}\\
u\geq 0, \; G(z)\leq 0,\; u^{\top}G(z)= 0,\label{VS-r3}\\
v\geq 0, \; g(z)\leq 0,\; v^{\top} g(z)= 0, \label{VS-r4}\\
w\geq 0, \; g_0(z, y)\leq 0,\; w^{\top} g_0(z, y)= 0,\label{KS-r3}
\end{eqnarray}
\end{subequations}
where $\nabla_z \ell(z, y, w)$ represents the gradient of the function $(x, y, z, w) \mapsto \ell(x, y, z, w)$ w.r.t. $(x, y)$ at $(x, y, y, w)$ and similarly, $\nabla_3 g_0(x, y, y)$ represents the Jacobian of the function $(x, y, \varsigma) \mapsto g_0(x, y, \varsigma)$ w.r.t. $\varsigma$ at $(x, y, y)$.
\end{theorem}
\begin{proof}The proof technique is well-known, see, e.g., \cite{DempeDuttaMordukhovichNewNece,DempeZemkohoGenMFCQ,YeZhuOptCondForBilevel1995}, but we establish the result here for the sake of completeness. Start by observing that since $S_0$ \eqref{S-v} is inner semicontinuous at $(z, y)$ and the MFCQ \eqref{LL regularity} holds at $(z, y)$, then $\varphi$ \eqref{varphi} is then locally Lipschitz continuous around $z$; { {see, e.g., \cite{MordukhovichNamPhanVarAnalMargBlP}}}. Problem \eqref{OPQVIC-LLVF-K} is therefore a locally Lipschitz continuous optimization problem.

 Next, note that under the partial calmness condition, it follows from Theorem \ref{equivalen} that we can find a number $\lambda>0$ such that $z$ is a local optimal solution of problem \eqref{OPQVIC-LLVF-PEN}. Applying the Lagrange multiplier rule for locally Lipschitz optimization  on the latter problem, it follows that, as the upper-level regularity condition \eqref{MFCQ-UL} holds at $z$, there exist $u\in \mathbb{R}^p$ and $v\in \mathbb{R}^q$ such that \eqref{VS-r3} and \eqref{VS-r4} hold together with
\begin{equation}\label{KKTCondMan}
\nabla F(z) + \nabla G(z)^\top u + \nabla g(z)^\top v + \lambda \nabla f(z) \in \lambda\bar\partial \varphi(z)
\end{equation}
given that $\bar\partial(-\varphi)(z)= -\bar\partial\varphi(z)$, since $\varphi$  is  locally Lipschitz continuous around $z$. Furthermore, under the { {inner semicontinuity of $S_0$ at $(z, y)$}} and fulfillment of the upper-level regularity condition \eqref{MFCQ-UL} at the same point,
\begin{equation}\label{PartSub}
\bar\partial \varphi(z)\;\, \subseteq \;\,  \left\{\nabla_z\ell( z,y,w) \left|\;\, w\in \Lambda(x, y, y) \right.\right\}
\end{equation}
{ {with  $\Lambda(x, y, y)$ defined in \eqref{Lambda(x,y)}}}. Subsequently, combining \eqref{KKTCondMan} and \eqref{PartSub}, we get \eqref{KS-r1}--\eqref{KS-r2} and \eqref{KS-r3}. \qed
\end{proof}
A few comments on this result are in order. First, the lower- and upper-level regularity conditions as MFCQ can easily be verified. As for the inner semicontinuity assumption of the set-valued mapping $S_0$ \eqref{S-v}, see relevant discussion at the end of the previous section, including Proposition \ref{InnerSemS}. The only thing to add here with regards to relaxing the inner semicontinuity condition by using the inner semicompactness is that it would lead to a more complicated system of optimality conditions (see relevant discussion in \cite{DempeZemkohoGenMFCQ}), which would create new challenges for the Newton scheme to be introduced in the next section. As for the partial calmness condition, see Figure \ref{Figure2} for a suitable framework ensuring that it holds. Finally, looking at the optimality conditions \eqref{KS-r1}--\eqref{KS-r3},
it is clear, similarly to the discussion from the previous section, that they do not involve second order
 information as it is the case in \cite{MOC07}, for example. Additionally, as shown in the series of papers \cite{FischerZemkohoZhouSemismooth2019,FliegeTinZemkoho2021,ZemkohoZhou2021,TinZemkohoLevenberg2022}, optimality resulting from the value function approach have the potential to have a better numerical behaviour compared to the ones from GE-type reformulation. This was one of the motivations to study the version of the semismooth Newton method introduced in the next subsection.
 
\subsection{Semismooth Newton-type method}\label{Newton method for the auxiliary equation}
From here on, we propose and study a method to solve the system \eqref{KS-r1}--\eqref{KS-r3}. To proceed, the first thing to address is whether these conditions fully represent the optimality conditions of problem \eqref{OPQVIC} based on reformulation \eqref{OPQVIC-LLVF-K}. This comes as the underlying inner semicontinuity assumption in Theorem \ref{KN stationarity partial calm 0} presumes that $y\in S_0(z)$. It is important to point out that this inclusion is redundant if we impose the convexity of $g_0$ w.r.t. $\varsigma$ in Theorem \ref{KN stationarity partial calm 0}, as it will guarantee that the combination of \eqref{KS-r2} and \eqref{KS-r3} implies that $y\in S_0(z)$.

Secondly, we can easily check that the system \eqref{KS-r1}--\eqref{KS-r3} is a \emph{nonsquare} system of equations. Using a trick introduced in \cite{FischerZemkohoZhouSemismooth2019}, we can get a \emph{square} system by adding the new variable $\xi\in \mathbb{R}^m$. To proceed, consider the Lagrangian-type function $L^{\lambda}$ defined by 
\begin{eqnarray*}\label{Upper Lagrangian}
\begin{array}{l}
L^{\lambda}( z,\varsigma,u,v,w):=\mathcal{L}^\lambda(z,u,v) - \lambda \ell (z,\varsigma,w),
\end{array}
\end{eqnarray*}
where $\ell (z,\varsigma,w)$  is given in (\ref{ell}) and $\mathcal{L}^\lambda(z,u,v)$ is  defined by
\begin{eqnarray*}\label{upper-level lagrangian}
\mathcal{L}^{\lambda}(z,u,v)& := &  F(z) +u^\top G(z) + v^\top g(z) + \lambda  f(z){  {.}}
\end{eqnarray*}
{ {Based on these tools, we can easily check that conditions \eqref{KS-r1}--\eqref{KS-r3} can be rewritten as
\begin{eqnarray}\label{Phi-lambda}
\Phi^{\lambda}(\zeta ):=\left[
\begin{array}{l}
  \nabla L^{\lambda}(z,\varsigma,u,v,w)\\
 \sqrt{G(z)^2 + u^2} + G(z) - u\\
 \sqrt{g(z)^2 + v^2} + g(z) - v\\
 \sqrt{g_0(z, \varsigma)^2 + w^2} + g_0(z, \varsigma) -w
\end{array}
\right]=0,
\end{eqnarray}
while setting $\varsigma=y$. Note that here, the square root is understood vector-wise, $\zeta : =(z,\varsigma,u,v,w)$, and   $\nabla L^{\lambda}$ represents the gradient of the function  $L^\lambda$ w.r.t. $(z, \varsigma)$.
\eqref{Phi-lambda} is actually a relaxation of the system \eqref{KS-r1}-\eqref{KS-r3}, as the condition $y=\varsigma$ is not explicitly included. It is easy to show that if $S_0(z)=\{y\}$, then $\varsigma=y$ and this system reduces to \eqref{KS-r1}-\eqref{KS-r3}. Most importantly, solving \eqref{Phi-lambda} is much easier as a square system, and numerical experiments show solving this system enables the proposed method to possess excellent performance in terms of finding local/global solutions. Note that without replacing $y$ by $\varsigma$, the system \eqref{KS-r1}-\eqref{KS-r3} is overdetermined; for interested readers, Gauss-Newton and Levenberg–Marquardt-type methods for similar systems resulting from  optimistic bilevel programs can be found in \cite{FliegeTinZemkoho2021,JolaosoMehlitzZemkoho2024,TinZemkohoLevenberg2022}. Such methods are out of the scope of this paper.

It might be important to emphasize that $\lambda$ is treated here as a parameter to our system; namely, it is the exact penalty parameter appearing in problem \eqref{OPQVIC-LLVF-PEN}. If the partial calmness condition is replaced in Theorem \ref{KN stationarity partial calm 0} by the calmness of the indicated set-valued mapping $\bar{\Psi}$, then a version of the system  \eqref{Phi-lambda} could be obtained with $\lambda$ being a variable. Interested readers are refers to the paper \cite{JolaosoMehlitzZemkoho2024} for an analysis in a similar context for optimistic bilevel optimization problems.}}

Obviously, \eqref{Phi-lambda} is a system of $n+2m+p+2q$ equations with  $n+2m+p+2q$ variables in $\zeta $. This, therefore, allows for a natural extension of standard versions of the semismooth Newton method (see, e.g., \cite{DeLuca1996,FischerASpecial1992,QiJiangSemismooth1997,QiSun1999,QiSunANonsmoothVersion1993}) to the bilevel optimization setting. In order to take full advantage of the structure of the function $\Phi^{\lambda}$ \eqref{Phi-lambda}, we will use the following globalized version of the semismooth Newton method developed by De Luca et al. \cite{DeLuca1996}.
Recall that there are various other classes of functions generally known as NCP (nonlinear complementarity problem) functions that have been used in the literature to reformulate complementarity conditions into equations; see \cite{Galantai2012} and references therein for an extended list and related properties.

{ {Note that in the sequel, we need the merit function 
\[
\Psi^{\lambda}(\zeta) := \frac{1}{2}\|\Phi^\lambda(\zeta)\|^2,
\]
which is continuously differentiable, thanks to the square root-based reformulation of the complementarity conditions \eqref{VS-r3}--\eqref{KS-r3} used in the system of equations \eqref{Phi-lambda}. This continuous differentiability property is crucial for the globalization of the method presented here.}}

\begin{algorithm}
\caption{Semi-smooth Newton Method for Bilevel Optimization}
\label{algorithm-SSNBO}
\begin{algorithmic}
 \STATE \textbf{Step 0}: Choose $\lambda,  \epsilon, M>0,\rho\in(0,1), \sigma\in(0,1/2), t>2$, $\zeta ^o:=( z^o, \varsigma^o,  u^o, v^o, w^o)$ and set $k:=0$.
 \STATE \textbf{Step 1}: If $\|\Phi^{\lambda} (\zeta ^k)\|<\epsilon$ or $k\leq M$, then stop.
 \STATE \textbf{Step 2}: Choose $W^k\in \partial_B \Phi^{\lambda}(\zeta ^k)$ and find the solution $d^k$ of the system
  $$W^k d=-\Phi^{\lambda}(\zeta ^k).$$
\hspace{1.25cm}If the above system is not solvable or if the condition
  \[
   \nabla\Psi^{\lambda}(\zeta ^k)^\top  d^k  \leq -\rho \Vert d^k\Vert^t 
  \]
\hspace{1.25cm}is not satisfied, set $d^k=-\nabla\Psi^{\lambda}(\zeta ^k)$.
 \STATE \textbf{Step 3}: Find the smallest nonnegative integer $s_k$ such that
$$\Psi^{\lambda}(\zeta^k + \rho^{s_k}d^k)   \leq   \Psi^{\lambda}(\zeta^k) + 2\sigma\rho^{s_k}\nabla\Psi^{\lambda}(\zeta ^k)^\top    d^k.$$
\hspace{1.25cm}Then set $ \alpha_k :=\rho^{s_k}, \zeta ^{k+1}:=\zeta ^k + \alpha_k d^k$, $k:=k+1$ and go to \textbf{Step 1}.
\end{algorithmic}
\end{algorithm}

It is important to clarify that the only difference between Algorithm \ref{algorithm-SSNBO} and the original one in \cite{DeLuca1996} is that in Step 0, we also have to provide the penalization parameter $\lambda$ in \eqref{OPQVIC-LLVF-PEN}. Also recall that in Step 2, $\partial_B \Phi^{\lambda}$ denotes the B-subdifferential \eqref{B-Subdifferential}. Obviously, equation $W^k d =-\Phi^{\lambda}(\zeta^k)$ has a solution if  $W^k$ is a nonsingular  matrix. The latter holds in particular if the function $\Phi^{\lambda}$ is  {BD-regular}. $\Phi^{\lambda}$ is said to be BD-regular at a point $\zeta$ if each element of  $\partial_B \Phi^{\lambda}(\zeta)$ is nonsingular.
Using this property, the convergence of Algorithm \ref{algorithm-SSNBO} can be established as follows (see \cite{DeLuca1996}), { {where, based on the definitions in Subsection \ref{Nonsmooth functions and normal cones}, problem \eqref{OPQVIC} is said to be SC$^1$  (resp. LC$^2$) if the functions $F$, $G_i$ with $i=1, \ldots, p$, $f_0$, and $g_{0j}$ with $j=1, \ldots, q$ are all SC$^1$  (resp. LC$^2$). Also note that  \eqref{OPQVIC}  being SC$^1$ (resp. LC$^2$) guarantees that  $\Phi^\lambda$ is semismooth (resp. strongly semismooth); cf. \cite{QiJiangSemismooth1997}.}}
\begin{theorem}\label{convergence result}  Let problem \eqref{OPQVIC} be SC$^1$ and  $\bar \zeta:=(\bar  z, \bar \varsigma, \bar u, \bar v, \bar w)$ an accumulation point of a sequence generated by Algorithm \ref{algorithm-SSNBO}  for some parameter $\lambda > 0$. Then $\bar \zeta$ is a stationary point of the problem of minimizing $\Psi^\lambda$, i.e., $\nabla \Psi^\lambda(\bar \zeta)=0$. If $\bar \zeta$ solves $\Phi^\lambda (\zeta)=0$ and the function $\Phi^\lambda$ is BD-regular at $\bar \zeta$, then the algorithm converges to  $\bar \zeta$ superlinearly. { {If problem \eqref{OPQVIC} is LC$^2$, then under the same assumptions, the convergence is quadratic.}}
\end{theorem}
 Results closely related to Theorem \ref{convergence result} are developed in \cite{FischerASpecial1992,QiJiangSemismooth1997} and many other references therein.
Observe that we have imposed BD-regularity in this theorem. We can replace it with the stronger  {CD-regularity}, which refers to the non-singularity of all matrices in  $\partial \Phi^{\lambda}(\bar \zeta)$. In the next section, we focus our attention on the derivation of conditions ensuring that CD-regularity holds.

\subsubsection{CD-regularity}\label{CD-regularity}
To provide sufficient conditions guaranteeing that CD-regularity holds for $\Phi^\lambda$ \eqref{Phi-lambda}, we first construct an upper estimate of the generalized Jacobian of $\Phi^\lambda$.

\begin{theorem}\label{Jacobian Phi}Let the functions $F$, $G$, $f_0$, and $g_0$ be twice continuously differentiable at the point $\bar\zeta:=(\bar z,\bar \varsigma,\bar u,\bar v,\bar w)$.  If $\lambda>0$, then $\Phi^\lambda$ is semismooth at $\bar \zeta$ and any $W^\lambda\in \partial \Phi^\lambda (\bar \zeta)$ can take the form
{\small{
\begin{equation}\label{Wlambda}
   W^\lambda =\left[
\begin{array}{ccccc}
\nabla^2_{zz}\mathcal{L}^\lambda(\bar \zeta)-\lambda \nabla^2_{z z}  \ell(\bar \zeta)    & -\lambda \nabla^2_{z \varsigma}  \ell(\bar \zeta)^{\top}    & \nabla  G(\bar z)^{\top} & \nabla  g(\bar z)^{\top} & -\lambda \nabla_{1} g_0(\bar z,\bar \varsigma)^{\top} \\
-\lambda \nabla^2_{z \varsigma}  \ell(\bar \zeta) &  -\lambda \nabla^2_{\varsigma\varsigma}  \ell(\bar \zeta) & O & O & -\lambda \nabla_{2} g_0(\bar z,\bar \varsigma)^{\top}  \\
  \Lambda_1\nabla  G(\bar z) &   O   & \Gamma_1 & O & O \\
  \Lambda_2 \nabla  g(\bar z) &   O   & O & \Gamma_2 & O \\
  \Lambda_3 \nabla_{1}g_0(\bar z,\bar \varsigma) &   \Lambda_3 \nabla_{2} g_0(\bar z,\bar \varsigma)   & O & O & \Gamma_3
\end{array}
\right],
\end{equation}
}}
where
 $\Lambda_i :={\rm diag} (a^i) $ and $\Gamma_i :={\rm diag}(b^i)$, $i=1, 2, 3$, are such that
 \begin{equation}\label{ab definition}
    (a^i_j,b^i_j)\left\{\begin{array}{lll}
                  =&(0,-1) & \mbox{ if } \;j\in  \eta^i, \\
                  =&(1,0) & \mbox{ if } \;j\in \nu^i, \\
                  \in& \{(\alpha, \beta): \; (\alpha-1)^2 + (\beta+1)^2\leq 1\} & \mbox{ if }\; j\in \theta^i,
                \end{array}
\right.
 \end{equation}
where $\eta^i$, $\nu^i$, and $\theta^i$ with $i=1, 2, 3$ are defined in \eqref{multiplier sets} and \eqref{nu2nu3}.
\end{theorem}

In the next result, we provide conditions ensuring that the function $\Phi^{\lambda}$ is CD-regular. To proceed,  analogously to \eqref{LICQ}, the  LICQ will be said to hold at  $\bar z $ for problem \ref{OPQVIC-LLVF-PEN} if the family of vectors
\begin{equation}\label{LICQ-leader}
\left\{\nabla G_i(\bar z):\; i\in I^1\right\} \cup  \left\{\nabla g_j(\bar z):\; j\in I^2\right\}
\end{equation}
is linearly independent. Furthermore, let us introduce the { {\textit{critical subspace}}}
\[
Q(\bar z,\bar \varsigma ) :=  \left\{{ {d^{12}:=\left(d^1,\,d^2\right)\in \mathbb{R}^{n+m}\times \mathbb{R}^m}}\left|\,
\begin{array}{rll}
\nabla G_j(\bar z)^\top d^1 &=0, &\; j\in \nu^1 \\
\nabla g_j(\bar z)^\top d^1&=0, & \;j\in \nu^2 \\
\nabla g_{0j}(\bar z,\bar \varsigma)^\top d^{12} &=0, & \;j\in \nu^3
\end{array}
\right.\right\},
\]
where the index sets $\nu^i$ for $i=1,\, 2,\, 3$ are defined as in \eqref{multiplier sets}--\eqref{nu2nu3}.   The last condition is {  {related to the strict}} complementarity condition (SCC). A point $(\bar z, \bar \varsigma, \bar w)$ will be said to satisfy the SCC if it holds that
\begin{equation}\label{SCC}
\theta^3 := \theta^h(\bar z, \bar \varsigma, \bar w)=\emptyset.
\end{equation}
Furthermore, for the next result, we use the notation
\[
\Box^2_{z\varsigma}\ell(\zeta):= \left[\begin{array}{cc}
   \nabla^2_{zz}\ell(\zeta)  &  \nabla^2_{\varsigma z}\ell(\zeta)  \\[0.75ex]
  \nabla^2_{z\varsigma}\ell(\zeta)    & \nabla^2_{\varsigma\varsigma}\ell(\zeta) 
\end{array} \right].
\]
\begin{theorem}\label{SOSSC-Theorem 1-1} Assume that problem \eqref{OPQVIC} is SC$^1$ and let the point $\bar \zeta:=(\bar z, \bar \varsigma, \bar u, \bar v, \bar w)$ satisfy the optimality conditions \eqref{KS-r1}--\eqref{KS-r3} 
 for some $\lambda >0$. Suppose that the LICQ \eqref{LICQ-leader} and  \eqref{LICQ} hold at $ \bar z$ and $(\bar z, \bar \varsigma)$, respectively. If  additionally, for all $d^{12}\in Q(\bar z, \bar \varsigma)\setminus \{0\}$, we have
  \begin{equation}\label{SOSSC}
  \begin{array}{l}
 (d^1)^{\top} \nabla^2_{zz}\mathcal{L}^{\lambda}(\bar \zeta) d^1 > \lambda (d^{12})^{\top} \Box^2_{z\varsigma}\ell(\bar \zeta) d^{12}
  \end{array}
  \end{equation}
and the SCC \eqref{SCC} is also satisfied at $(\bar z, \bar \varsigma, \bar w)$, then  $\Phi^{\lambda}$ is CD-regular at $\bar \zeta$.
\end{theorem}
\begin{proof}
Let $W^\lambda\in \partial \Phi^{\lambda}(\bar \zeta)$. Then from \eqref{Wlambda}--\eqref{ab definition}, it follows that for any $d:=(d^1, d^2, d^3, d^4, d^5)$ with $d^1\in \mathbb{R}^{n+m}$, $d^2\in \mathbb{R}^m$, $d^3\in \mathbb{R}^p$, $d^4\in \mathbb{R}^q$ and $d^5\in \mathbb{R}^q$ such that $W^\lambda d =0$, we have
\begin{subequations}
\begin{eqnarray}
\left[\nabla^2_{zz}  \mathcal{L}^\lambda (\bar\zeta)-\lambda \nabla^2_{z z}  \ell(\bar \zeta) \right]d^1 - \lambda \nabla^2_{z\varsigma} \ell (\bar\zeta)^\top d^2 + \;\nabla G(\bar z)^\top d^3 + \nabla  g(\bar z)^\top d^4  -\lambda \nabla_1 g_0(\bar z, \bar \varsigma)^\top d^5=0,\label{p11}\\
  -\;\lambda\nabla^2_{z \varsigma}\ell (\bar \zeta)d^1 - \lambda \nabla^2_{\varsigma\varsigma} \ell (\bar \zeta)d^2 -\lambda \nabla_2 g_0(\bar z, \bar \varsigma)^\top d^5=0,\label{p13}\\
  \forall j=1, \ldots, p, \; a^1_j \nabla G_j(\bar z)^\top d^1 + b^1_j d^3_j =0,\label{p14}\\
    \forall j=1, \ldots, q, \; a^2_j \nabla g_j(\bar z)^\top d^1 + b^2_j d^4_j =0,\label{p15}\\
     \forall j=1, \ldots, q, \; a^3_j \nabla h_j(\bar z, \bar \varsigma)^\top d^{12} + b^3_j d^5_j =0.\label{p16}
\end{eqnarray}
\end{subequations}
{ {Multiplying \eqref{p11} (resp. \eqref{p13}) from the left by $(d^1)^\top$ (resp. $(d^2)^\top$) and adding the resulting equations together,}} 
\begin{equation}\label{Quad-term}
\begin{array}{l}
    (d^1)^\top\nabla^2_{zz}  \mathcal{L}^\lambda (\bar\zeta)d^1 - 
\lambda (d^{12})^\top\Box^2_{z\varsigma}\ell (\bar\zeta)d^{12} \\[2ex]
  \qquad \qquad  +    (d^3)^\top \nabla G(\bar z)  d^1 +  (d^4)^\top \nabla  g(\bar z)  d^1 -\lambda  (d^5)^\top \nabla  g_0(\bar z, \bar \varsigma)  d^{12}=0.
\end{array}
\end{equation}

Recall that $p$ and $q$ represent the number of constraint functions $G$ and $g$. For $i=1, 2, 3$, let $p_1:=p$, $p_2:=q$ and $p_3:=q$. Then define 9 index sets by
\begin{eqnarray}
\label{P1i}P^i_1&:=&\{j\in\{1 \ldots, p_i\}~|~a^i_j>0,~b^i_j<0\}, \\
\label{P2i}P^i_2&:=&\{j\in\{1 \ldots, p_i\}~|~a^i_j=0,~b^i_j=-1\}\hspace{5mm}\supseteq \;\;\;\eta^i,\\
\label{P3-nu} P^i_3&:=&\{j\in\{1 \ldots, p_i\}~|~a^i_j=1,~b^i_j=0\}\hspace{8mm}\supseteq \;\;\;\nu^i,
\end{eqnarray}
where the relations `$\supseteq$' can be verified easily from (\ref{ab definition}). For example, for any $j\in \eta^i$, we have $a^i_j=0$, $b^i_j=-1$ from \eqref{ab definition}. Thus, $j\in P^i_2$. It follows from \eqref{p14}--\eqref{p16} that for $j\in P^1_2$, $j\in P^2_2$, and $j\in P^3_2$,
\begin{equation}\label{yes1}
    d^3_j=0, \; d^4_j=0,\, \mbox{ and }d^5_j=0,
\end{equation}
respectively due to $a^i_j=0,~b^i_j=-1$. As for $j\in P^1_3$, $j\in P^2_3$ and $j\in P^3_3$, we respectively get
\begin{equation}\label{yes2}
    \nabla G_j (\bar z)^\top d^{1}=0,\; \nabla g_j (\bar z)^\top d^{1}=0, \, \mbox{ and }\, \nabla h_j (\bar z,\bar \varsigma)^\top d^{12}=0.
\end{equation}
due to $a^i_j=1,~b^i_j=0$. Now observe that under the SCC \eqref{SCC},  $\theta^3=\emptyset$. Hence, from the corresponding counterpart of \eqref{ab definition}, it follows that $P^3_1=\emptyset$. We can further check that  for $j\in P^1_1$ and $j\in P^2_1$,
\begin{equation}\label{yes3}
    \nabla G_j (\bar z)^\top d^{1 }=c^1_j d^3_j \,\mbox{ and }\,  \nabla g_j (\bar z)^\top d^{1 }=c^2_j d^4_j,
\end{equation}
where $c^1_j:=- b^1_j/a^1_j $ and $c^2_j:=- b^2_j/a^2_j $, respectively. Since $P^3_1=\emptyset$, $d^5_j=0$ for $j\in P^3_{2}$ by \eqref{yes1} and $\nabla h_j(\bar z,\bar \varsigma)^{\top} d^{12}=0$ for $j\in P^3_{3}$ by \eqref{yes2}, it follows that
\begin{equation}\label{SCS-applied}
 (d^5)^\top\nabla h(\bar z,\bar \varsigma) d^{12} =  \sum_{j\in P^3_{2}}d^5_j\nabla h_j(\bar z,\bar \varsigma)^{\top} d^{12} + \sum_{j\in P^3_{3}}d^5_j\nabla g_{0j}(\bar z,\bar \varsigma)^{\top} d^{12}=0
\end{equation}
Inserting \eqref{SCS-applied} into \eqref{Quad-term} while taking into account \eqref{yes1}--\eqref{yes3}, we get
 \begin{eqnarray*}
0&=&\underset{=:\Delta}{\underbrace{(d^{1})^\top \nabla^2_{zz} \mathcal{L}^\lambda (\bar \zeta) d^{1} - \lambda (d^{12})^\top \Box^2_{z\varsigma}\ell (\bar\zeta) d^{12}}}
  + \sum_{j\in P^1_1}c^1_j (d^3_j)^2 + \sum_{j\in P^2_1}c^2_j (d^4_j)^2 \geq \Delta,
\end{eqnarray*}
where the inequality results from  $c^1_j >0$ for $j\in P^1_1$ and $c^2_j >0$ for $j\in P^2_1$. Again, as \eqref{yes2} and  $\nu^i \subseteq P^i_3$ hold for $i=1, 2, 3$ (from \eqref{P3-nu}), it results that $d^{12}\in Q(\bar z, \bar \varsigma)$. If $d^{12}\neq0$, then it holds that $\Delta>0$, considering \eqref{SOSSC}. This contradicts the above inequality. Therefore, we have $d^{12}=0$. Hence, $d^{1}=0$, $d^{2}=0$, and $d^3_j=0$ for $j\in P^1_1$ and $d^4_j=0$ for $j\in P^2_1$.
Substituting these  into \eqref{p11}--\eqref{p13}, it follows from \eqref{yes1} and $P^3_1=\emptyset$ that
 \begin{eqnarray}
 \sum_{j\in P^1_3}d^3_j\nabla G_j(\bar z) + \sum_{j\in P^2_3}d^4_j\nabla  g_j(\bar z) + \sum_{j\in P^3_3}(-\lambda d^5_j)\nabla_1 g_{0j}(\bar z,\bar \varsigma)=0,\label{p111}\\
  \sum_{j\in P^3_3}d^5_j\nabla_2 g_{0j}(\bar z,\bar \varsigma)=0.\label{p133}
\end{eqnarray}
Since the LICQ \eqref{LICQ} is satisfied at $(\bar z, \bar \varsigma)$ and $P^3_3 \subseteq I^3$ holds, it follows from \eqref{p133} that  $d^5_j=0$ for $j\in P^3_{3}$. Therefore $\sum_{j\in P^1_3}d^3_j\nabla G_j(\bar z) + \sum_{j\in P^2_3}d^4_j\nabla  g_j(\bar z)=0$. This together with the LICQ \eqref{LICQ-leader}  at $ \bar z $, and  $P^i_3 \subseteq I^i$ for $i=1, 2$, we have $d^3_j=0$ for $j\in P^1_3$ and $d^4_j=0$ for $j\in P^2_3$. Hence, $d=0$. \qed
\end{proof}

\subsubsection{Numerical experiments}\label{Numerical experiments}
{\bf a) Implementation details.} We code Algorithm \ref{algorithm-SSNBO} in   {MATLAB (R2023b) on a \textsf{ThinkPad} laptop of 32GB memory and \textsf{Intel(R) Core(TM) i9-13900H   2.60 GHz}}. Before presenting the the results, we first outline the implementation details. Recall that $\lambda$ is a parameter for problem \eqref{OPQVIC-LLVF-PEN}. Hence, one of the main difficulties in implementing Algorithm \ref{algorithm-SSNBO} is to be able to suitably select values of this parameter that can generate stationary points that could be optimal for problem \eqref{OPQVIC-LLVF-K}. It therefore goes without saying that we cannot expect the same value of $\lambda$ to perform well for all the examples under consideration, as solutions obtained depend on this parameter.   {Therefore this is a shortcoming of the proposed algorithm, as in general, $\lambda$ needs be to be carefully chosen in advance. In the sequel, we take} $\lambda$ from a small set of constants. Precisely, we select from the set $\tilde{\Lambda}:=\{3^{-2},3^{-1},3^0,3^{1},3^2\}$ for all the examples used in our experiments if there are no extra explanations.

As starting point, we choose  $z^0:=(x^0, \, y^0)$ that is feasible for $G$ and  $g_0$ (or $g$), i.e., $G(z^0)\leq 0$ and $g(z^0)\leq 0$, respectively. If it is not easy to find such a point, we just set $x^0=0$ and $y^0=0$. After initializing $(x^0,\, y^0)$, we choose $\varsigma^o=y^0$, $u^0= \left(\left|G_1(z^0)\right|, \ldots, \left|G_p(z^0)\right|\right)^\top$, $v^0= \left(\left|g_1(z^0)\right|, \ldots, \left|g_q(z^0)\right|\right)^\top$, and $w^0=v^0$. The choices of $\varsigma^0$, $u^0$, $v^0$, and $w^0$ are not exclusive, as different initial points would render different solutions.

For the input parameters, we set $\epsilon=10^{-6}$, $\beta=10^{-8}$, $t=2.1$, $\rho=0.5$, $M=1000$ and $\sigma=10^{-4}$.   We choose  $s^k$ from $\{-1,0,1,2,\cdots\},$ which means that the starting step length of Newton direction is $\rho^{-1}=2$. This contributes to making our algorithm render a desirable average performance over all examples.    {For the stopping criteria}, apart from setting $\|\Phi^{\lambda}(\zeta^k)\|<\epsilon$ or $k\leq M$, we also terminate the algorithm if   $\|\Phi^{\lambda}(\zeta^k)\|$ does not change significantly within 100 steps, e.g., the variance of the vector $\left(\|\Phi^{\lambda}(\zeta^{k-100})\|,\|\Phi^{\lambda}(\zeta^{k-99})\|,\cdots, \|\Phi^{\lambda}(\zeta^k)\|\right)$ is less than $10^{-6}$. The latter termination criterion was used for a very small portion of testing examples; resulting in a sequence $\{\zeta^k\}$ converging  (to a point $\zeta^*$), but which still did not satisfy $\|\Phi^{\lambda}(\zeta^*)\|<\epsilon$. A potential reason to explain such a phenomenon is that the corresponding selection of $\lambda$ is not appropriate.   {For instance, for example {\tt	NieWangYe2017Ex34}, our method obtained the global optimal solution for $\lambda\in\{3^1,3^0,3^{-1},3^{-2}\}$ but with  $\|\Phi^{\lambda}(\zeta^*)\|>\epsilon$.}
Therefore, to avoid such a phenomenon, which might result in unnecessary calculations, it is reasonable to terminate the algorithm if   $\|\Phi^{\lambda}(\zeta^k)\|$ does not significantly vary within 100 steps. Note that in order to pick elements from the generalized Jacobian of $\Phi^{\lambda}$, cf. Step 2 of Algorithm \ref{algorithm-SSNBO}, we extend the a corresponding algorithm from \cite{DeLuca1996} to our context.\\

\noindent{\bf b) Solving two OPQVI examples from the literature.} Before going to our main set of examples that we apply Algorithm \ref{algorithm-SSNBO} to, we first test it on two examples from \cite{MOC07}, which is one of the publications that motivated the work in this paper. The constraints in the first example describe a Nash game of two players \cite{H91} and in the second one, they describe an oligopolistic market equilibrium  \cite{FK92,MSS82}.
\begin{figure}[htp]
\centering
  \includegraphics[width=.75\linewidth]{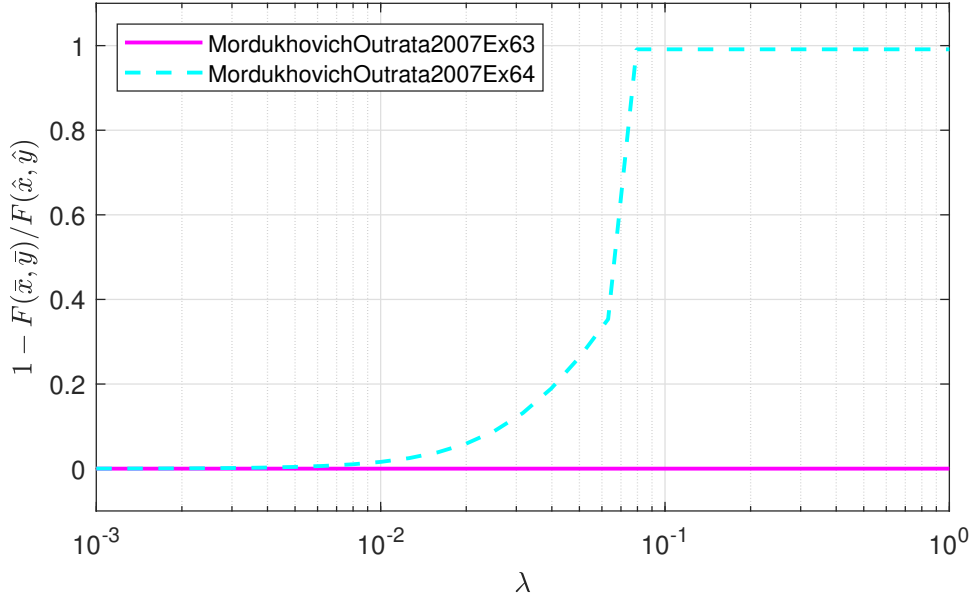}\vspace{-2mm}
\caption{The relative difference   {between the objective value obtained in our experiments
and the best known objective function value from the literature}, i.e., $1-F(\bar{x},\bar{y})/F(\hat{x},\hat y)$, under a selection of  $\lambda\in[0.001,\,1]$.}\vspace{-3mm}
\label{impact-lambda}
\end{figure}
\begin{example}\label{EX1} (\texttt{MordukhovichOutrata2007Ex63} \cite{H91,MOC07}) Consider a version of problem \eqref{OPQVIC} with  
\begin{eqnarray*}
  F(x,y)&:=&x-3y_1- 11y_2/3+ (y_1-9)^2/2,\\
  G(x,y)&:= &\left[
  \begin{array}{l}
  -1-x\\
  -1+x
  \end{array}
  \right],\\
  f_0(x,y)&:= &\left[
  \begin{array}{l}
  -34+2y_1+8y_2/3\\
  - 97/4+5y_1/4+2y_2
  \end{array}
  \right],\\
   g_0(x,y,\varsigma)&:= &\left[
  \begin{array}{l}
   y_1+\varsigma_2-15-x\\
   y_2+\varsigma_1-15-x
  \end{array}
  \right]{  {.}}
\end{eqnarray*}
  {A local optimal solution of this problem  given in \cite{MOC07}  is $\bar{x}=0$, $\bar{y}=(9,6)^\top$ with $F(\bar{x},\, \bar{y})=-49$. One can verify that this solution is optimal.} \qed
\end{example}

\begin{example}\label{EX2} (\texttt{MordukhovichOutrata2007Ex64} \cite{FK92,MOC07,MSS82}) Consider the version of  \eqref{OPQVIC} described by
\begin{eqnarray*}
  F(x,y)&:=&(\max\{x-135,0\})^2+0.6(y_1-34)^2+0.6(y_2-16)^2,\\
  f_0(x,y)&:= &\left[
  \begin{array}{r}
  -76+2y_1+ y_2\\
  -72+y_1+2y_2
  \end{array}
  \right],\\
   g_0(x,y,\varsigma)&:= &\left[
  \begin{array}{r}
   -0.333x+\varsigma_1+\varsigma_2\\
   -\varsigma_1\\
   -\varsigma_2
  \end{array}
  \right],
\end{eqnarray*}
without any constraint of the form $G(x, y)\leq 0$. An   {approximate optimal solution to this problem given in \cite{MOC07} is $\bar{x}'=135.15$, $\bar{y}'=(30.95,14.1)^\top$ with $F(\bar{x}' ,\bar{y}')=7.77$.  A more accurate solution than  $\bar{x}',\bar{y}'$  is $\bar{x}=135.488$, $\bar{y}=(31.56,13.56)^\top$ with $F(\bar{x},\,\bar{y})=7.39$.} \qed 
\end{example}

\begin{figure}[H]
\centering
  \includegraphics[width=.425\linewidth]{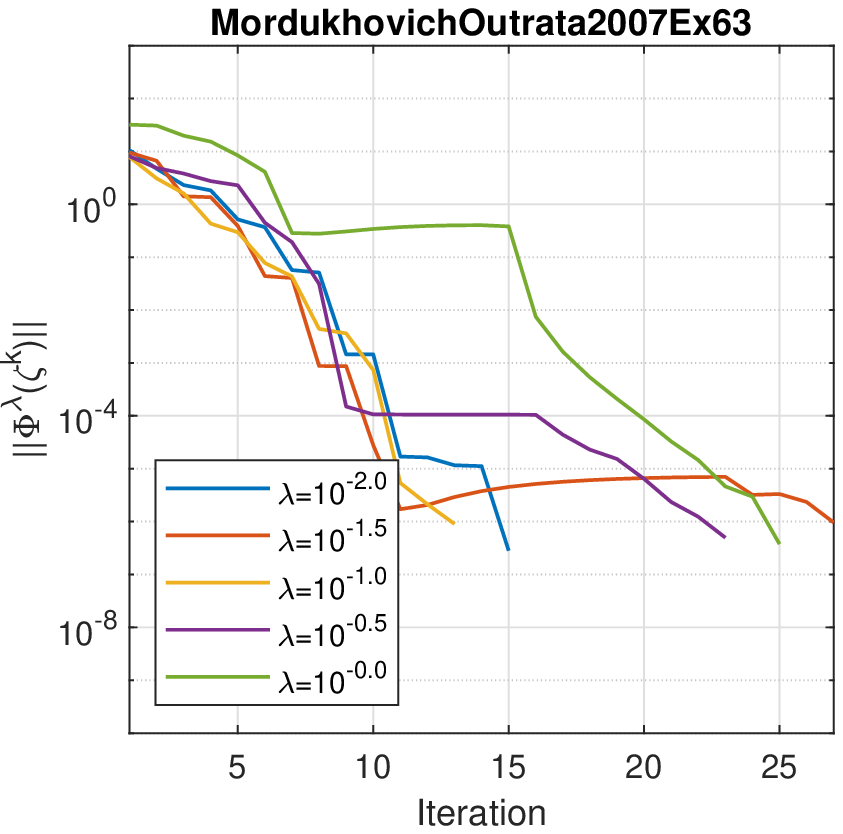}\qquad
    \includegraphics[width=.425\linewidth]{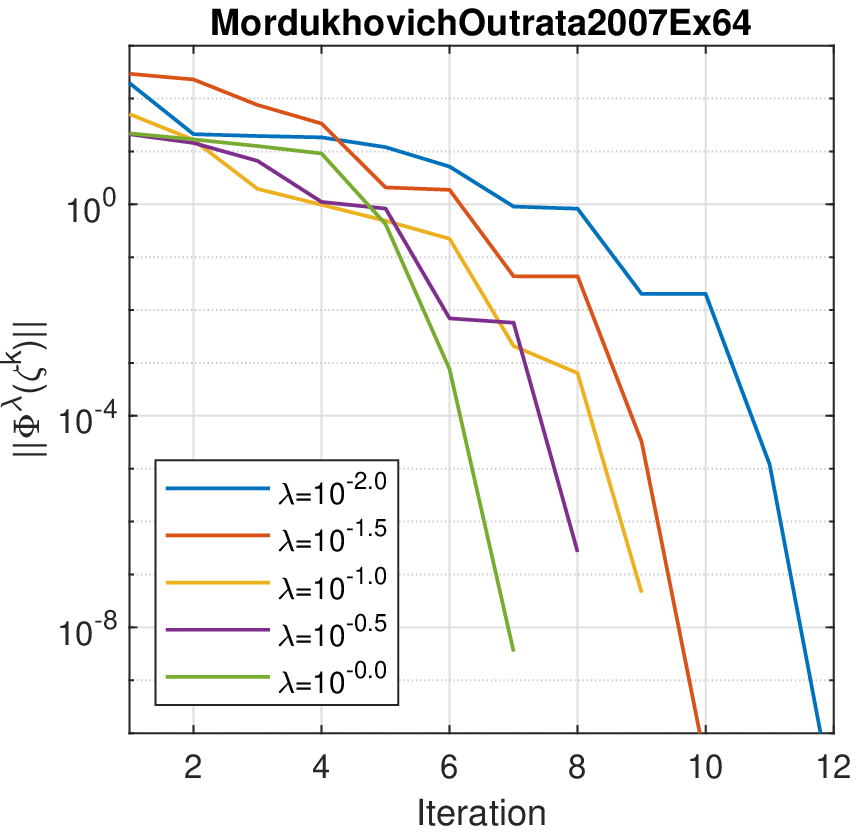}\\
      \includegraphics[width=.425\linewidth]{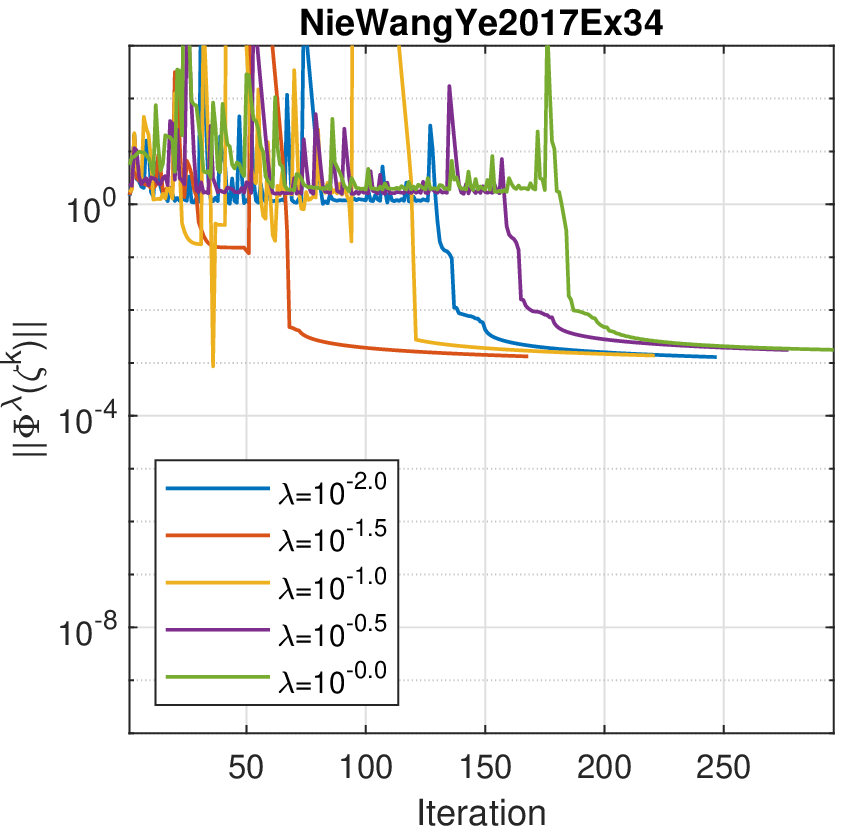}\qquad
    \includegraphics[width=.425\linewidth]{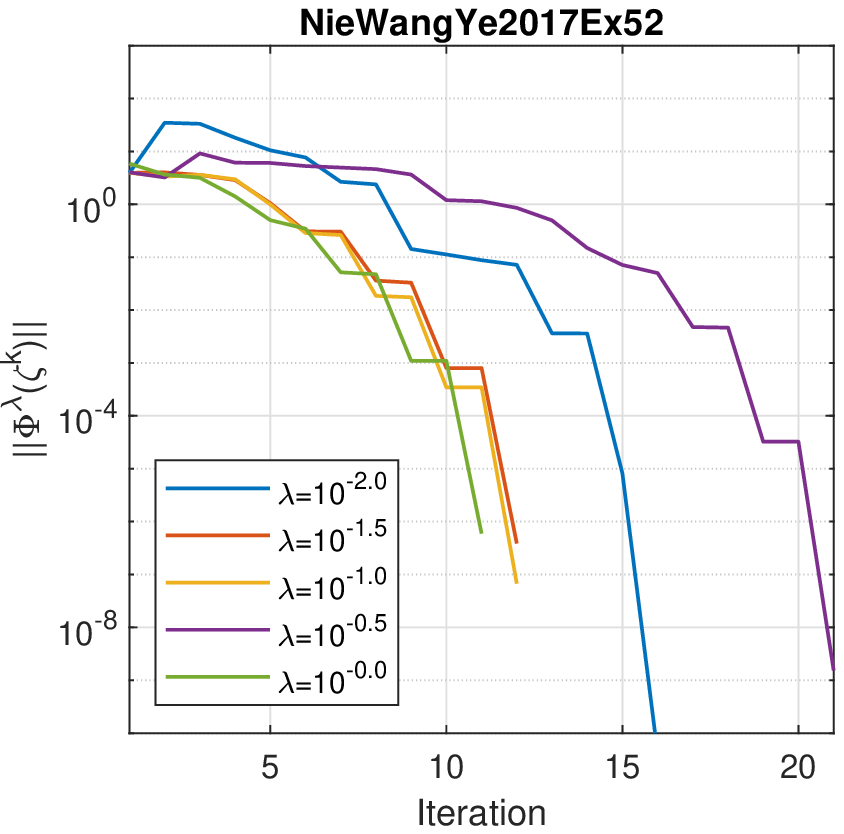}\\
\caption{  {Variation of the values of $\|\Phi^{\lambda}(\zeta^k)\|$ at each iteration.}}
\label{stopping}
\end{figure}

  {The starting point used for these two examples is $x^0=0$, $y^0=(0,0)^\top$.} 
To see the impact of $\lambda$ on the performance of Algorithm \ref{algorithm-SSNBO}, we solve Examples \ref{EX1} and \ref{EX2} by choosing $\lambda\in\{0.001,0.002,\cdots,1\}$. The results are reported in Figure \ref{impact-lambda}, where $(\hat{x}, \hat y)$ denotes the solution obtained from Algorithm \ref{algorithm-SSNBO}. For Example  \ref{EX1}, Algorithm \ref{algorithm-SSNBO}    {yields the optimal solutions for all $\lambda$ due to  due to $1-F(\bar{x} ,\bar{y})/F(\hat{x},\hat y)=0$. However, one can verify that the LICQ does not hold at $(\hat{x},\hat y)$ which fails the assumptions of Theorem 4.6}.

 For Example  \ref{EX2},   {$F(\hat{x},\hat{y})$ increases and tends to be stable eventually with the rising of $\lambda$. We observe that $F(\hat{x} ,\hat{y})$ is identical to $F(\bar{x} ,\bar{y})$ when  $\lambda<0.0018$. 
By fixing $\lambda=0.001$, the algorithm obtained solution $\hat{x}=\bar{x},\;\hat y=\bar{y}, \;\hat \varsigma=(0,45.112)^\top,\; \hat v=(2.93,0,0)^\top, \;\hat w=(13.32,13.99,0)^\top$. In the sequel, we check that the assumptions of Theorem \ref{SOSSC-Theorem 1-1} hold at this solution. Let $\hat z:=(\hat x^\top, \hat y^\top)^\top$. It is easy to obtain that
\begin{eqnarray*}
\begin{array}{lllll}
I^1=\emptyset, &  I^2=\{1\}, & I^3=\{1,2\},\;\nu^1=\emptyset, &\nu^2=\{1\}, & \nu^3=\{1,2\}, \; \theta^3=\emptyset.
\end{array}
\end{eqnarray*}
Therefore,  both $\nabla g_1(\hat z)=(-0.333,1,1)^\top$  and $\{\nabla_\varsigma g_j(\hat z, \hat \varsigma),\; j\in I^3\}=\{(1,1)^\top,\; (-1,0)^\top\}$ are linearly independent. Namely, the LICQ \eqref{LICQ-leader} and  \eqref{LICQ} hold at $ \hat z$ and $(\hat z, \hat \varsigma)$, respectively.  Moreover, the SCC \eqref{SCC} is also satisfied at $(\hat z, \hat \varsigma, \hat w)$ due to $\theta^3=\emptyset$. The remaining is to verify condition \eqref{SOSSC}. Based on the definition, we have
  \begin{equation*} 
  \begin{array}{lll}
Q(\hat z,\hat\varsigma ) =  \left\{(d^1, \,d^2)\Bigg|
\begin{array}{rll}
-0.333 d^1_1+d^1_2+d^1_3 &=0 \\
-0.333 d^1_1 + d^2_1 + d^2_2&=0\\
-d^2_1&=0
\end{array}
\right\}=  \left\{(d^1_1, d^1_2, d^1_3, 0, d^2_2)^\top \left|
\begin{array}{rll}
0.333d^1_1&=&d^2_2\\
 d^1_2+d^1_3 & =& d^2_2  
\end{array}
\right.\right\}
\end{array}
  \end{equation*}
  and for any $d^{12}\in Q(\hat z, \hat \varsigma)\setminus \{0\}$,  
  \begin{eqnarray*}
  (d^1)^\top \nabla^2\mathcal{L}^{\lambda}(\hat \zeta) d^1 - \lambda (d^{12})^{\top} \nabla^2\ell(\hat \zeta) d^{12} 
&=&
\left[ \begin{array}{c}
 d^1_1\\
  d^1_2\\ d^1_3
 \end{array}\right]^\top
\left[ \begin{array}{ccc}
 2&0&0\\
 0&1.2+4\lambda&4\lambda\\
 0&4\lambda&1.2+4\lambda\\
 \end{array}\right]
 \left[ \begin{array}{c}
 d^1_1\\
  d^1_2\\ d^1_3
 \end{array}\right]- \lambda \left[ \begin{array}{c}
 d^1_1\\
  d^1_2\\ 
  d^1_3\\
  0\\
  d^2_2
 \end{array}\right]^\top
\left[ \begin{array}{ccccc}
 0&0&0&0&0\\
 0&0&0&2&1\\
 0&0&0&1&2\\
 0&2&1&0&0\\
 0&1&2&0&0\\
 \end{array}\right]
 \left[ \begin{array}{c}
 d^1_1\\
  d^1_2\\
    d^1_3\\
  0\\
  d^2_2
 \end{array}\right]\\[1ex]
 &=& 2(d^1_1)^2  + (1.2+\lambda)(d^1_2)^2+  (1.2-\lambda)(d^1_3)^2+\lambda (d^2_2)^2>0
 \end{eqnarray*}
  due to $\lambda=0.001$. Overall, the assumptions of Theorem \ref{SOSSC-Theorem 1-1} hold at this solution. From the above calculation, we can observe that too large  $\lambda$ may degrade the performance of the algorithm since the larger  $\lambda$ results in the higher possibility of violating condition \eqref{SOSSC}. }

  {We also examined the local convergence behavior of Algorithm \ref{algorithm-SSNBO}. As depicted in Figure \ref{stopping}, we observe a rapid local convergence towards the end of a run for examples {\tt {MordukhovichOutrata2007Ex64}} and {\tt NieWangYe2017Ex52}, but not for examples {\tt {MordukhovichOutrata2007Ex63}} and {\tt NieWangYe2017Ex34}. One possible explanation for this phenomenon is that the CD-regularity condition holds at the final solution for the former two examples but not for the latter two examples. In other words, if the CD-regularity condition holds at the final point for a testing problem, then the local fast convergence might be observed for this problem.}\\

\noindent{\bf c) Solving bilevel optimization examples.} {We now focus our attention to the set of problems from the Bilevel Optimization LIBrary (BOLIB) \cite{BOLIB2017}. Recall that a bilevel   {optimization problem} can take the form
\begin{equation*}\label{Bilevel}
\underset{x, y}\min~F(x,y)\;\,  \mbox{ s.t. } \;\, G(x,y)\leq 0, \;\, y\in \arg\underset{y}\min~\left\{\hat{f}(x,y)\left|\; \hat{g}(x,y)\leq 0 \right.\right\},
\end{equation*}
where $\hat{f} :\mathbb{R}^n\times \mathbb{R}^m \rightarrow \mathbb{R}$ and $\hat{g}:\mathbb{R}^n\times \mathbb{R}^m \rightarrow \mathbb{R}^q$. This problem can be cast as \eqref{OPQVIC} with $f_0(x,y):= \nabla_2 \hat{f}(x,y)$, provided that the lower-level problem is convex; i.e., the functions $\hat{f}(x, .)$ and $\hat{g}_{j}$, $j=1, \ldots, q$ are convex for all $x\in \mathbb{R}^n$.    {The only reason why we consider the BOLIB test set for our experiments here is because it provides access to a large test base. The convexity assumption here has no bearing in our analysis, as it is just needed to put problem  \eqref{Bilevel} in the format of \eqref{OPQVIC}. The convexity assumption is not satisfied for all our examples and is not relevant in the  analysis of the performance of the algorithm. More precisely, the BOLIB library has 124 examples, where 72 of them (i.e., with No 1-72 in Table \ref{numerical-results})  have convex lower-level problems, while the lower-level problem is nonconvex for the remaining 52 examples (i.e., with No 73-124 in Table \ref{numerical-results}). We will solve them under choices of  $\lambda\in\tilde{\Lambda}$.
Results of all examples are listed in Table \ref{numerical-results}, where $\hat{F}:=F(\hat{x},\hat{y})$ represents our optimal value of the leader's objective function $F$, $\bar{F}:=F(\bar{x},\bar{y})$ 
 represents the known/true optimal value of $F$, and $\|\hat{\Phi}\|:=\|\Phi^\lambda(\hat{\zeta})\|$ with $\hat{\zeta}$ being the final solution obtained by Algorithm \ref{algorithm-SSNBO}.  The complete results and related technical details can be found in the supplementary material in \cite{ZhouZemkohoDetailed2018}, where it can be seen that for each example we run Algorithm  \ref{algorithm-SSNBO}  under two starting points. However, to save space, results listed in Table \ref{numerical-results} are generated under one of these two starting points only. Note that most of the starting points are  taken from the existing literature, as the BOLIB library is essentially made of a collection from published paper or technical reports.}

In the third column of Table \ref{numerical-results}, NaN stands for examples with unknown solutions. Notice that from the original reference, some examples are associated with a parameter that should be provided by the user; i.e., for example, \texttt{CalamaiVicente1994a}, \texttt{HenrionSurowiec2011}, and \texttt{IshizukaAiyoshi1992a}.  The first one is associated with $\rho\geq 1$, which separates the problem into 4 cases: (i) $\rho=1$, (ii) $1<\rho<2$, (iii) $\rho=2$, and (iv) $\rho>2$. The results presented in Table \ref{numerical-results} correspond to case (i). For the other three cases, our method still produced the true global optimal solutions.  Problem \texttt{HenrionSurowiec2011} has a unique global optimal solution $-0.5c(1,1)$, where $c$ is the parameter to be provided by the user. The results presented in Table \ref{numerical-results} are for $c=0$. We also tested our method when $c=\pm1$, and obtained the unique optimal solutions as well. Problem \texttt{IshizukaAiyoshi1992a} contains a parameter $M>1$, and the results presented in Table \ref{numerical-results} correspond to $M=4$.

 \begin{sidewaystable}
 \vspace{17cm}
  \centering \caption{Upper-level objective function values at the solution for  different selections of the penalty parameter $\lambda\in\tilde{\Lambda}$.}\label{numerical-results}  
{\small  \begin{tabular}{p{0.6cm}p{4.8cm}rrrrrrrrrrrrrrr}\hline
No	&	Problems	&		&	\multicolumn{2}{c}{$\lambda=3^{-2}$}			&	&	\multicolumn{2}{c}{$\lambda=3^{-1}$}			&	&	\multicolumn{2}{c}{$\lambda=3^0$}			&	&	\multicolumn{2}{c}{$\lambda=3^1$}			&	&	\multicolumn{2}{c}{$\lambda=3^2$}			\\\cline{4-5}\cline{7-8}\cline{10-11}\cline{13-14}\cline{16-17}
	&		&	$\bar{F}$	&	$\hat{F}$	&	$\|\hat{\Phi}\|$	&	&	$\hat{F}$	&	$\|\hat{\Phi}\|$	&	&	$\hat{F}$	&	$\|\hat{\Phi}\|$	&	&	$\hat{F}$	&	$\|\hat{\Phi}\|$	&	&	$\hat{F}$	&	$\|\hat{\Phi}\|$	\\\hline
1	&	\texttt{AiyoshiShimizu1984Ex2}	&	5.00 	&	5.00 	&	2.08E-07	&	&	5.00 	&	6.60E-12	&	&	10.00 	&	1.21E-08	&	&	0.00 	&	6.35E-12	&	&	11.57 	&	7.39E-02	\\\rowcolor[gray]{.9}
2	&	\texttt{AllendeStill2013}	&	-1.50 	&	-8.49 	&	2.41E-08	&	&	-8.15 	&	8.95E-08	&	&	-7.14 	&	7.30E-07	&	&	-6.25 	&	7.29E-07	&	&	-6.25 	&	6.80E-14	\\
3	&	\texttt{Bard1988Ex1}	&	17.00 	&	2.00 	&	2.53E-08	&	&	7.68 	&	1.97E-09	&	&	17.00 	&	1.14E-09	&	&	30.59 	&	5.81E-07	&	&	48.42 	&	6.50E-09	\\\rowcolor[gray]{.9}
4	&	\texttt{Bard1988Ex2}	&	-6600.00 	&	-6600.00 	&	3.46E-11	&	&	-6600.00 	&	2.68E-07	&	&	-6600.00 	&	1.15E-08	&	&	-5852.78 	&	1.69E-07	&	&	-5772.97 	&	4.58E-11	\\
5	&	\texttt{Bard1991Ex1}	&	2.00 	&	2.00 	&	2.06E-08	&	&	2.00 	&	1.72E-09	&	&	2.00 	&	6.76E-10	&	&	2.00 	&	5.77E-12	&	&	2.00 	&	3.71E-07	\\\rowcolor[gray]{.9}
6	&	\texttt{BardBook1998}	&	NaN	&	0.00 	&	3.06E-12	&	&	0.00 	&	4.98E-09	&	&	0.00 	&	1.66E-12	&	&	11.11 	&	1.72E-07	&	&	11.11 	&	2.91E-11	\\
7	&	\texttt{CalamaiVicente1994a}	&	0.00 	&	0.00 	&	6.21E-07	&	&	0.00 	&	1.40E-11	&	&	0.00 	&	8.22E-07	&	&	0.00 	&	2.76E-08	&	&	0.00 	&	7.13E-06	\\\rowcolor[gray]{.9}
8	&	\texttt{CalamaiVicente1994b}	&	NaN	&	0.06 	&	9.77E-07	&	&	0.31 	&	9.97E-08	&	&	0.50 	&	6.38E-09	&	&	1.25 	&	9.74E-07	&	&	1.44 	&	9.05E-12	\\
9	&	\texttt{CalamaiVicente1994c}	&	0.31 	&	0.06 	&	5.22E-08	&	&	0.31 	&	7.30E-07	&	&	0.50 	&	4.62E-07	&	&	0.50 	&	6.36E-07	&	&	1.40 	&	2.48E-07	\\\rowcolor[gray]{.9}
10	&	\texttt{ClarkWesterberg1990a}	&	5.00 	&	22.23 	&	1.07E+00	&	&	1.70 	&	3.35E-10	&	&	5.00 	&	1.22E-08	&	&	7.43 	&	1.01E-11	&	&	5.00 	&	1.28E-11	\\
11	&	\texttt{Colson2002BIPA1}	&	250.00 	&	250.00 	&	1.34E-08	&	&	250.00 	&	4.51E-11	&	&	250.00 	&	1.83E-11	&	&	250.00 	&	4.28E-13	&	&	250.00 	&	5.03E-09	\\\rowcolor[gray]{.9}
12	&	\texttt{Colson2002BIPA2}	&	17.00 	&	2.00 	&	3.07E-12	&	&	17.00 	&	6.62E-08	&	&	17.00 	&	5.47E-13	&	&	30.59 	&	1.20E-11	&	&	25.00 	&	1.26E-08	\\
13	&	\texttt{Colson2002BIPA3}	&	2.00 	&	2.00 	&	3.79E-07	&	&	2.00 	&	8.68E-08	&	&	2.00 	&	3.12E-07	&	&	2.00 	&	2.39E-13	&	&	2.00 	&	1.42E-10	\\\rowcolor[gray]{.9}
14	&	\texttt{Colson2002BIPA5}	&	2.75 	&	0.20 	&	2.18E-07	&	&	0.89 	&	1.08E-08	&	&	2.96 	&	8.51E-03	&	&	2.75 	&	7.08E-09	&	&	2.75 	&	3.00E-07	\\
15	&	\texttt{Dempe1992b}	&	31.25 	&	4.39 	&	1.31E-12	&	&	6.69 	&	3.40E-07	&	&	18.45 	&	6.07E-07	&	&	31.25 	&	6.72E-07	&	&	31.25 	&	7.18E-08	\\\rowcolor[gray]{.9}
16	&	\texttt{DempeDutta2012Ex24}	&	0.00 	&	0.00 	&	1.05E+00	&	&	0.00 	&	2.43E-01	&	&	0.00 	&	9.44E-02	&	&	47.94 	&	1.57E-01	&	&	14.05 	&	2.33E-02	\\
17	&	\texttt{DempeDutta2012Ex31}	&	-1.00 	&	-2.15 	&	1.30E-01	&	&	-1.12 	&	3.33E-10	&	&	-0.46 	&	8.57E-07	&	&	-0.16 	&	1.52E-10	&	&	0.00 	&	7.84E-06	\\\rowcolor[gray]{.9}
18	&	\texttt{DempeEtal2012}	&	-1.00 	&	-1.00 	&	2.08E-10	&	&	-1.00 	&	5.31E-12	&	&	-1.00 	&	8.60E-08	&	&	-1.00 	&	2.87E-10	&	&	0.00 	&	3.05E-14	\\
19	&	\texttt{DempeFranke2011Ex41}	&	5.00 	&	-0.88 	&	9.69E-07	&	&	-0.22 	&	1.15E-07	&	&	2.56 	&	1.49E-09	&	&	4.52 	&	2.68E-09	&	&	5.00 	&	8.49E-07	\\\rowcolor[gray]{.9}
20	&	\texttt{DempeFranke2011Ex42}	&	2.13 	&	-0.93 	&	7.12E-07	&	&	-0.62 	&	6.56E-07	&	&	0.61 	&	5.46E-07	&	&	2.13 	&	3.48E-07	&	&	3.13 	&	4.61E-12	\\
21	&	\texttt{DempeFranke2014Ex38}	&	-1.00 	&	-3.00 	&	7.55E-08	&	&	-3.00 	&	1.01E-11	&	&	-2.00 	&	2.93E-09	&	&	-1.00 	&	6.45E-03	&	&	-1.00 	&	2.49E-09	\\\rowcolor[gray]{.9}
22	&	\texttt{DempeLohse2011Ex31a}	&	-6.00 	&	-5.99 	&	8.72E-07	&	&	-5.94 	&	2.32E-10	&	&	-5.50 	&	9.77E-07	&	&	-5.50 	&	9.15E-08	&	&	0.00 	&	4.02E-09	\\
23	&	\texttt{DempeLohse2011Ex31b}	&	-12.00 	&	-12.00 	&	9.71E-07	&	&	-12.00 	&	8.87E-07	&	&	-12.00 	&	7.77E-07	&	&	-0.91 	&	1.11E+00	&	&	1.33 	&	1.93E-01	\\\rowcolor[gray]{.9}
24	&	\texttt{DeSilva1978}	&	-1.00 	&	-1.48 	&	4.09E-09	&	&	-1.28 	&	1.34E-09	&	&	-1.00 	&	6.53E-07	&	&	-1.00 	&	2.25E-09	&	&	-1.00 	&	7.17E-07	\\
25	&	\texttt{EdmundsBard1991}	&	0.00 	&	5.00 	&	6.20E-08	&	&	5.00 	&	4.22E-11	&	&	0.00 	&	5.76E-11	&	&	10.00 	&	6.71E-07	&	&	90.00 	&	2.89E-08	\\\rowcolor[gray]{.9}
26	&	\texttt{FalkLiu1995}	&	-2.25 	&	-3.98 	&	1.18E-11	&	&	-3.72 	&	1.24E-11	&	&	-2.25 	&	4.14E-07	&	&	-2.25 	&	5.05E-12	&	&	-2.25 	&	8.92E-09	\\
27	&	\texttt{GumusFloudas2001Ex1}	&	2250.00 	&	1652.20 	&	1.92E-13	&	&	1863.78 	&	6.65E-13	&	&	2018.86 	&	1.65E-08	&	&	2250.00 	&	3.32E-12	&	&	2250.00 	&	1.47E-09	\\\rowcolor[gray]{.9}
28	&	\texttt{GumusFloudas2001Ex4}	&	9.00 	&	1.00 	&	3.58E-08	&	&	4.84 	&	2.98E-09	&	&	9.00 	&	2.22E-08	&	&	9.00 	&	8.18E-07	&	&	9.00 	&	9.22E-12	\\
29	&	\texttt{HatzEtal2013}	&	0.00 	&	0.00 	&	6.47E-13	&	&	0.00 	&	9.97E-11	&	&	0.00 	&	6.83E-11	&	&	7414.39 	&	2.77E-01	&	&	359.99 	&	1.08E-01	\\\rowcolor[gray]{.9}
30	&	\texttt{HendersonQuandt1958}	&	-3266.67 	&	-3266.67 	&	4.93E-07	&	&	-3266.67 	&	4.06E-13	&	&	-3266.67 	&	4.82E-11	&	&	-3266.67 	&	5.03E-11	&	&	-3266.67 	&	4.43E-11	\\
31	&	\texttt{HenrionSurowiec2011}	&	0.00 	&	0.00 	&	0.00E+00	&	&	0.00 	&	5.23E-17	&	&	0.00 	&	1.57E-16	&	&	0.00 	&	0.00E+00	&	&	0.00 	&	0.00E+00	\\\rowcolor[gray]{.9}
32	&	\texttt{IshizukaAiyoshi1992a}	&	0.00 	&	0.00 	&	2.89E-05	&	&	0.00 	&	8.96E-07	&	&	0.00 	&	1.69E-05	&	&	0.00 	&	9.25E-07	&	&	0.00 	&	1.48E-09	\\
33	&	\texttt{LamparielloSagratella2017Ex31}	&	1.00 	&	1.00 	&	7.04E-07	&	&	1.00 	&	8.38E-12	&	&	1.00 	&	6.20E-07	&	&	1.00 	&	6.20E-07	&	&	1.00 	&	6.20E-07	\\\rowcolor[gray]{.9}
34	&	\texttt{LamparielloSagratella2017Ex32}	&	0.50 	&	0.50 	&	3.14E-16	&	&	0.50 	&	0.00E+00	&	&	0.50 	&	0.00E+00	&	&	0.50 	&	6.28E-16	&	&	0.50 	&	6.28E-16	\\
35	&	\texttt{LamparielloSagratella2017Ex33}	&	0.50 	&	0.50 	&	9.13E-10	&	&	0.50 	&	1.13E-07	&	&	0.50 	&	4.90E-07	&	&	0.50 	&	5.94E-08	&	&	0.50 	&	7.39E-08	\\\rowcolor[gray]{.9}
36	&	\texttt{LamparielloSagratella2017Ex35}	&	0.80 	&	0.00 	&	6.15E-10	&	&	0.03 	&	1.89E-07	&	&	0.80 	&	2.78E-09	&	&	0.80 	&	2.68E-07	&	&	1.25 	&	3.79E-09	\\
37	&	\texttt{LucchettiEtal1987}	&	0.00 	&	0.50 	&	2.44E-11	&	&	0.00 	&	6.93E-10	&	&	0.00 	&	2.37E-07	&	&	0.78 	&	3.94E-09	&	&	0.00 	&	2.20E-08	\\\rowcolor[gray]{.9}
38	&	\texttt{LuDebSinha2016c}	&	1.12 	&	1.12 	&	3.05E-09	&	&	1.12 	&	3.05E-09	&	&	1.12 	&	3.05E-09	&	&	1.97 	&	1.13E+00	&	&	1.12 	&	6.84E-08	\\
39	&	\texttt{LuDebSinha2016f}	&	NaN	&	-18.54 	&	4.20E-08	&	&	-160.00 	&	6.09E-07	&	&	-160.00 	&	2.55E-09	&	&	-160.00 	&	3.17E-13	&	&	-160.00 	&	3.17E-13	\\\rowcolor[gray]{.9}
40	&	\texttt{MacalHurter1997}	&	81.33 	&	81.33 	&	3.56E-15	&	&	81.33 	&	3.55E-15	&	&	81.33 	&	3.55E-15	&	&	81.33 	&	3.57E-15	&	&	81.33 	&	6.40E-14	\\
41	&	\texttt{MitsosBarton2006Ex38}	&	0.00 	&	0.00 	&	6.83E-07	&	&	0.00 	&	1.54E-07	&	&	0.00 	&	2.79E-07	&	&	0.00 	&	3.50E-11	&	&	0.00 	&	1.58E-10	\\\rowcolor[gray]{.9}
42	&	\texttt{MitsosBarton2006Ex313}	&	-1.00 	&	-2.00 	&	1.41E-12	&	&	-2.00 	&	1.31E-09	&	&	-2.00 	&	1.21E-11	&	&	0.30 	&	7.96E-12	&	&	0.00 	&	3.84E-12	\\
43	&	\texttt{MitsosBarton2006Ex317}	&	0.19 	&	0.00 	&	1.38E-09	&	&	0.00 	&	3.26E-07	&	&	0.00 	&	1.03E-09	&	&	0.19 	&	1.09E-07	&	&	0.19 	&	2.60E-08	\\\rowcolor[gray]{.9}
44	&	\texttt{MorganPatrone2006a}	&	-1.00 	&	-1.50 	&	4.28E-09	&	&	-1.50 	&	5.25E-11	&	&	-1.00 	&	7.08E-12	&	&	-1.00 	&	5.72E-11	&	&	-1.00 	&	3.33E-07	\\
45	&	\texttt{MorganPatrone2006b}	&	-1.25 	&	-1.50 	&	8.47E-10	&	&	-1.50 	&	7.34E-10	&	&	-1.25 	&	2.91E-07	&	&	-0.75 	&	5.07E-07	&	&	0.53 	&	2.90E-09	\\\rowcolor[gray]{.9}
46	&	\texttt{MorganPatrone2006c}	&	-1.00 	&	-3.00 	&	4.13E-11	&	&	-3.00 	&	4.32E-11	&	&	1.00 	&	4.14E-12	&	&	1.00 	&	2.81E-10	&	&	1.00 	&	6.99E-07	\\

		  \end{tabular}}
\end{sidewaystable}

 \begin{sidewaystable} 
 \vspace{17cm}
  \centering  
{\small  \begin{tabular}{p{0.6cm}p{5.5cm}rrrrrrrrrrrrrrr}\hline
 No	&	Problems	&		&	\multicolumn{2}{c}{$\lambda=3^{-2}$}			&	&	\multicolumn{2}{c}{$\lambda=3^{-1}$}			&	&	\multicolumn{2}{c}{$\lambda=3^0$}			&	&	\multicolumn{2}{c}{$\lambda=3^1$}			&	&	\multicolumn{2}{c}{$\lambda=3^2$}			\\\cline{4-5}\cline{7-8}\cline{10-11}\cline{13-14}\cline{16-17}
	&		&	$\bar{F}$	&	$\hat{F}$	&	$\|\hat{\Phi}\|$	&	&	$\hat{F}$	&	$\|\hat{\Phi}\|$	&	&	$\hat{F}$	&	$\|\hat{\Phi}\|$	&	&	$\hat{F}$	&	$\|\hat{\Phi}\|$	&	&	$\hat{F}$	&	$\|\hat{\Phi}\|$	\\\hline
47	&	\texttt{MuuQuy2003Ex1}	&	-2.08 	&	-3.68 	&	1.05E-09	&	&	-3.01 	&	4.57E-09	&	&	-2.24 	&	7.76E-11	&	&	1.30 	&	3.80E-01	&	&	-2.08 	&	2.55E-08	\\\rowcolor[gray]{.9}
48	&	\texttt{MuuQuy2003Ex2}	&	0.64 	&	0.64 	&	4.00E-08	&	&	0.64 	&	1.70E-10	&	&	0.64 	&	5.86E-07	&	&	0.64 	&	1.12E-11	&	&	0.64 	&	2.79E-11	\\
49	&	\texttt{NieWangYe2017Ex54}	&	-0.44 	&	0.00 	&	2.40E-06	&	&	0.00 	&	4.55E-06	&	&	0.00 	&	8.32E-07	&	&	0.00 	&	5.16E-06	&	&	-0.44 	&	1.35E-07	\\\rowcolor[gray]{.9}
50	&	\texttt{Outrata1990Ex1a}	&	-8.92 	&	-11.15 	&	4.60E-08	&	&	-9.01 	&	8.33E-07	&	&	-8.92 	&	5.04E-09	&	&	-8.92 	&	4.35E-08	&	&	-8.92 	&	7.96E-07	\\
51	&	\texttt{Outrata1990Ex1b}	&	-7.56 	&	-11.21 	&	2.30E-07	&	&	-9.25 	&	9.43E-10	&	&	-7.58 	&	7.66E-09	&	&	-7.58 	&	6.00E-07	&	&	-7.58 	&	3.65E-07	\\\rowcolor[gray]{.9}
52	&	\texttt{Outrata1990Ex1c}	&	-12.00 	&	-12.00 	&	9.11E-13	&	&	-12.00 	&	7.10E-07	&	&	-12.00 	&	7.11E-10	&	&	0.00 	&	1.98E-07	&	&	0.00 	&	5.46E-07	\\
53	&	\texttt{Outrata1990Ex1d}	&	-3.60 	&	-6.89 	&	1.57E-08	&	&	-3.88 	&	3.99E-11	&	&	-3.60 	&	3.94E-08	&	&	-0.37 	&	3.97E-03	&	&	-0.38 	&	7.88E-08	\\\rowcolor[gray]{.9}
54	&	\texttt{Outrata1990Ex1e}	&	-3.15 	&	-6.36 	&	5.32E-10	&	&	-4.21 	&	5.66E-07	&	&	-3.92 	&	8.22E-08	&	&	-3.79 	&	1.09E-08	&	&	-3.79 	&	3.93E-07	\\
55	&	\texttt{Outrata1990Ex2a}	&	0.50 	&	0.50 	&	7.64E-10	&	&	0.50 	&	6.31E-08	&	&	0.50 	&	9.41E-07	&	&	0.50 	&	1.34E-10	&	&	0.50 	&	1.46E-08	\\\rowcolor[gray]{.9}
56	&	\texttt{Outrata1990Ex2b}	&	0.50 	&	0.50 	&	8.48E-05	&	&	0.50 	&	2.34E-04	&	&	0.50 	&	3.84E-07	&	&	0.50 	&	9.72E-07	&	&	0.50 	&	7.47E-07	\\
57	&	\texttt{Outrata1990Ex2c}	&	1.86 	&	0.51 	&	5.35E-08	&	&	1.86 	&	1.00E-07	&	&	1.86 	&	3.20E-08	&	&	1.86 	&	5.50E-09	&	&	9.62 	&	2.84E+00	\\\rowcolor[gray]{.9}
58	&	\texttt{Outrata1993Ex32}	&	3.21 	&	2.82 	&	1.54E-07	&	&	3.01 	&	1.10E-07	&	&	3.12 	&	2.07E-07	&	&	8.06 	&	7.48E-08	&	&	3.20 	&	1.52E-08	\\
59	&	\texttt{Outrata1994Ex31}	&	3.21 	&	2.82 	&	6.23E-09	&	&	3.01 	&	7.24E-12	&	&	3.12 	&	2.59E-07	&	&	3.18 	&	4.18E-11	&	&	3.20 	&	2.73E-11	\\\rowcolor[gray]{.9}
60	&	\texttt{OutrataCervinka2009}	&	0.00 	&	0.00 	&	6.63E-07	&	&	0.00 	&	9.98E-07	&	&	0.00 	&	6.35E-07	&	&	102.79 	&	5.47E-01	&	&	0.00 	&	5.94E-07	\\
61	&	\texttt{PaulaviciusAdjiman2017b}	&	-2.00 	&	-2.00 	&	4.93E-11	&	&	-2.00 	&	9.20E-09	&	&	-2.00 	&	3.34E-10	&	&	-2.00 	&	1.07E-12	&	&	0.70 	&	2.05E-10	\\\rowcolor[gray]{.9}
62	&	\texttt{SahinCiric1998Ex2}	&	5.00 	&	0.36 	&	4.55E-12	&	&	1.70 	&	1.91E-13	&	&	4.58 	&	3.77E-09	&	&	5.00 	&	6.92E-08	&	&	13.00 	&	1.96E-10	\\
63	&	\texttt{ShimizuAiyoshi1981Ex1}	&	100.00 	&	66.33 	&	1.52E-10	&	&	82.00 	&	6.34E-08	&	&	92.60 	&	4.71E-11	&	&	97.33 	&	3.75E-10	&	&	99.09 	&	1.25E-09	\\\rowcolor[gray]{.9}
64	&	\texttt{ShimizuAiyoshi1981Ex2}	&	225.00 	&	113.12 	&	6.68E-12	&	&	118.06 	&	4.84E-08	&	&	175.00 	&	3.72E-10	&	&	225.00 	&	6.45E-14	&	&	225.00 	&	1.84E-11	\\
65	&	\texttt{ShimizuEtal1997a}	&	NaN	&	25.00 	&	6.40E-14	&	&	25.00 	&	2.95E-10	&	&	14.04 	&	1.60E-07	&	&	30.59 	&	7.40E-07	&	&	48.42 	&	2.64E-09	\\\rowcolor[gray]{.9}
66	&	\texttt{ShimizuEtal1997b}	&	2250.00 	&	1652.20 	&	4.60E-07	&	&	1863.78 	&	4.76E-12	&	&	2018.86 	&	4.55E-11	&	&	2250.00 	&	4.23E-09	&	&	2329.21 	&	1.05E-12	\\
67	&	\texttt{SinhaMaloDeb2014TP6}	&	-1.21 	&	-1.21 	&	5.77E-10	&	&	-1.21 	&	2.55E-09	&	&	-1.21 	&	8.26E-07	&	&	3.10 	&	7.28E-02	&	&	-1.21 	&	1.50E-11	\\\rowcolor[gray]{.9}
68	&	\texttt{SinhaMaloDeb2014TP8}	&	0.00 	&	0.17 	&	1.53E-09	&	&	4.63 	&	2.24E-10	&	&	25.00 	&	3.33E-12	&	&	73.47 	&	5.93E-08	&	&	0.00 	&	1.62E-04	\\
69	&	\texttt{TuyEtal2007}	&	22.50 	&	0.00 	&	4.62E-07	&	&	0.03 	&	1.11E-08	&	&	0.28 	&	1.41E-09	&	&	22.50 	&	9.24E-07	&	&	24.50 	&	5.11E-08	\\\rowcolor[gray]{.9}
70	&	\texttt{WanWangLv2011}	&	10.62 	&	4.00 	&	5.98E-07	&	&	4.00 	&	3.21E-11	&	&	12.28 	&	3.51E-01	&	&	11.25 	&	9.81E-07	&	&	11.25 	&	1.46E-07	\\
71	&	\texttt{Yezza1996Ex31}	&	1.50 	&	2.00 	&	1.46E-12	&	&	1.00 	&	5.86E-10	&	&	1.50 	&	3.56E-08	&	&	1.50 	&	7.86E-12	&	&	3.78 	&	7.77E-01	\\\rowcolor[gray]{.9}
72	&	\texttt{Yezza1996Ex41}	&	0.50 	&	0.04 	&	3.18E-11	&	&	0.18 	&	1.44E-12	&	&	0.50 	&	8.14E-10	&	&	0.50 	&	1.37E-07	&	&	4.00 	&	8.39E-09	\\
73	&	\texttt{AnEtal2009}	&	2251.55 	&	2251.55 	&	6.37E-07	&	&	2251.55 	&	9.25E-07	&	&	2251.55 	&	1.98E-05	&	&	2250.81 	&	2.00E-01	&	&	2254.33 	&	2.00E+00	\\\rowcolor[gray]{.9}
74	&	\texttt{Bard1988Ex3}	&	-12.68 	&	-10.36 	&	4.44E-09	&	&	-10.36 	&	8.15E-11	&	&	-12.79 	&	4.05E-07	&	&	-12.68 	&	1.45E-07	&	&	-10.36 	&	4.44E-07	\\
75	&	\texttt{CalveteGale1999P1}	&	-29.20 	&	-58.00 	&	8.34E-09	&	&	-58.00 	&	1.44E-10	&	&	-58.00 	&	1.04E-08	&	&	-13.00 	&	1.41E-12	&	&	-29.20 	&	6.38E-13	\\\rowcolor[gray]{.9}
76	&	\texttt{Colson2002BIPA4}	&	88.79 	&	62.52 	&	2.56E-10	&	&	69.90 	&	1.00E-09	&	&	80.08 	&	2.82E-07	&	&	86.94 	&	7.32E-07	&	&	88.79 	&	4.32E-07	\\
77	&	\texttt{Dempe1992a}	&	NaN	&	-2.25 	&	4.04E-09	&	&	-0.81 	&	9.26E-12	&	&	-0.33 	&	3.33E-08	&	&	0.00 	&	2.52E-11	&	&	0.00 	&	8.02E-09	\\\rowcolor[gray]{.9}
78	&	\texttt{FloudasZlobec1998}	&	1.00 	&	0.00 	&	4.22E-07	&	&	0.00 	&	4.64E-07	&	&	-1.00 	&	7.92E-07	&	&	100.00 	&	7.86E-07	&	&	0.93 	&	9.58E-11	\\
79	&	\texttt{GumusFloudas2001Ex3}	&	-29.20 	&	-58.00 	&	5.21E-07	&	&	-29.20 	&	3.87E-11	&	&	-58.00 	&	1.07E-08	&	&	-58.00 	&	7.33E-12	&	&	-58.00 	&	8.26E-12	\\\rowcolor[gray]{.9}
80	&	\texttt{GumusFloudas2001Ex5}	&	0.19 	&	0.19 	&	3.96E-07	&	&	0.19 	&	3.92E-07	&	&	0.19 	&	6.10E-09	&	&	0.19 	&	2.01E-07	&	&	0.19 	&	4.96E-08	\\
81	&	\texttt{KleniatiAdjiman2014Ex3}	&	-1.00 	&	-1.86 	&	2.78E-09	&	&	-1.00 	&	1.14E-12	&	&	-1.00 	&	1.08E-07	&	&	-1.00 	&	2.30E-08	&	&	-0.45 	&	9.22E-07	\\\rowcolor[gray]{.9}
82	&	\texttt{KleniatiAdjiman2014Ex4}	&	-10.00 	&	-5.04 	&	7.75E-12	&	&	-3.73 	&	2.10E-03	&	&	-9.00 	&	1.41E+00	&	&	-7.00 	&	1.65E-08	&	&	-3.39 	&	5.23E-07	\\
83	&	\texttt{LamparielloSagratella2017Ex23}	&	-1.00 	&	-1.00 	&	4.64E-09	&	&	-1.00 	&	7.56E-09	&	&	-1.00 	&	2.73E-10	&	&	-1.00 	&	6.04E-11	&	&	-1.00 	&	3.80E-10	\\\rowcolor[gray]{.9}
84	&	\texttt{LuDebSinha2016a}	&	1.94 	&	1.11 	&	1.90E-08	&	&	1.11 	&	1.90E-08	&	&	1.11 	&	1.90E-08	&	&	1.11 	&	1.90E-08	&	&	1.11 	&	1.90E-08	\\
85	&	\texttt{LuDebSinha2016b}	&	0.00 	&	0.07 	&	1.47E-07	&	&	0.07 	&	1.47E-07	&	&	0.07 	&	1.45E-07	&	&	0.07 	&	1.46E-07	&	&	0.07 	&	1.82E-07	\\\rowcolor[gray]{.9}
86	&	\texttt{LuDebSinha2016d}	&	NaN	&	-192.00 	&	6.30E-08	&	&	-16.27 	&	3.24E-01	&	&	19.91 	&	5.00E-01	&	&	-19.56 	&	2.32E-08	&	&	-155.04 	&	1.02E-09	\\
87	&	\texttt{LuDebSinha2016e}	&	NaN	&	8.35 	&	1.10E-01	&	&	5.44 	&	4.08E-07	&	&	28.65 	&	4.99E-01	&	&	7.85 	&	8.74E-01	&	&	9.75 	&	1.08E+00	\\\rowcolor[gray]{.9}
88	&	\texttt{Mirrlees1999}	&	0.01 	&	0.01 	&	7.91E-14	&	&	0.01 	&	7.42E-08	&	&	0.01 	&	7.31E-08	&	&	0.01 	&	8.93E-14	&	&	0.01 	&	8.65E-14	\\
89	&	\texttt{MitsosBarton2006Ex39}	&	-1.00 	&	-1.00 	&	8.42E-08	&	&	-1.00 	&	5.66E-11	&	&	-0.33 	&	1.33E-07	&	&	-0.06 	&	6.71E-07	&	&	-0.02 	&	8.63E-07	\\\rowcolor[gray]{.9}
90	&	\texttt{MitsosBarton2006Ex310}	&	0.50 	&	0.50 	&	7.22E-09	&	&	0.50 	&	1.17E-11	&	&	0.50 	&	1.78E-10	&	&	-0.25 	&	1.69E-09	&	&	0.50 	&	8.70E-12	\\
91	&	\texttt{MitsosBarton2006Ex311}	&	-0.80 	&	-0.80 	&	9.21E-10	&	&	-0.80 	&	1.76E-09	&	&	0.50 	&	8.35E-08	&	&	0.50 	&	1.79E-09	&	&	-0.50 	&	3.85E-07	\\\rowcolor[gray]{.9}
92	&	\texttt{MitsosBarton2006Ex312}	&	0.00 	&	-1.02 	&	7.25E-10	&	&	-1.01 	&	1.93E-10	&	&	-1.00 	&	3.00E-08	&	&	-1.00 	&	3.13E-08	&	&	-1.00 	&	1.55E-08	\\
		  \end{tabular}}
\end{sidewaystable}

 \begin{sidewaystable} 
 \vspace{17cm}
  \centering 
{\small  \begin{tabular}{p{0.6cm}p{5.5cm}rrrrrrrrrrrrrrr}\hline
 No	&	Problems	&		&	\multicolumn{2}{c}{$\lambda=3^{-2}$}			&	&	\multicolumn{2}{c}{$\lambda=3^{-1}$}			&	&	\multicolumn{2}{c}{$\lambda=3^0$}			&	&	\multicolumn{2}{c}{$\lambda=3^1$}			&	&	\multicolumn{2}{c}{$\lambda=3^2$}			\\\cline{4-5}\cline{7-8}\cline{10-11}\cline{13-14}\cline{16-17}
	&		&	$\bar{F}$	&	$\hat{F}$	&	$\|\hat{\Phi}\|$	&	&	$\hat{F}$	&	$\|\hat{\Phi}\|$	&	&	$\hat{F}$	&	$\|\hat{\Phi}\|$	&	&	$\hat{F}$	&	$\|\hat{\Phi}\|$	&	&	$\hat{F}$	&	$\|\hat{\Phi}\|$	\\\hline
93	&	\texttt{MitsosBarton2006Ex314}	&	0.25 	&	0.00 	&	3.69E-10	&	&	0.03 	&	3.08E-11	&	&	0.06 	&	2.79E-10	&	&	0.06 	&	2.47E-08	&	&	0.06 	&	3.56E-13	\\\rowcolor[gray]{.9}
94	&	\texttt{MitsosBarton2006Ex315}	&	0.00 	&	-2.00 	&	8.54E-12	&	&	-2.00 	&	6.57E-13	&	&	-2.00 	&	5.22E-08	&	&	-2.00 	&	4.47E-11	&	&	0.00 	&	1.10E-11	\\
95	&	\texttt{MitsosBarton2006Ex316}	&	-2.00 	&	-3.00 	&	9.74E-09	&	&	-3.00 	&	1.57E-07	&	&	-2.00 	&	7.32E-07	&	&	-0.25 	&	8.06E-08	&	&	-0.54 	&	2.70E-11	\\\rowcolor[gray]{.9}
96	&	\texttt{MitsosBarton2006Ex318}	&	-0.25 	&	-1.00 	&	6.95E-10	&	&	-1.00 	&	1.45E-09	&	&	-1.00 	&	1.02E-08	&	&	-1.00 	&	1.02E-08	&	&	-1.00 	&	1.02E-08	\\
97	&	\texttt{MitsosBarton2006Ex319}	&	-0.26 	&	-1.20 	&	1.98E-08	&	&	-0.69 	&	1.83E-10	&	&	-0.26 	&	2.23E-10	&	&	-0.26 	&	4.21E-08	&	&	-0.26 	&	2.24E-08	\\\rowcolor[gray]{.9}
98	&	\texttt{MitsosBarton2006Ex320}	&	0.31 	&	0.00 	&	8.28E-08	&	&	0.00 	&	2.39E-11	&	&	0.02 	&	1.98E-09	&	&	0.04 	&	9.82E-08	&	&	0.05 	&	1.56E-08	\\
99	&	\texttt{MitsosBarton2006Ex321}	&	0.21 	&	0.00 	&	8.07E-10	&	&	0.00 	&	3.56E-11	&	&	0.21 	&	1.02E-09	&	&	0.21 	&	1.70E-09	&	&	0.21 	&	2.32E-11	\\\rowcolor[gray]{.9}
100	&	\texttt{MitsosBarton2006Ex322}	&	0.21 	&	0.01 	&	3.13E-10	&	&	0.01 	&	1.54E-09	&	&	0.21 	&	4.98E-07	&	&	0.21 	&	3.62E-07	&	&	0.21 	&	1.20E-07	\\
101	&	\texttt{MitsosBarton2006Ex323}	&	0.18 	&	0.18 	&	1.69E-07	&	&	0.25 	&	5.53E-07	&	&	0.05 	&	2.59E-01	&	&	0.18 	&	2.55E-12	&	&	0.18 	&	1.64E-07	\\\rowcolor[gray]{.9}
102	&	\texttt{MitsosBarton2006Ex324}	&	-1.76 	&	-1.76 	&	8.98E-09	&	&	-1.75 	&	3.21E-11	&	&	-1.75 	&	3.11E-10	&	&	-1.75 	&	7.47E-07	&	&	-1.75 	&	2.01E-08	\\
103	&	\texttt{MitsosBarton2006Ex325}	&	-1.00 	&	0.00 	&	2.34E-04	&	&	0.00 	&	4.35E-04	&	&	-0.18 	&	2.25E-08	&	&	0.00 	&	1.00E-01	&	&	0.00 	&	1.71E-05	\\\rowcolor[gray]{.9}
104	&	\texttt{MitsosBarton2006Ex326}	&	-2.35 	&	-2.00 	&	1.91E-08	&	&	-2.04 	&	4.76E-08	&	&	-1.67 	&	1.76E-01	&	&	-2.35 	&	5.55E-08	&	&	-2.00 	&	8.01E-07	\\
105	&	\texttt{MitsosBarton2006Ex327}	&	2.00 	&	0.00 	&	9.77E-07	&	&	0.00 	&	1.93E-07	&	&	0.00 	&	7.42E-07	&	&	0.02 	&	1.44E-07	&	&	0.20 	&	4.02E-07	\\\rowcolor[gray]{.9}
106	&	\texttt{MitsosBarton2006Ex328}	&	-10.00 	&	-10.00 	&	6.48E-09	&	&	-5.94 	&	1.70E-04	&	&	-2.74 	&	5.59E-06	&	&	-4.24 	&	6.91E-01	&	&	-2.61 	&	1.54E-03	\\
107	&	\texttt{NieWangYe2017Ex34}	&	2.00 	&	2.00 	&	1.60E-03	&	&	3.13 	&	1.18E+00	&	&	2.00 	&	2.01E-03	&	&	2.00 	&	2.04E-03	&	&	2.49 	&	1.58E+01	\\\rowcolor[gray]{.9}
108	&	\texttt{NieWangYe2017Ex52}	&	-1.71 	&	-1.41 	&	4.85E-08	&	&	-0.07 	&	6.59E-07	&	&	-2.14 	&	5.94E-07	&	&	0.35 	&	2.13E-10	&	&	-0.44 	&	1.32E-08	\\
109	&	\texttt{NieWangYe2017Ex57}	&	-2.00 	&	0.00 	&	1.80E-05	&	&	0.00 	&	2.37E-05	&	&	-2.04 	&	8.89E-12	&	&	0.00 	&	8.17E-01	&	&	0.00 	&	1.28E-03	\\\rowcolor[gray]{.9}
110	&	\texttt{NieWangYe2017Ex58}	&	-3.49 	&	-3.54 	&	2.74E-09	&	&	-3.53 	&	4.78E-08	&	&	0.00 	&	1.46E-08	&	&	0.00 	&	8.46E-09	&	&	0.00 	&	3.25E-07	\\
111	&	\texttt{NieWangYe2017Ex61}	&	-1.02 	&	-4.16 	&	1.99E-10	&	&	0.00 	&	2.92E-07	&	&	-1.03 	&	3.44E-06	&	&	-0.11 	&	1.48E-05	&	&	0.00 	&	3.36E-07	\\\rowcolor[gray]{.9}
112	&	\texttt{Outrata1990Ex2d}	&	0.92 	&	0.35 	&	4.72E-08	&	&	0.85 	&	1.47E-10	&	&	0.85 	&	1.40E-10	&	&	0.85 	&	2.66E-12	&	&	0.85 	&	2.22E-07	\\
113	&	\texttt{Outrata1990Ex2e}	&	0.90 	&	0.66 	&	5.82E-10	&	&	0.90 	&	1.70E-07	&	&	0.90 	&	2.02E-09	&	&	6.26 	&	4.16E-01	&	&	6.96 	&	4.16E-01	\\\rowcolor[gray]{.9}
114	&	\texttt{Outrata1993Ex31}	&	1.56 	&	0.88 	&	4.57E-09	&	&	1.56 	&	1.10E-11	&	&	1.56 	&	4.13E-10	&	&	1.56 	&	8.91E-10	&	&	12.98 	&	1.59E+00	\\
115	&	\texttt{PaulaviciusAdjiman2017a}	&	0.25 	&	0.00 	&	7.69E-16	&	&	0.00 	&	7.69E-16	&	&	0.00 	&	7.69E-16	&	&	0.00 	&	8.74E-16	&	&	0.00 	&	7.69E-16	\\\rowcolor[gray]{.9}
116	&	\texttt{SinhaMaloDeb2014TP3}	&	-18.68 	&	-26.98 	&	5.09E-09	&	&	-24.64 	&	3.08E-07	&	&	-18.79 	&	7.25E-07	&	&	-18.68 	&	4.17E-07	&	&	-5.27 	&	2.73E-08	\\
117	&	\texttt{SinhaMaloDeb2014TP7}	&	-1.96 	&	-1.98 	&	5.13E-07	&	&	-1.98 	&	2.34E-10	&	&	-1.98 	&	1.45E-09	&	&	-1.98 	&	5.20E-11	&	&	-1.98 	&	2.44E-07	\\\rowcolor[gray]{.9}
118	&	\texttt{SinhaMaloDeb2014TP9}	&	0.00 	&	0.00 	&	7.38E-17	&	&	0.00 	&	2.21E-16	&	&	0.00 	&	6.64E-16	&	&	0.00 	&	1.99E-15	&	&	0.00 	&	5.97E-15	\\
119	&	\texttt{SinhaMaloDeb2014TP10}	&	0.00 	&	0.00 	&	7.38E-18	&	&	0.00 	&	2.21E-17	&	&	0.00 	&	6.64E-17	&	&	0.00 	&	1.99E-16	&	&	0.00 	&	5.97E-16	\\\rowcolor[gray]{.9}
120	&	\texttt{Vogel2012}	&	0.00 	&	0.00 	&	4.22E-07	&	&	0.00 	&	8.32E-14	&	&	0.00 	&	9.60E-07	&	&	0.00 	&	6.29E-08	&	&	0.00 	&	7.76E-11	\\
121	&	\texttt{YeZhu2010Ex42}	&	1.00 	&	1.00 	&	2.98E-07	&	&	1.00 	&	2.37E-07	&	&	1.00 	&	2.46E-07	&	&	1.00 	&	2.96E-07	&	&	1.00 	&	2.26E-07	\\\rowcolor[gray]{.9}
122	&	\texttt{YeZhu2010Ex43}	&	1.25 	&	1.00 	&	1.63E-09	&	&	1.00 	&	2.89E-08	&	&	1.00 	&	1.95E-07	&	&	1.00 	&	6.59E-09	&	&	1.00 	&	2.12E-10	\\
123	&	\texttt{Zlobec2001a}	&	-1.00 	&	0.00 	&	4.00E-02	&	&	-0.11 	&	1.12E-01	&	&	-1.00 	&	9.59E-09	&	&	0.00 	&	1.34E-03	&	&	-1.00 	&	5.94E-09	\\\rowcolor[gray]{.9}
124	&	\texttt{Zlobec2001b}	&	NaN	&	0.08 	&	3.19E-03	&	&	0.00 	&	6.91E-08	&	&	1.00 	&	4.01E-12	&	&	1.00 	&	2.46E-08	&	&	0.13 	&	8.06E-04	\\

		  \end{tabular}}
\end{sidewaystable}

\begin{figure}[th]
\begin{center}
  \includegraphics[width=.49\linewidth]{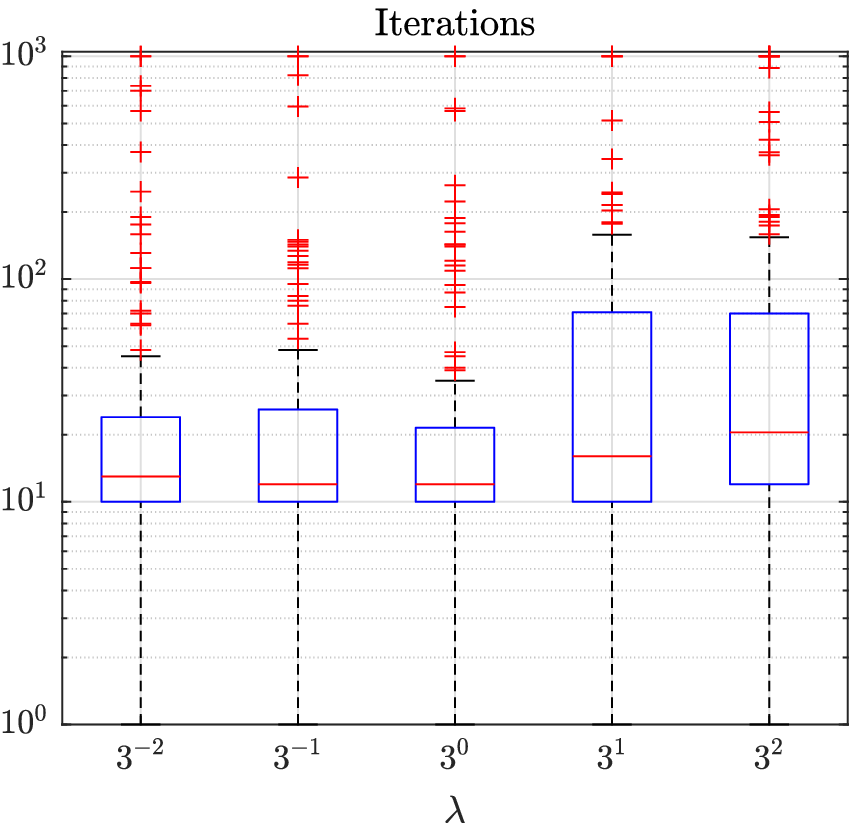}
  \includegraphics[width=.49\linewidth]{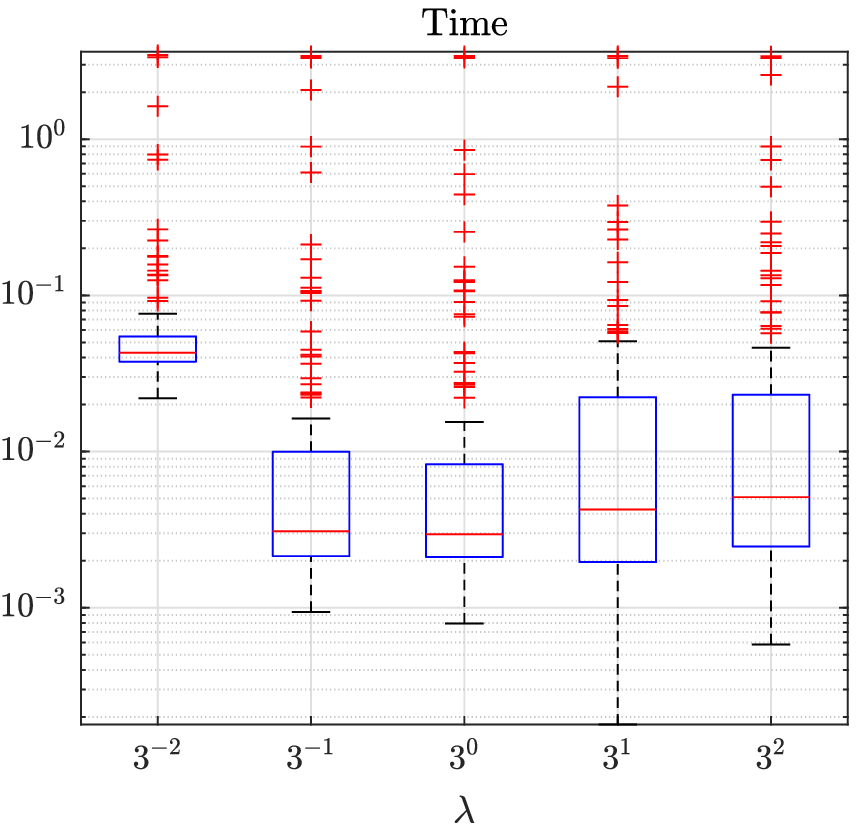}
  \captionof{figure}{Box plot of iterations and time.}
  \label{iter}
\end{center}
\end{figure}
 
  {As the true solutions for many of the examples in our test set are known (see supplementary material in \cite{ZhouZemkohoDetailed2018}),  we can count that our method closely achieved $\bar{F}$ for all problems except for $21$ ones. They are problems number 2,    17,    42,    46,    54,    63,    70,    78,    82,    84,    85,    92,    93,    96,    98,   103,   105,   108,   112,   115, and 122 (cf. numbering in the first column of Table \ref{numerical-results}). 
 However, it is worth mentioning that not closely achieving $\bar{F}$ does not mean that our method did not find a better solution because some known best solutions may not be the global optimal ones. For instance, the known best solutions of problems with numbers 17 and 54 might not be optimal, which indicates that the solutions produced by our methods such that  $\hat{F}<\bar{F}$ could be better than the known ones for these two problems.   
 Also observing that for 32 examples out of the 124 examples (i.e., problems number 5,     7,    11,    13,    30-35,    40,    41,    48,    55,    56,    80,    83-85,    88, 92,    93,    96,    98,   102,   115, and   117-122), we obtained the same values for $\hat{F}$ for all choices of $\lambda$. }

\newpage
 For each value of $\lambda\in \tilde{\Lambda}$, we draw the box plot of the number of iterations in the left sub-figure of Figure \ref{iter}.  In each box, the central mark (red
line) indicates the median, the bottom, and top edges of the
box indicate the 25th and 75th percentiles, respectively. The
outliers are plotted individually using ‘+’ symbols.   {In general, the larger $\lambda$ is, the more the average number of iterations is.} Moreover, from Table \ref{wrap-tab:1}, for each $\lambda\in \tilde{\Lambda}$,   {over $112/124\approx 90.32\%$ problems were solved by using less than 200 iterations. Similarly, we also plot the CPU time (in seconds) in the right sub-figure of Figure \ref{iter}. In terms of the average values,   it seems that $\lambda=3^{0}$ allowed the algorithm to run the fastest.  A majority of problems were solved within 1 second and as presented in  Table \ref{wrap-tab:1}, $119$ (resp. $119, 120, 121$, and $119$) problems were solved within 0.5 second for $\lambda=3^{-2}$ (resp. $3^{-1}, 3^0, 3^1$, and $3^2$), which indicates that our algorithm can address $119/124\approx 95.96\%$ problems within 0.5 second for each $\lambda\in\tilde{\Lambda}$.}

 Considering the local version of Algorithm \ref{algorithm-SSNBO} (i.e., the corresponding version without the step size), it might be useful to evaluate how often it terminates with $\alpha^k=1$. As shown in Table \ref{wrap-tab:1}, it turns out that for $\lambda = 3^{-2}$ (resp.  $\lambda=3^{-1}$, $\lambda=3^0$,  $\lambda=3^{1}$, and $\lambda=3^{2}$),   {our algorithm solved $106$ (resp. $108$, $101	$, $100	$, and $99$)} problems and stopped with $\alpha^k=1$. This confirms that Algorithm \ref{algorithm-SSNBO} also   {has very good} local behavior in practice. 

 \begin{table}[H]
\centering\caption{Number of problems solved by Algorithm \ref{algorithm-SSNBO} under each $\lambda\in\tilde{\Lambda}$.}\label{wrap-tab:1}\vspace{-5mm}
\renewcommand\baselinestretch{1.15}\selectfont
\begin{tabular}{lccccc}\\\hline &$\lambda=3^{-2}$&$\lambda=3^{-1}$&$\lambda=3^0$&$\lambda=3^{1}$&$\lambda=3^{2}$\\\hline
Iteration $< 100$& 111	&	110	&	109	&	98	&	97	\\
Iteration $< 200$& 117	&	119	&	118	&	115	&	112	\\
Time $< 0.1$ &   110	&	112	&	112	&	115	&	109	\\
Time $< 0.5$ &   119	&	119	&	120	&	121	&	119	\\
$\alpha^k=1$ &     106	&	108	&	101	&	100	&	99	\\
$y^k\approx \varsigma^k$ &38	&	38	&	53	&	65	&	84	\\\hline
\end{tabular}
\end{table}
 
 Finally, we assess how often we get $y=\varsigma$. As mentioned in Section \ref{First order necessary optimality conditions},  $y=\varsigma$ generates a stationary point in the sense of Theorem \ref{KN stationarity partial calm 0}. As for $v$ and $w$, they are not necessarily related, as $w$ is a multiplier associated with the lower-level constraint function but only in connection to the lower-level problem; $v$ is associated with the same function, $g$, but in connection to the upper-level problem. But considering the similarity in the complementarity systems \eqref{VS-r4} and \eqref{KS-r3}, we now look a bit closely to see how often they could coincide. To proceed, we recorded how many examples in Table \ref{wrap-tab:1}, for each $\lambda$, were solved with $y^k\approx \varsigma^k$ (at the final iteration), where symbol ``$\approx$'' used here refers to a very small relative error, i.e., $\|y^k-\varsigma^k\|/\max\{1,\|y^k\|+\|\varsigma^k\|\}<0.05$.   {For $\lambda=3^{-2},3^{-1},3^0,3^1$, and $3^2$ we respectively had $38, 38, 53, 65$, and $84$ examples solved with  $y^k\approx \varsigma^k$. }\\

\noindent{ {{\bf d) Solving a discretized bilevel optimal control (BOC) program 
 with larger sizes.} We point out that all examples in the BOLIB library are on a small scale, with dimensions satisfying $\max\{n, m, p, q\} \leq20$. Therefore, to see the performance of Algorithm \ref{algorithm-SSNBO} on solving problems on larger scales, we take advantage of a BOC program from \cite{MehlitzGerd16}. This is a quadratic program and its dimensions $\{n,m,p,q\}$ can be altered. The optimization model is decided by the following functions
$$
\begin{array}{rll}
  F(x,y)& := &\frac{1}{2}[(
  							y^1;0)-c]^\top
  						  D[(
  							y^1;
  							  0)-c] -d^\top x,\\
  G(x,y)& := &\left( - x_1+x_2-1;-x \right),\\	
  f_1(x,y)& := &\frac{1}{2}(Cy^1-Px)^\top U(Cy^1-Px)+\frac{\sigma}{2}(y^2-Qx)^\top V(y^2-Qx),\\
  g_1(x,y)& := & \left(
                       y^2-u;
                      -y^2+l;
                     Ay;
                       - Ay
     \right),
\end{array}
$$
where $x\in\mathbb{R}^2$,  $y=( y^1;y^2)$ with $y^i\in\mathbb{R}^{m_i}$, $D\in\mathbb{R}^{m_1\times m_1}$, $d\in\mathbb{R}^{n}$, $c\in\mathbb{R}^{m}$, $C\in\mathbb{R}^{s\times m_1}$, $P\in\mathbb{R}^{s\times n}$, $U\in\mathbb{R}^{s\times s}$, $Q\in\mathbb{R}^{m_2\times n}$, $V\in\mathbb{R}^{m_2\times m_2}$, $u\in\mathbb{R}^{m_2}$, $l\in\mathbb{R}^{m_2}$, $A\in\mathbb{R}^{t\times m}$ are given data. Here, $(a;b)=(a^\top~b^\top)^\top.$ 
Similarly to  \cite{MehlitzGerd16}, we set $m_1=m_2=s=t$ and thus the dimension is  $(n, m, p, q)=(2, 2t, 3, 4t)$. Results are presented in Table \ref{wrap-tab:2}. As the sizes of $t$ increase, it is challenging for the algorithm to meet the stopping criterion, i.e., $\|\hat{\Phi}\|<10^{-6}.$ This is reasonable since such a stopping criterion is relaxed for small-scale $\hat{\Phi}$ but relatively crucial for large-scale $\hat{\Phi}$. It seems that the proposed algorithm achieved similar $\hat{F}$ for $\lambda\in\{3^{-1},3^{0},3^{1},3^{2}\}$.

 \begin{table}[htp]
\centering
\caption{Results of Algorithm \ref{algorithm-SSNBO} solving BOC under each $\lambda\in\tilde{\Lambda}$.}\label{wrap-tab:2}\vspace{-5mm}
\renewcommand\baselinestretch{1.15}\selectfont
\begin{tabular}{llrrrrr}\\\hline 
$(n,m,p,q)$	&		&	$\lambda=3^{-2}$	&	$\lambda=3^{-1}$	&	$\lambda=3^{0}$	&	$\lambda=3^{1}$	&	$\lambda=3^{2}$	\\\hline
\multirow{4}{*}{$(2,100,3,200)$}	&	$\hat{F}$	&	0.428	&	0.515	&	0.515	&	0.515	&	0.515	\\
	&	$\|\hat{\Phi}\|$	&	9.47E-07	&	6.90E-07	&	6.85E-07	&	4.71E-07	&	8.51E-07	\\
	&	Iterations	&	86 	&	16 	&	16 	&	43 	&	52 	\\
	&	Time	&	0.468	&	0.059	&	0.052	&	0.125	&	0.199	\\\hline
\multirow{4}{*}{$(2,500,3,1000)$}	&	$\hat{F}$	&	0.537	&	0.597	&	0.599	&	0.599	&	0.599	\\
	&	$\|\hat{\Phi}\|$	&	1.56E-03	&	1.61E-04	&	1.62E-04	&	5.56E-07	&	7.59E-07	\\
	&	Iterations	&	123 	&	121 	&	119 	&	45 	&	71 	\\
	&	Time	&	11.602	&	10.833	&	11.114	&	3.376	&	5.536	\\\hline
\multirow{4}{*}{$(2,1000,3,2000)$}	&	$\hat{F}$	&	0.536	&	0.585	&	0.587	&	0.588	&	0.587	\\
	&	$\|\hat{\Phi}\|$	&	2.31E-03	&	4.24E-04	&	6.35E-04	&	6.63E-04	&	4.62E-04	\\
	&	Iterations	&	139 	&	127 	&	126 	&	160 	&	240 	\\
	&	Time	&	43.736	&	40.578	&	39.567	&	50.667	&	73.022	\\\hline
\end{tabular}
\end{table}
}}

  {To summarize, Algorithm \ref{algorithm-SSNBO} enables solving a wide range of problems (regardless of the convexity) with desirable accuracy and computational speed. However, as illustrated earlier, the efficacy of the algorithm relies on the selection of parameters $\lambda$ and initial points. Moreover, in scenarios where testing problems exhibit poor singularity conditions, the local convergence behavior of  Algorithm \ref{algorithm-SSNBO} may suffer significant degradation.}

\section{Conclusion and final comments}
We have developeded a new upper estimate for the coderivative of the solution set-valued mapping of a parametric QVI, and subsequently, we used it to derive sufficient conditions ensuring that this mapping is Lipschitz-like in the sense of Aubin (\cite{Aubin1984}); cf. Section \ref{Coderivative and robust stability of solution maps}. Unlike in the existing literature (see, e.g., \cite{MOC07}), we do not require the convexity assumption on the feasible set of the QVI. This was possible thanks to the value function reformulation that we borrowed from the bilevel optimization literature (see, e.g., \cite{DempeZemkohoBook} and references therein), but which seem not to have been exploited in the QVI literature. This value function reformulation comes with the additional advantage that it allows our analysis to be conducted without requiring second order information as it is often the case when the GE reformulation is used for a similar analysis. However, a drawback of the value function reformulation is that it does not enable tractable dual qualification conditions of the basic--type constraint qualification in \eqref{BCQ}, for example. Hence, we use the calmess condition (for set-valued mappings) to avoid the failure of basic-type constraint qualifications such as the Mangasarian-Fromovitz-type, for example. 

As second major contribution of the paper, we consider an optimization with quasi-variational inequality constraint (OPQVI), and construct necessary optimality conditions, which then served as base for a semismooth Newton-type method used to solve the problem; cf. Section \ref{Optimization problems with a QVI constraint}. As we also use the value function reformulation introduced for the stability analysis of the QVI, we exploit all the relevant similarity to facilitate the study the OPQVI without much duplication from Section \ref{Coderivative and robust stability of solution maps}. One of the challenges of our algorithm is the need for the partial penalization parameter ($\lambda$), which as for exact penality methods, is difficult to select. However, the experiments on over 125 examples are only based on five different penalization parameters, and seem to show a remarkable performance of our method, which could inspire further work for similar problem classes.

\section*{Data availability statement} The BOLIB test set used for the numerical experiments can be found at \href{https://biopt.github.io/bolib/}{https://biopt.github.io/bolib/}



\begin{thebibliography}{plain}


\bibitem{AdamHenrionOutrata2018} Adam L,  Henrion R, Outrata J (2018) On M-stationarity conditions in MPECs and the associated qualification conditions. Math Program 168:229-259

\bibitem{Andreani-2012} Andreani R, Haeser, G,  Schuverdt, ML, Silva, PJS (2012) A relaxed constant positive linear dependence constraint qualification and applications. Math Program 135(1):255-273

\bibitem{Aubin1984} Aubin JP (1984) Lipschitz behavior of solutions to convex minimization problem. Math Oper Res 9:87-111

\bibitem{bank-1983} Bank B, Guddat J, Klatte D, Kummer B, Tammer K (1983) Nonlinear parametric optimization. Birkh\"{a}user

\bibitem{MehlitzMinchenko2020-100} Bednarczuk EM, Minchenko LI, Rutkowski KE (2019) On Lipschitz-like continuity of a class of set-valued mappings. Optimization 69(12):2535-2549

\bibitem{BonnansShapiro2000}
Bonnans JF, Shapiro A (2000) Perturbation Analysis of Optimization Problems. Springer

\bibitem{ChanPang1982} Chan D, Pang JS (1982) The generalized quasi-variational inequality problem. Math Oper Res 7(2):211-222 

\bibitem{ClarkeBook1983} Clarke FH (1994) Optimization and nonsmooth analysis. SIAM Classics in Applied Mathematics 5

\bibitem{DeLuca1996} De Luca T, Facchinei F,  Kanzow C (1996) A semismooth equation approach to the solution of nonlinear
  complementarity problems. Math Program 75(3):407-439
  
\bibitem{DempeFoundations2002} Dempe S (2002) Foundations of bilevel programming. Kluwer Academic Publishers

\bibitem{DempeDuttaMordukhovichNewNece} Dempe S, Dutta J, Mordukhovich, BS (2007) New necessary optimality conditions in optimistic bileve programming. Optimization 56(5-6):577-604

\bibitem{DempeZemkohoBook} Dempe S, Zemkoho AB (Eds)  (2020) Bilevel Optimization: Advances and Next Challenges. Springer

\bibitem{DempeZemkohoGenMFCQ} Dempe S,  Zemkoho, AB (2011) The generalized Mangasarian-Fromowitz constraint qualification
  and optimality conditions for bilevel programs. J Optim Theory Appl 148:46-68
  
\bibitem{Dreves2016} Dreves A (2016) Uniqueness for quasi-variational inequalities. Set-Valued Var Anal 24:285-297

\bibitem{Fiacco1983}
Fiacco AV (1983) Introduction to Sensitivity and Stability Analysis in Nonlinear Programming. Academic Press Inc.

\bibitem{FischerASpecial1992}  Fischer A (1992) A special newton-type optimization method. Optimization 24(3-4):269-284

\bibitem{FischerZemkohoZhouSemismooth2019} Fischer A, Zemkoho AB, Zhou S (2022) Semismooth Newton-type method for bilevel optimization: Global convergence and extensive numerical experiment. Optim Methods Softw 37(5):1770-1804

\bibitem{FK92} Flam SD, Kummer B (1992) Great fish wars and Nash equilibria, Technical Report WP-0892, Department of Economics, University of Bergen, Norway

\bibitem{FliegeTinZemkoho2021} Fliege J, Tin A, Zemkoho AB (2021) Gauss-Newton-type methods for bilevel optimization. Comput Optim Appl 78(3):793-824

\bibitem{Galantai2012} Gal{\'a}ntai A (2012) Properties and construction of NCP functions. Comput Optim Appl 52(3):805-824

\bibitem{GfrererYe2020NewConst} Gfrerer H, Ye JJ (2017) New constraint qualifications for mathematical programs with equilibrium constraints via variational analysis. SIAM J Optim 27(2):842-865

\bibitem{GfrererYe2020NewNecc} Gfrerer H,  Ye JJ (2020) New sharp necessary optimality conditions for mathematical programs with equilibrium constraints.  Set-Valued Var Anal  28(2):395-426

\bibitem{H91} Harker PT (1991) Generalized Nash games and quasi-variational inequalities. European J Oper Res 54(1):81-94 

\bibitem{HenrionJouraniOutrataCalmness2002} Henrion R, Jourani A, Outrata JV (2002) On the calmness of a class of multifunctions. SIAM J Optim 13:603-618 

\bibitem{hen-out01} Henrion R, Outrata JV (2001) A subdifferential condition for calmness of multifunctions. J Math Anal Appl 258(1):110-130 

\bibitem{HenrionOutrataSurowiec2012} Henrion R, Outrata JV, Surowiec T (2012) Analysis of M-stationary points to an EPEC modeling oligopolistic competition in an electricity spot market. ESAIM Control Optim Calc Var 18(2):295-317

\bibitem{HenrionSurowiecCalmness2011} Henrion R,  Surowiec T (2011) On calmness conditions in convex bilevel programming. Appl  Anal 90(6):951-970

\bibitem{HobbsPang2007} Hobbs BF, Pang JS (2007) Nash-Cournot Equilibria in Electric Power Markets with Piecewise Linear Demand Functions and Joint Constraints. Oper Res 55(1):113-127

\bibitem{JolaosoMehlitzZemkoho2024} 
Jolaoso LO, Mehlitz P, Zemkoho AB (2024) A fresh look at nonsmooth Levenberg–Marquardt methods with applications to bilevel optimization. Optimization \href{https://doi.org/10.1080/02331934.2024.2313688}{https://doi.org/10.1080/02331934.2024.2313688}

\bibitem{QiJiangSemismooth1997} Jiang H,  Qi L (1997) Semismooth karush-kuhn-tucker equations and convergence analysis of newton and quasi-newton methods for solving these equations. Math Oper Res 22(2):301-325

\bibitem{MehlitzMinchenko2020} Mehlitz P, Minchenko  LI (2022) R-regularity of set-valued mappings under the relaxed constant positive linear dependence constraint qualification with applications to parametric and bilevel optimization. Set-Valued Var Anal 30:179-205

\bibitem{MehlitzGerd16}
 Mehlitz P, Wachsmuth G (2016) Weak and strong stationarity in generalized bilevel programming and bilevel optimal control.  Optimization 65:907-935.

\bibitem{MehlitzZemkoho2021}
 Mehlitz P, Zemkoho AB (2021) Sufficient optimality conditions in bilevel programming. Math Oper Res 

\bibitem{Mifflin1977} Mifflin R (1977) Semismooth and semiconvex functions in constrained optimization. SIAM J Control Optim 15(6):959-972

\bibitem{MordukhovichGeneralizedDifferential1994} Mordukhovich BS (1994) Generalized Differential Calculus for Nonsmooth and Set-Valued Mappings. J Math Anal Appl 183(1):250-288 

\bibitem{MordukhovichBook2006} Mordukhovich BS (2006) Variational analysis and generalized differentiation. Springer, Berlin

\bibitem{MordukhovichNamPhanVarAnalMargBlP} Mordukhovich  BS, Nam MN, Phan HM (2012) Variational analysis of marginal functions with applications to bilevel programming. J Optim Theory Appl 152:557-586

\bibitem{MOC07} Mordukhovich BS, Outrata  JV (2007) Coderivative analysis of quasi-variational inequalities
with applications to stability and optimization. SIAM J Optim 18:389-412

\bibitem{MSS82} Murphy FH, Sherali HD, Soyster AL (1982) A mathematical programming approach for determining oligopolistic market equilibrium. Math Program 24:92-106 

\bibitem{OutrataZowe1995} Outrata JV, Zowe  J (1995) A Newton method for a class of quasi-variational
inequalities. Comput Optim Appl 4(1):5-21

\bibitem{PangFukushima2005} Pang JS, Fukushima  M (2005) Quasi-variational inequalities, generalized Nash equilibria,
and multi-leader-follower games. Comput Manag Sci 2(1):21-56

\bibitem{Qi1993convergence} Qi L (1993) Convergence analysis of some algorithms for solving nonsmooth equations. Math Oper Res 18(1):227-244

\bibitem{QiSunANonsmoothVersion1993} Qi L, Sun J (1993) A nonsmooth version of Newton's method. Math Program 58:353-367

\bibitem{QiSun1999} Qi L, Sun D (1999) A Survey of Some Nonsmooth Equations and Smoothing Newton Methods. In: Eberhard A, Hill  R, Ralph D, Glover  BM (Eds) Progress in Optimization, pp. 121-146. Springer US, Boston, MA

\bibitem{Robinson1981} Robinson  SM (1981) Some continuity properties of polyhedral multifunctions.  Math Program Study 14:206-2014

\bibitem{RockafellarWetsBook1998} Rockafellar RT, Wets RJ-B (1998) Variational Analysis. Springer, Berlin


\bibitem{TinZemkohoLevenberg2022} Tin A, Zemkoho AB (2023) Levenberg-Marquardt method and partial exact penalty parameter selection in bilevel optimization. Optim Eng 24:1343-1385 

\bibitem{Ye1998New} Ye  JJ (1998) New uniform parametric error bounds. J Optim Theory Appl 98(1):197-219

\bibitem{YeZhuOptCondForBilevel1995} Ye JJ, Zhu DL (1995) Optimality conditions for bilevel programming problems. Optimization 33:9-27

\bibitem{zemkoho2017estimates} Zemkoho  AB (2022) Estimates of generalized Hssians for optimal value functions in
  mathematical programming. Set-Valued Var Anal 30:847-871
  
\bibitem{BOLIB2017} Zhou  S, Zemkoho  AB, Tin A (2020) BOLIB: Bilevel Optimization Library of Test Problems. In: Dempe S, Zemkoho AB, Bilevel Optimization, pp. 563-580. Springer, Cham

\bibitem{ZhouZemkohoDetailed2018}Dutta J, Lahoussine L, Zemkoho AB, Zhou S (2023) Detailed numerical experiment for {\em Nonconvex quasi-variational inequalities: stability analysis and application to numerical optimization} \href{https://github.com/abzemkoho/QVI}{https://github.com/abzemkoho/QVI}

\bibitem{ZemkohoZhou2021} Zemkoho AB, Zhou S (2021) Theoretical and numerical comparison of the Karush-Kuhn-Tucker and value function reformulations in bilevel optimization. Comput Optim Appl 78(2):625-674

\bibitem{Huang} Huang HY (2009) Unified approach to quadratically convergent algorithms for function minimization. J Optim Theory Appl 5:354-405

\bibitem{Miele} Miele A (Ed) (1965) Theory of Optimum Aerodynamic Shapes. Academic Press, New York (1965)


\bibitem{Nocedal}  Nocedal J, Wright SJ (2000) Numerical Optimization. Springer, New York 

\bibitem{Ralston} Ralston A (1960) Numerical integration methods for the solution of ordinary differential equations. In: Ralston A, Wilf HS (Eds): Mathematical Methods for Digital Computers 1, pp. 95-109. Wiley, New York

\bibitem{Yang} Yang TL (1969) Optimal control for a rocket in a three-dimensional central force field. Technical Memorandum TM 69-2011-2, Bellcomm


\end{thebibliography}
\end{document}